\documentclass[11pt]{amsart}

\usepackage{verbatim}
\usepackage{amsmath,amsthm,amssymb,amscd}
\usepackage{mathrsfs}
\usepackage[all]{xy}
\usepackage{xspace}
\usepackage{ulem}
\usepackage{enumerate}
\usepackage{todonotes}
\usepackage{xcolor}
\usepackage{hyperref}

\newcommand{\R}{\ensuremath{\mathbb{R}}\xspace}
\newcommand{\C}{\ensuremath{\mathbb{C}}\xspace}

\newcommand{\nt}{\ensuremath{\mathbb{N}}\xspace}
\newcommand{\inc}{\ensuremath{\hookrightarrow}}
\newcommand{\id}{\vartriangleleft}

\newcommand{\pr}[2]{\ensuremath{\langle {#1},{#2}\rangle}}
\newcommand{\norm}[1]{\ensuremath{\|#1\|}}

\newcommand{\ma}{\ensuremath{\mathcal{A}}\xspace}
\newcommand{\mb}{\ensuremath{\mathcal{B}}\xspace}
\newcommand{\mc}{\ensuremath{\mathcal{C}}\xspace}

\newcommand{\mh}{\ensuremath{\mathcal{H}}\xspace}
\newcommand{\mi}{\ensuremath{\mathcal{I}}\xspace}

\newcommand{\mk}{\ensuremath{\mathcal{K}}\xspace}
\newcommand{\ml}{\ensuremath{\mathcal{L}}\xspace}

\newcommand{\mn}{\ensuremath{\mathcal{N}}\xspace}

\newcommand{\ms}{\ensuremath{\mathcal{S}}\xspace}

\renewcommand{\r}[1]{\ensuremath{\big|_{#1}}\xspace}

\newcommand{\bts}[3]{\ensuremath{{#1}\bigotimes_{#2}{#3}}\xspace}
\newcommand{\bats}[2]{\ensuremath{{#1}\bigodot{#2}}\xspace}
\newcommand{\tmin}[2]{\ensuremath{{#1}\bigotimes_{\min}{#2}}\xspace}
\newcommand{\tmax}[2]{\ensuremath{{#1}\bigotimes_{\max}{#2}}\xspace}
\newcommand{\chal}{\overline{\alpha}}

\newcommand{\al}{\ensuremath{\alpha}\xspace} 
\newcommand{\del}{\ensuremath{\delta}\xspace} 

\newcommand{\la}{\ensuremath{\lambda}\xspace}

\newcommand{\La}{\ensuremath{\Lambda}\xspace}

\newcommand{\Del}{\ensuremath{\Delta}\xspace}

\newcommand{\ov}[1]{\overline{#1}}

\newcommand{\cs}{$C^*$-algebra\xspace}
\newcommand{\css}{$C^*$-algebras\xspace}

\newcommand{\scss}{$C^*$-subalgebras\xspace}
\newcommand{\ct}{$C^*$-tring\xspace}
\newcommand{\cts}{$C^*$-trings\xspace}
\newcommand{\hm}{homomorphism\xspace}
\newcommand{\hms}{homomorphisms\xspace}

\newcommand{\rep}{representation\xspace}
\newcommand{\reps}{representations\xspace}

\newcommand{\rrs}{regular \reps}

\newcommand{\ap}{approximation property\xspace}
\newcommand{\aps}{approximation properties\xspace}

\newcommand{\ilt}{inductive limit topology\xspace}
\newcommand{\fb}{Fell bundle\xspace}
\newcommand{\fbs}{Fell bundles\xspace}
\newcommand{\es}{exact sequence\xspace}

\newcommand{\fela}{\ensuremath{\ma=(A_t)_{t\in G}}\xspace}

\newcommand{\felbh}{\ensuremath{\mb=(B_s)_{s\in H}}\xspace}

\newcommand{\adj}[1]{\mathcal{L}({#1})}

\newcommand{\resp}{\hspace*{1cm}}

\newcommand{\bc}{\begin{center}}
\newcommand{\ec}{\end{center}}
\newcommand{\be}{\begin{enumerate}}
\newcommand{\ee}{\end{enumerate}}
\newcommand{\bi}{\begin{itemize}}
\newcommand{\ei}{\end{itemize}}
\newcommand{\bd}{\begin{description}}
\newcommand{\ed}{\end{description}}
\newcommand{\beq}{\begin{equation}}
\newcommand{\eeq}{\end{equation}}
\newcommand{\beqa}{\begin{eqnarray}}
\newcommand{\eeqa}{\end{eqnarray}}
\newcommand{\bfr}{\begin{flushright}}
\newcommand{\efr}{\end{flushright}}
\newcommand{\bfl}{\begin{flushleft}}
\newcommand{\efl}{\end{flushleft}}
\newcommand{\bp}{\begin{picture}}
\newcommand{\ep}{\end{picture}}

\DeclareMathOperator{\gen}{span}

\DeclareMathOperator{\supp}{supp}

\newcommand{\cgen}{\ensuremath{\ov{\gen}}}

\title{Tensor Products of Fell Bundles over groups}

\address{Centro de Matem\'atica\\ 
         Facultad de Ciencias\\
         Universidad de la Rep\'ublica\\
         Igu\'a 4225\\
         CP 11400\\
         Montevideo--URUGUAY\\
       }
       \email{fabadie@cmat.edu.uy}

\theoremstyle{plain}
\newtheorem{thm}{Theorem}[section]
\newtheorem{prop}[thm]{Proposition}
\newtheorem{lem}[thm]{Lemma}
\newtheorem{cor}[thm]{Corollary}          

\theoremstyle{definition}
\newtheorem{df}[thm]{Definition}

\theoremstyle{remark}
\newtheorem{rk}[thm]{Remark}

\date{}

\begin{document}

\author[F. Abadie]{Fernando Abadie}
\begin{abstract}
We extend the theory of tensor products of C*-algebras to the larger
category of Fell bundles over locally compact groups. We prove that, like in the case of C*-algebras, there exist maximal and minimal tensor products. Given two Fell bundles, we compare the tensor products of their cross-sectional algebras with the cross-sectional algebras of their tensor products. As
applications we prove that, under certain conditions, the
cross-sectional C*-algebra of a Fell bundle is nuclear or exact whenever 
so is its fiber over the unit element of the group. 
\end{abstract}

\maketitle
\setcounter{tocdepth}{1}
\setcounter{secnumdepth}{3}
\tableofcontents
\section{Introduction}\label{sec:intro}
\par The original motivation for the present work was to study
nuclearity and exactness of crossed products by partial
actions, both important properties of \css related with tensor products.  
\par The best way to define and study crossed products by partial
actions is through the theory of $C^*$-algebraic bundles, today also called 
Fell bundles (for a comprehensive treatment of such theory see  
\cite{fd}). According to \cite{extwist}, given a partial action $\al$ 
of the locally compact group $G$ on the \cs $A$, or even a twisted
partial action, a Fell bundle $\mb_{\al}$ over $G$ is associated to~$\al$. The cross-sectional algebra of $\mb_{\al}$ is called the
crossed-product of $A$ by the partial action $\al$, and it is denoted
by $A\rtimes_{\al}G$. Similarly, the reduced cross-sectional algebra
of $\mb_{\al}$ is called the reduced crossed-product of $A$ by the
partial action $\al$, and it is denoted by $A\rtimes_{\al,r}G$ (in
Section \ref{sec:cross} we recall the definition of the reduced
cross-sectional algebra of a Fell bundle; for additional information
the reader is referred to \cite{exeng} and~\cite{env}). On the other hand, Fell bundles are closely related to partial actions, since not only many of them can be described as associated to twisted partial actions (\cite{extwist}), but any Fell bundle carries a natural partial action of its underlying group on the spectrum of its unit fiber (\cite{bafa}) and, in a sense, it is equivalent to the Fell bundle associated to a partial action (see \cite{fellequiv}, \cite{camila} and \cite{abf}).
\par A point 
exploited in this paper is that some properties of the cross-sectional
algebras of a \fb are in part just consequences of properties of the fibers of the bundle
itself, which in turn are many times directly related to those of the unit fiber. Moreover, some constructions with these algebras are better 
understood when they are made directly on the bundle. In particular
this viewpoint applies to tensor products. Thus we were led to define
and study tensor products of \fbs.
So posed in terms of Fell bundles, what we are interested in studying are the tensor products of cross-sectional algebras of Fell bundles, and the strategy we follow is to permute the order in which we consider such constructions, i.e., first define the tensor products of Fell bundles and then consider the cross-sectional algebras of the resulting bundles. In fact, what we will show is that these constructions “commute”, in the sense that, starting from two Fell bundles, the result is independent of the order in which we take the tensor product and the cross-sectional algebras. Furthermore, we will see that there is a perfect harmony in relation to the type of construction we choose in each case, i.e., maximal tensor products and full cross-sectional algebras, or spatial tensor products and reduced cross-sectional algebras (see below).

\par Let us describe briefly the contents and structure of the
paper. 
\par Since the fibers of a \fb are C*-ternary rings (\cts for short), the study carried out in \cite{trings} (in particular Section~5.2) can be considered as a preliminary step in the direction of studying tensor products of Fell bundles.  In the present work we will make considerable use of the results of \cite{trings}, so, for the reader's convenience, in the next section we will recall and expand on some of the aspects that interest us most. Also briefly discussed in this section will be the possibility of extending a C*-norm on the unit fiber of a *-algebraic bundle to a C*-norm on the entire bundle, which will lead to consideration of the notion of positive *-algebraic bundle.  
\par In the third section we deal with tensor products of
\fbs. If $\fela$ and $\felbh$ are \fbs over the locally compact groups 
$G$ and $H$ respectively, then a tensor product $\bts{\ma}{\al}{\mb}$
will be a \fb over $G\times H$, with fibers $\bts{A_t}{\al}{B_s}$. As
in the case of \css and \cts there are a maximal and a minimal tensor
products, which correspond respectively to the maximal and minimal
tensor products of the corresponding unit fibers of the bundles.
First we consider bundles over discrete groups, and show that the algebraic tensor product $\bats{\ma}{\mb}$ is a positive *-algebraic bundle. Then we define the tensor product $\bts{\ma}{\alpha}{\mb}$ for any C*-norm on $\bats{\ma}{\mb}$. Finally, we topologize $\bts{\ma}{\alpha}{\mb}$ for the case the base groups of $\ma$ and $\mb$ are general locally compact groups. We
end the section by generalizing some results on \reps of tensor
products of \css to the case of \fbs. 
\par The fourth section is devoted to comparing the cross-sectional
algebras of tensor products. Let $C^*(\mb)$ and $C^*_r(\mb)$ be the full and the reduced cross-sectional algebras of the \fb \mb respectively. On one hand we prove that $C^*(\tmax{\ma}{\mb})\cong\tmax{C^*(\ma)}{C^*(\mb)}$, and in the other hand we show that also 
$C^*_r(\tmin{\ma}{\mb})\cong\tmin{C^*_r(\ma)}{C^*_r(\mb)}$, which
reflects the harmony between universal constructions on one hand and
between spatial ones on the other hand. Perhaps it is appropriate to comment here that, in reality, the most useful results we obtain refer to these two norms. However, we have tried to develop the theory in general, which could be useful for example if a theory of “exotic tensor products”, in the style of exotic crossed products, were to be developed in the future. 
\par In the final section we consider some applications. We consider \fbs with
certain \aps and we prove that these \aps are preserved by taking
tensor products. We show how to apply our results to prove the nuclearity or exactness of cross-sectional C*-algebras of Fell bundles under suitable conditions.
\par This paper corresponds to the first part of \cite{favthesis}, and is an expanded version of the previous work ``\textit{Tensor products  
of \fbs over discrete groups}''
(\texttt{http://xxx.if.usp.br/abs/funct-an/9712006}), which circulated
as a preprint, and where only Fell bundles over discrete groups were considered. It should also be mentioned that in his 2017 book \cite{pds}, Exel developed a minimal theory of tensor products between C*-algebras and Fell bundles over discrete groups. 

\section{C*-trings and Fell bundles}\label{sec:prelim}
\par In the first two parts of this section we will recall from \cite{z}, \cite{env} and \cite{trings} some aspects of the theory of C*-ternary rings and their tensor products that will be needed later. Since in \cite{trings} the context is more general than that of tensor products, we have tried to outline the proofs concerning to our setting, mainly those leading to Theorem~\ref{thm:correspondence}. The occasion will also serve to prove some new results and to introduce some of the notation to be used later. In the third part of the section we will begin the preparation for defining tensor products of Fell bundles in the next section.  
\subsection{C*-trings and the functors of Zettl}\label{subsec:z}
\par A $*$-ternary ring, or $*$-tring for short, is a complex vector space $E$ with a map (called $*$-ternary product)
$\mu :E\times E\times E\to E$, which is linear in the odd variables and conjugate linear in the second one, and satisfies: 
$\mu\big(\mu(x,y,z),u,v\big) =\mu\big(x,\mu(u,z,y),v\big) =\mu\big(x,y,\mu(z,u,v)\big)$, $\forall x,y,z,u,v\in E$. A C*-seminorm on $E$ is a seminorm that satisfies $\norm{\mu(x,y,z)}\leq\norm{x}\norm{y}\norm{z}$, and $\norm{\mu(x,x,x)}=\norm{x}^3$ $\forall x,y,z\in E$. A $*$-tring $E$ with a C*-norm making it a Banach space is called a C*-ternary ring, or just a C*-tring. In general we write just $(x,y,z)$ instead of $\mu(x,y,z)$. Note that if $(E,\mu)$ is a C*-tring, its opposite $E^{\textrm{op}}:=(E,-\mu)$ also is a C*-tring.
\par In \cite{z} Zettl proved that if $E$ is a C*-tring, then there exist a C*-algebra $E^r$ (unique up to isomorphism) and an $E^r$-valued sesquilinear map $\pr{\ }{\ }_r:E\times E\to E^r$ such that $E$ is a right $E^r$-module and $\pr{\ }{\ }_r$ satisfies all the properties of a right inner product except possibly that of positivity, with $(x,y,z)=x\pr{y}{z}_r$, and $\norm{x}^2=\norm{\pr{x}{x}_r}$ $\forall  x,y,z\in E$, and in addition $\textrm{span}\{\pr{y}{z}_r:y,z\in E\}$ is dense in $E^r$. Moreover, he showed that, if $E_+:=\{x\in E: \pr{x}{x}_r\in E^{r,+}\}$ and $E_-:=\{x\in E: -\pr{x}{x}_r\in E^{r,+}\}$ (here $E^{r,+}$ is the positive cone of the C*-algebra $E^r$), then $E_+$ and $E_-$ are sub-C*-trings of $E$ such that $\pr{E_+}{E_-}=0$ and $E=E_+\oplus E_-$ as C*-trings. When $E=E_+$ we say that $E$ is a positive C*-tring (so in this case the sesquilinear map $\pr{\ }{\ }_r$ is an inner product). When $E=E_-$, so $E$ is the opposite of a positive C*-tring, we say that $E$ is a negative C*-tring. Besides, $(E_+,\pr{\ }{\ }_r)$ and $(E_-,-\pr{\ }{\ }_r)$ are full right Hilbert modules over $(E_+)^r$ and $(E_-)^r$ respectively, and $E^r=(E_+)^r\oplus(E_-)^r$ as C*-algebras. Note that, conversely, Hilbert modules provide examples of C*-trings: if $(F,\pr{\ }{\ })$ is a right Hilbert module, and we define $\mu(x,y,z):=x\pr{y}{z}$, then both $(F,\mu)$ and $(F,-\mu)$ are C*-trings, the former positive.     
\par Actually, C*-trings are the objects of a category, which we denote $\mathsf{Ct}$, in which the morphisms are linear maps $\pi:E\to F$ that preserve the ternary product, that is $\pi(x,y,z)=(\pi x,\pi y,\pi z)$, $\forall x,y,z\in E$. As shown in \cite{env}, in this case there exists a unique homomorphism $\pi^r:E^r\to F^r$ such that
\begin{equation}\label{eqn:piar}
  \pi^r(\pr{x}{y}_r)=\pr{\pi x}{\pi y}\quad\forall x,y\in E,
\end{equation}
so the the correspondence $E\mapsto E^r$ is in fact the object part of a functor from the category $\mathsf{Ct}$ of C*-trings to the category $\mathsf{C}$ of C*-algebras. In particular, if $E$ is a full right Hilbert module over the C*-algebra $A$, and we define on $E$ the ternary product $(x,y,z):=x\pr{y}{z}_A$ as above, then we have an isomorphism $E^r\cong A$, such that $\pr{x}{y}_r\mapsto \pr{x}{y}_A$, $\forall x,y\in E$. It is easily seen that, as is the case with homomorphisms of C*-algebras, morphisms of C*-trings are automatically contractive and have closed range, and are isometric exactly when they are injective. In passing, we note that a C*-algebra is also a C*-tring with the $*$-ternary product given by $(x,y,z):=xy^*z$. Then any homomorphism of C*-algebras is also a morphism or $*$-ternary rings, so the category of C*-algebras embedds into the category of C*-trings.          
\par Finally, let us mention that, just as we have a Zettl functor $\big(E\stackrel{\pi}{\to}F\big)\mapsto \big(E^r\stackrel{\pi^r}{\to}F^r\big)$ on the right, we also have one on the left: $\big(E\stackrel{\pi}{\to}F\big)\mapsto \big(E^l\stackrel{\pi^l}{\to}F^l\big)$. Of course, here $E^l$ is a C*-algebra and we have an $E^l$-sesquilinear map $\pr{\ }{\ }_l:E\times E\to E^l$ such that $E$ is a left $E^l$-module and $\pr{\ }{\ }_l$ satisfies all the properties of a left inner product, except possibly that of positivity, with $(x,y,z)=\pr{x}{y}_lz$, $\forall  x,y,z\in E$, and also $E^l=\overline{\textrm{span}}\{\pr{y}{z}_r:y,z\in E\}$. Combining both Zettl functors we conclude that a positive C*-tring $E$ is an $E^l-E^r$ Morita-Rieffel equivalence bimodule. In fact, if $E=E_+\oplus E_-$ is the Zettl's decomposition of $E$, then $E^p:=E_+\oplus (E_-)^{\textrm{op}}$ is a positive C*-tring, and we have $E^r=(E^p)^r$ and $E^l=(E^p)^l$, so an arbitrary C*-tring $E$ is close to being an equivalence bimodule, and in any case its associated C*-algebras $E^l$ and $E^r$ are Morita-Rieffel equivalent. For this reason, many properties of these C*-algebras can be considered as inherited from the C*-tring. This is the case of nuclearity and exactness for example, as shown in \cite{trings}.

\subsection{Tensor products of C*-ternary rings and of Hilbert modules}\label{subsec:tenstrings}
\par Suppose that $E$ and $F$ are right Hilbert modules over the C*-algebras $A$ and $B$ respectively. Then one can form its \textit{exterior tensor product} $\bts{E}{}{F}$, which is a right Hilbert module over the C*-algebra $\bts{A}{}{B}$, where the latter is the spatial tensor product of $A$ and $B$ (see \cite{l}). However, as shown below, it is possible to make the same construction using other tensor products between $A$ and $B$ and without major modifications.. 
\par In what follows we denote by $\ms\mn(E)$ and by $\mn(E)$ the sets of C*-seminorms and C*-norms respectively on the $*$-tring or $*$-algebra $E$. Note that $\ms\mn(E)$ is a partially ordered set with the pointwise order: $\gamma_1\leq \gamma_2\iff \gamma_1(x)\leq \gamma_2(x)$ $\forall x\in E$.  
\par Recall that if $(G,\norm{\ })$ is a seminormed space, and $N:=\{x\in G:\norm{x}=0\}$, the Hausdorff completion of $G$ is the completion of the quotient space $G/N$ with respect to the quotient norm $\norm{x+N}:=\norm{x}$.
\par Suppose that $A$ is a $*$-algebra and that $\alpha$ is a C*-seminorm on $A$. Then the Hausdorff completion of $A$ is a C*-algebra, which we denote by $A_\alpha$. Let $p_\alpha:A\to A_\alpha$ be the canonical map. If $A_\alpha^+$ is the set of positive elements of the C*-algebra $A_\alpha$, the set $p_\alpha^{-1}(A_\alpha^+)$ is a cone in $A$, whose elements will be called $\alpha$-positive elements of $A$. Note that if $\alpha\geq \beta$ are C*-seminorms, then the identity $id:(A,\alpha)\to (A,\beta)$ is continuous, so it defines a surjective homomorphism of C*-algebras $\sigma^\alpha_\beta:A_\alpha\to A_\beta$ such that $p_\beta=\sigma^\alpha_\beta p_\alpha$. Thus any $\alpha$-positive element is also a $\beta$-positive element.
\begin{df}\label{df:positiveelements}
Let $A$ be a $*$-algebra. We define the set of positive elements of $A$ to be the set $A^+:=\cap_{\alpha\in\ms\mn(A)}p_{\alpha}^{-1}(A_\alpha)$, where $A_\alpha$ is the Hausdorff completion of $A$ with respect to the C*-seminorm $\alpha$.    
\end{df}

\par Note that elements of the set $C_A:=\{\sum_{i=1}^na_i^*a_i: n\in\nt, a_1,\ldots,a_n\in A\}$ are $\alpha$-positive, $\forall \alpha\in\ms\mn(A)$. 

\begin{rk}\label{rk:positiveelements}
  If $\mn(A)\neq\emptyset$, then $A^+:=\cap_{\alpha\in\mn(A)}p_{\alpha}^{-1}(A_\alpha)$, that is, we only need C*-norms rather than C*-seminorms to determine the positive elements. To see this, let $\beta$ be any C*-seminorm on $A$, and $\alpha$ a C*-norm that it is supposed to exist on $A$. Then $\gamma:=\max\{\alpha,\beta\}\in\mn(A)$ and $\gamma\geq\beta$. Therefore, as observed before the definition, $p_{\gamma}^{-1}(A_\gamma)\subseteq p_{\beta}^{-1}(A_\beta)$. Then $\cap_{\alpha\in\mn(A)}p_{\alpha}^{-1}(A_\alpha)\subseteq\cap_{\beta\in\ms\mn(A)}p_{\beta}^{-1}(A_\beta)\subseteq \cap_{\alpha\in\mn(A)}p_{\alpha}^{-1}(A_\alpha)$, so they are equal. 
\end{rk}

\par We will need the following result, which is exactly \cite[Lemma~4.3]{l}, except that the C*-norm considered here is arbitrary, while Lance's version is only stated for the minimal norm. Since the proof is also the same, we omited it. 

\begin{lem}\label{lem:lance}
Let $A$ and $B$ be \css , and suppose that 
    $\mathfrak{a}=\big(a_{ij}\big)$, 
    $\mathfrak{c}=\big(c_{ij}\big)\in M_n(A)$,
    $\mathfrak{b}=\big(b_{ij}\big)$, 
    $\mathfrak{d}=\big(d_{ij}\big)\in M_n(B)$. Let \bts{A}{\alpha}{B} 
    be a $C^*$-tensor product of $A$ and $B$. Then: 
    \begin{enumerate}
    \item If $0\leq\mathfrak{a}\leq\mathfrak{c}$ and
      $0\leq\mathfrak{b}\leq\mathfrak{d}$, we have $0\leq \big(a_{ij}\otimes
      b_{ij}\big)\leq \big(c_{ij}\otimes d_{ij}\big)$ in 
      $M_n(\bts{A}{\alpha}{B})$. 
    \item If $\mathfrak{a}$, $\mathfrak{b}\geq 0$, then 
      $\sum_{i,j=1}^na_{ij}\otimes b_{ij}\geq 0$ in \bts{A}{\alpha}{B}. 
    \end{enumerate}
\end{lem}

Let $E$ and $F$ be right Hilbert modules over the C*-algebras $A$ and $B$, and let $E\odot F$ and $A\odot B$ their corresponding algebraic tensor products. Using the universal property of the algebraic tensor product, we easily see that \bats{E}{F} is a right module over \bats{A}{B} and that we have an \bats{A}{B}-valued sesquilinear form on \bats{E}{F}. On elementary tensors the action and the form are given by $(x\odot y)(a\odot b)=xa\odot yb$ and $\pr{x\odot y}{x'\odot y'}=\pr{x}{x'}_E\odot\pr{y}{y'}_F$. We want to see that this sesquilinear form is positive, that is, that $\pr{z}{z}\in (A\odot B)^+$ according to Definition~\ref{df:positiveelements}.

\begin{prop}\label{prop:postimesposispos}
  The sesquilinear map $\pr{\,}{}:(\bats{E}{F})\times(\bats{E}{F})\to\bats{A}{B}$, given by $\pr{z}{z'}=\sum_{i=1}^n\sum_{j=1}^m \pr{x_i}{x'_j}_E\odot\pr{y_i}{y'_j}_F$ for $z=\sum_{i=1}^nx_i\odot y_i$, $z'=\sum_{j=1}^mx'_j\odot y'_j$ is positive, that is  
$\pr{z}{z}\in (A\odot B)^+$ $\forall z\in \bats{E}{F}$.   
\end{prop}
\begin{proof}
  Since $\mn(\bats{A}{B})\neq\emptyset$, by Remark~\ref{rk:positiveelements} it is enough to show that $\pr{z}{z}\in (A\otimes_\alpha B)^+$ for every $\alpha\in \mn(\bats{A}{B})$ (in fact it would be enough to do so just for $\norm{\ }_{\max}$, but the proof is the same).
  \par So let $\alpha$ be any C*-norm on $\bats{A}{B}$, and $z=\sum_{i=1}^nx_i\odot y_i\in\bats{E}{F}$. By \cite[Lemma~4.2]{l} the Gramian matrices $X=(\pr{x_i}{x_j}_E)$ and $Y=(\pr{y_i}{y_j}_F)$ are positive elements of $M_n(A)$ and $M_n(B)$ respectively. Therefore $\pr{z}{z}=\sum_{i,j=1}^n\pr{x_i}{x_j}_E\odot \pr{y_i}{y_j}_F$ is a positive element of $\bts{A}{\alpha}{B}$ by (2) of Lemma~\ref{lem:lance}, which ends the proof.    
\end{proof}

Let $\alpha$ be a C*-norm on \bats{A}{B}. Since the sesquilinear map just defined $\pr{\,}{}:(\bats{E}{F})\times(\bats{E}{F})\to\bats{A}{B}$ is positive, we can perform the double completion process described in \cite[top of page 5]{l} to obtain a Hilbert module \bts{E}{\tilde{\alpha}}{F}, which is the completion of \bats{E}{F} with respect to the norm $\tilde{\alpha}:\bats{E}{F}\to\R$ given by

\begin{equation}\label{eqn:tildeh}
\tilde{\alpha}(z):=\sqrt{\alpha(\pr{z}{z})},\quad \forall z\in \bats{E}{F}.
\end{equation}
\begin{rk}
  Lance proves along \cite{l} that for $z$ as in Proposition~\ref{prop:postimesposispos} we have $z=0$ in $\bats{E}{F}$ if and only if $\pr{z}{z}=0$ in $\bats{A}{B}$. This shows that the sesquilinear maps above are actually inner products. 
  \end{rk}
\medskip
\begin{df}\label{df:alphaext}
  We call the right Hilbert $\bts{A}{\alpha}{B}$-module \bts{E}{\tilde{\alpha}}{F} the $\alpha$-exterior product corresponding to the C*-norm $\alpha\in\mn(\bats{E}{F})$.
\end{df}
\par Note that \bts{E}{\tilde{\alpha}}{F} is full whenever $E$ and $F$ are full Hilbert modules.
\medskip
\par We turn again to the C*-trings perspective. Suppose that $E$ and $F$ are positive C*-trings, so they are full right Hilbert modules over the C*-algebras $E^r$ and $F^r$ respectively. Note that $\bats{E}{F}$ has a structure of $*$-tring with the ternary product given by $(x\odot y,x'\odot y',z\odot z'):=(x,y,z)\odot(x',y',z')$ on elementary tensors, which in terms of our just defined sesquilinear form can be written as $(x\odot y,x'\odot y',z\odot z')=(x\odot y)\pr{x'\odot y'}{x''\odot y''}$. So we have just seen that every C*-norm $\alpha$ on $\bats{E^r}{F^r}$ defines a C*-norm $\tilde{\alpha}$ on the $*$-tring \bats{E}{F} (given by \eqref{eqn:tildeh}), whose completion is the positive C*-tring \bts{E}{\tilde{\alpha}}{F}, and $(\bts{E}{\tilde{\alpha}}{F})^r$ turns out to be $\bts{E^r}{\alpha}{F^r}$ (recall \eqref{eqn:piar} and subsequent comments). 
\par Suppose conversely that $\gamma$ is a C*-norm on the $*$-tring \bats{E}{F}, and let $\bts{E}{\gamma}{F}$ be the corresponding completion, which is a C*-tring. Let $E_0^r:=\textrm{span}\{\pr{x}{x'}_E:x,x'\in E\}$ and $F_0^r:=\textrm{span}\{\pr{y}{y'}_F:y,y'\in F\}$. Then $E_0^r$ and $F_0^r$ are dense two-sided ideals of $E^r$ and $F^r$ respectively. Let $z=\sum_{i=1}^nx_i\odot y_i\in\bats{E}{F}$ and $c:=\sum_{j=1}^m\pr{x'_j}{x''_j}_E\odot\pr{y'_j}{y''_j}_F\in\bats{E^r}{F^r}$.  
We have  
\begin{gather*}
  zc=\sum_{j=1}^m\sum_{i=1}^nx_i\pr{x'_j}{x''_j}_E\odot y_i\pr{y'_j}{y''_j}_F
            =\sum_{j=1}^m\sum_{i=1}^n(x_i,x'_j,x''_j)_E\odot (y_i,y'_j,y''_j)_F \\ 
            =\sum_{j=1}^m(\sum_{i=1}^n(x_i\odot y_i,x'_j\odot y'_j,x''_j\odot y''_j)
            =\sum_{j=1}^m(z,x'_j\odot y'_j,x''_j\odot y''_j).
\end{gather*}            
So, since $\gamma$ is a C*-norm: $\gamma(zc)=\gamma(\sum_{j=1}^m(z,x'_j\odot y'_j,x''_j\odot
y''_j) )
\leq \sum_{j=1}^m\gamma(x'_j\odot y'_j)\gamma(x''_j\odot y''_j)\gamma(z)$, $\forall z\in\bats{E}{F}$. Therefore the action of multiplication by $c$ is $\gamma$-bounded on $\bats{E}{F}$, and hence it extends to a bounded operator on $\bts{E}{\gamma}{F}$.
In fact, recalling that $(x\odot y,x'\odot y',x''\odot y'')=(x\odot y)\pr{x'\odot y'}{x''\odot y''}$, where the latter is the inner product that Zettl's associates to the C*-tring $\bts{E}{\gamma}{F}$, we can continue our computations above, and get:

\begin{gather*}
  zc=\sum_{j=1}^m(z,x'_j\odot y'_j,x''_j\odot y''_j)
    =\sum_{j=1}^mz\pr{x'_j\odot y'_j}{x''_j\odot y''_j}\\
    =z\big(\sum_{j=1}^m\pr{x'_j\odot y'_j}{x''_j\odot y''_j}\big).             
  \end{gather*}

Thus we get an injective homomorphism of $*$-algebras $\bats{E^r}{F^r}\to (\bts{E}{\gamma}{F})^r$, given by $c=\sum_{j=1}^m\pr{x'_j}{x''_j}_E\odot\pr{y'_j}{y''_j}_F\mapsto \sum_{j=1}^m\pr{x'_j\otimes y'_j}{x''_j\otimes y''_j}$.  
So we can define on $\bats{E^r_0}{F^r_0}$ the operator norm, namely $\gamma^r:\bats{E^r_0}{F^r_0}\to\R$ such that 
\begin{equation}\label{eqn:oph}
\gamma^r(c):=\sup\{\gamma(zc):z\in\bats{E}{F},\gamma(z)\leq 1\}. 
\end{equation}

Now observe that, since $E_0^r$ and $F_0^r$ are dense ideals in $E^r$ and $F^r$ respectively, this C*-norm uniquely extends to a C*-norm on \bats{E^r}{F^r} (because of \cite[Lemma~5.12]{trings}) by the same formula \eqref{eqn:oph}.   

In conclusion, given two positive C*-trings, we have two maps
\begin{gather*}
    \Psi_r:\mn(\bats{E^r}{F^r})\to \mn(\bats{E}{F})\textrm{ such that } \alpha\mapsto \tilde{\alpha} \textrm{ given by }\eqref{eqn:tildeh}
\\
  \Phi_r:\mn(\bats{E}{F})\to \mn(\bats{E^r}{F^r})\textrm{ such that } \gamma\mapsto \gamma^r \textrm{ given by }\eqref{eqn:oph}.
\end{gather*}
And these correspondences satisfy
\begin{equation}\label{eqn:errseytildes}
(\bts{E}{\tilde{\alpha}}{F})^r=\bts{E^r}{\alpha}{F^r}\qquad and  \qquad  \bts{E^r}{\gamma^r}{F^r}=(\bts{E}{\gamma}{F})^r.
\end{equation}
\par It is easily checked that $\Psi_r$ is order preserving and $\Psi_r\Phi_r$ is the identity on $\mn(\bats{E}{F})$. Moreover, due to the uniqueness of the C*-algebra $E^r$, $\Phi_r\Psi_r$ is the identity on $\mn(\bats{E^r}{F^r})$. Finally, $\Phi_r$ also is order preserving: if $\gamma_1\geq \gamma_2$ are C*-norms on $\bats{E}{F}$, then $id:(\bats{E}{F},\gamma_1)\to (\bats{E}{F},\gamma_2)$ is continuous, and therefore it extends by continuity to a \hm of C*-trings $\pi:\bts{E}{\gamma_1}{F}\to \bts{E}{\gamma_2}{F}$, which induces a \hm of C*-algebras $\pi^r:\bts{E^r}{\gamma_1^r}{F^r}\to \bts{E^r}{\gamma_2^r}{F^r}$, thus contractive; therefore $\gamma_1^r\geq \gamma_2^r$. In conclusion the maps $\Phi_r$ and $\Psi_r$ are mutually inverse isomorphisms between the posets $\mn(\bats{E}{F})$ and  $\mn(\bats{E^r}{F^r})$. We record this fact: 

\begin{thm}\label{thm:correspondence}
  Let $E$ and $F$ be positive C*-trings. Then the maps $\Psi_r:\mn(\bats{E^r}{F^r})\to \mn(\bats{E}{F})$ such that $\alpha\mapsto \tilde{\alpha}$, given by \eqref{eqn:tildeh}, and $\Phi_r:\mn(\bats{E}{F})\to \mn(\bats{E^r}{F^r})$ such that $\gamma\mapsto \gamma^r$, given by \eqref{eqn:oph}, are mutually inverse isomorphisms of partially ordered sets.
  Moreover, if $\alpha\in\mn(\bats{E^r}{F^r})$ and $\gamma\in\mn(\bats{E}{F})$, then $\bts{E}{\tilde{\alpha}}{F}$ and $\bts{E}{\gamma}{F}$ are full right Hilbert modules over $\bts{E^r}{\alpha}{F^r}$ and $\bts{E^r}{\gamma^r}{F^r}$ respectively, so $\big(\bts{E}{\tilde{\alpha}}{F}\big)^r\cong\bts{E^r}{\alpha}{F^r}$ and $\big(\bts{E}{\gamma}{F}\big)^r\cong\bts{E^r}{\gamma^r}{F^r}$, where the isomorphisms extend the map $\pr{x\odot y}{x'\odot y'}\mapsto\pr{x}{x'}\odot\pr{y}{y'}$, $\forall x,x'\in E$, $y,y'\in F$.   
\end{thm}
\par In fact in \cite{trings} it is proved that the correspondences above extend to isomorphisms between $\ms\mn(\bats{E}{F})$ and $\ms\mn(\bats{E^r}{F^r})$ for abitrary C*-trings.
\par In particular, since $\Psi_r$ and $\Phi_r$ are order preserving maps, and $\mn(\bats{E^r}{F^r})$ has a maximum and a minimum elements $\norm{\ }_{\max}$ and $\norm{\ }_{\min}$ respectively, we have: 

\begin{cor}\label{cor:maxmin}(cf \cite[Corollary~5.13]{trings}).
Let $E$ and $F$ be positive \cts. Then there exist a maximum $C^*$-norm 
$\norm{\cdot}_{\max}$ on \bats{E}{F}, and a minimum $C^*$-norm 
$\norm{\cdot}_{\min}$ on \bats{E}{F}, and  
\begin{gather*}  
  \big(\bts{E}{\max}{F}\big)^l=\bts{E^l}{\max}{F^l},
\hspace*{1cm} 
\big(\bts{E}{\max}{F}\big)^r=\bts{E^r}{\max}{F^r},
\\
\big(\bts{E}{\min}{F}\big)^l=\bts{E^l}{\min}{F^l}
\hspace*{1cm} 
\big(\bts{E}{\min}{F}\big)^r=\bts{E^r}{\min}{F^r}.
\end{gather*}
Recall that the minimum norm on the $*$-algebras \bats{E^r}{F^r} and \bats{E^l}{F^l} agrees with the so called spatial one.
\end{cor}
\begin{rk}\label{rk:crossnorms}
  Let $\alpha$ be a C*-norm on $\bats{E^r}{F^r}$. The well-known fact that $\alpha$ is cross-norm, that is $\alpha(a\odot b)=\norm{a}_{E^r}\norm{b}_{F^r}$ $\forall a\in E^r$ and $b\in F^r$, implies that $\tilde{\alpha}$ also is cross-norm, for if $x\in E$, $y\in F$:
\begin{gather*}
  \tilde{\alpha}(x\odot y)^2=\alpha(\pr{x\odot y}{x\odot y})
                            =\alpha(\pr{x}{x}\odot\pr{y}{y})\\ 
                            =\norm{\pr{x}{x}}_E\norm{\pr{y}{y}}_F
                            =\tilde{\alpha}(x)^2\tilde{\alpha}(y)^2. 
\end{gather*}
\end{rk}  

\medskip
\par To end this part of the section we prove the following two results, which will be needed later. 
 To prove the first of them, recall first from \cite[Proposition~4.1]{env} that if $\pi:E\to F$ is an injective homomorphism of C*-trings, then the induced \hm of C*-algebras $\pi^r:E^r\to F^r$ also is injective. We remark that the converse is easily proved as well: if $\pi(x)=0$, then $0=\pr{\pi(x)}{\pi(x)}=\pi^r(\pr{x}{x})$, so $\pr{x}{x}=0$ if $\pi^r$ is injective, and in this case $x=0$.   

\begin{prop}\label{prop:biminimal}
  Let $\pi_1:E_1\to F_1$ and $\pi_2:E_2\to F_2$ be homomorphisms of positive C*-trings. Then $\pi_1\odot\pi_2:\bats{E_1}{E_2}\to \bats{F_1}{F_2}$ is $\norm{\ }_{\textrm{min}}$-continuous, so it extends to a homomorphism $\tmin{\pi_1}{\pi_2}:\tmin{E_1}{E_2}\to \tmin{F_1}{F_2}$. Moreover, if $\pi_1$ and $\pi_2$ are injective, then so is $\tmin{\pi_1}{\pi_2}$. 
\end{prop}
\begin{proof}
  Applying the right Zettl functor we obtain homomorphisms $\pi_1^r:E_1^r\to F_1^r$ and $\pi_2^r:E_2\to F_2^r$, so by \cite[T.5.19]{w} we get a homomorphism $\tmin{\pi_1^r}{\pi_2^r}:\tmin{E_1^r}{E_2^r}\to \tmin{F_1^r}{F_2^r}$, which by definition extends $\bats{\pi_1^r}{\pi_2^r}:\bats{E_1^r}{E_2^r}\to \bats{F_1^r}{F_2^r}$. 
Now, if $z=\sum_{i=1}^nx_i\odot y_i\in \bats{E}{F}$: 
  \begin{gather*}
    \pr{(\pi_1\odot\pi_2)z}{(\pi_1\odot \pi_2)z}
    =\sum_{i,j=1}^n\pr{\pi_1(x_i)\odot \pi_2(y_i)}{\pi_1(x_j)\odot \pi_2(y_j)}\\
    =\sum_{i,j=1}^n\pr{\pi_1(x_i)}{\pi_1(x_j)}\odot\pr{\pi_2(y_i)}{\pi_2(y_j)}\\
    =\sum_{i,j=1}^n\pi_1^r(\pr{x_i}{x_j})\odot\pi_2^r(\pr{y_i}{y_j})
    =(\pi_1^r\odot\pi_2^r)(\pr{z}{z}).
  \end{gather*}
  Therefore, since $\norm{(\pi_1\odot\pi_2)z}_{\textrm{min}}^2
  =\norm{\pr{(\pi_1\odot\pi_2)z}{(\pi_1\odot \pi_2)z}}_{\textrm{min}}$, we get 
  \[\norm{(\pi_1\odot\pi_2)z}_{\textrm{min}}^2
    =\norm{(\pi_1^r\odot\pi_2^r)(\pr{z}{z})}_{\textrm{min}}
    =\norm{(\pi_1^r\otimes_{\textrm{min}}\pi_2^r)(\pr{z}{z})}_{\textrm{min}}
  \leq\norm{z}_{\textrm{min}}^2,\] which ends the proof of the first statement. As for the second one, it follows from the last assertion of \cite[T.6.9]{w} and the remark preceding the present Proposition.   
\end{proof}

  In the same way, but using \cite[T.6.9]{w} instead of \cite[T.5.19]{w}, we obtain 

  \begin{prop}\label{prop:bimaximal}
  Let $\pi_1:E_1\to F_1$ and $\pi_2:E_2\to F_2$ be homomorphisms of positive C*-trings. Then $\pi_1\odot\pi_2:\bats{E_1}{E_2}\to \bats{F_1}{F_2}$ is $\norm{\ }_{\textrm{max}}$-continuous, so it extends to a homomorphism $\tmax{\pi_1}{\pi_2}:\tmax{E_1}{E_2}\to \tmax{F_1}{F_2}$. 
\end{prop}

\par It follows from Propositions~\ref{prop:biminimal} and \ref{prop:bimaximal} that both the minimal and maximal tensor products are bifunctors $\mathsf{Ct}\times\mathsf{Ct}\to\mathsf{Ct}$, where $\mathsf{Ct}$ is the category of C*-trings. 

\begin{lem}\label{lem:lance2}
  Let $E$ and $F$ be full right Hilbert modules over the C*-algebras $A$ and $B$ respectively, and $S\in\ml(E)$, $T\in\ml(F)$. If $\gamma$ is a C*-norm on \bats{E}{F}, then the map $S\odot T:\bats{E}{F}\to \bats{E}{F}$ is $\gamma$-continuous, with $\norm{S\odot T}\leq\norm{S}\norm{T}$, and it extends to an adjointable map $S\otimes T\in\ml(\bts{E}{\gamma}{F})$, whose adjoint is $S^*\otimes T^*$.
\end{lem}
\begin{proof}
  First recall that the C*-norm $\gamma$ induces a C*-norm $\gamma^r$ on \bats{A}{B}, namely the operator norm
  \[\gamma^r(c)=\sup\{\gamma(zc):z\in\bats{E}{F}:\gamma(z)\leq 1\}.\]
  Let $z=\sum_{i=1}^nx_i\odot y_i\in\bats{E}{F}$, and consider the matrices $X=(\pr{x_i}{x_j})$, $X_S=(\pr{Sx_i}{Sx_j})$, $Y=(\pr{y_i}{y_j})$ and $Y_T=(\pr{Ty_i}{Ty_j})$. By \cite[Lemmas~4.1 and 4.2]{l}, we see that $0\leq X_S\leq\norm{S}^2 X$ and $0\leq Y_T\leq\norm{T^2}Y$. Then, using the last assertion of Lemma~\ref{lem:lance}, we get:
  \begin{gather*}
    \pr{(S\odot T) z}{(S\odot T) z}=\sum_{i,j=1}^n\pr{Sx_i\odot Ty_i}{Sx_j\odot Ty_j}\\
    =\sum_{i,j=1}^n\pr{Sx_i}{Sx_j}\odot\pr{Ty_i}{Ty_j}\\
    \leq \sum_{i,j=1}^n\norm{S}^2\pr{x_i}{x_j}\odot\norm{T}^2\pr{y_i}{y_j}
    =\norm{S}^2\norm{T}^2\pr{z}{z}.
  \end{gather*}
  Therefore:
  \begin{gather*}
  \gamma((S\odot T)z)^2=\gamma^r(\pr{(S\odot T) z}{(S\odot T) z})
  \leq\gamma^r(\norm{S}^2\norm{T}^2\pr{z}{z})\\
  =\norm{S}^2\norm{T}^2\gamma^r(\pr{z}{z})
  =\norm{S}^2\norm{T}^2\gamma(z)^2. 
\end{gather*}
We conclude that $S\odot T$ is bounded, with $\norm{S\odot T}\leq\norm{S}\norm{T}$ as claimed. Thus $S\odot T$ extends by continuity to a map $S\otimes T$. It is now easy to verify that $S\otimes T$ is adjointable, and that $(S\otimes T)^*=S^*\otimes T^*$.    
  \end{proof}

  \begin{cor}\label{cor:alfita}
    Let $E$ and $F$ be Hilbert modules, and $\gamma$ be a C*-norm on \bats{E}{F}. Then there exists a (unique) C*-norm $\overline{\gamma}$ on $\bats{\ml(E)}{\ml(F)}$ such that $\bts{\ml(E)}{\overline{\gamma}}{\ml(F)}$ is a C*-subalgebra of $\ml(\bts{E}{\gamma}{F})$. In case $\gamma=\norm{\ }_{\textrm{min}}$, also is
    $\overline{\gamma}=\norm{\ }_{\textrm{min}}$.   
  \end{cor}
  \begin{proof}
    It follows from Lemma~\ref{lem:lance2} that we have a $*$-\hm $\varphi_\gamma:\bats{\ml(E)}{\ml(F)}\to \ml(\bts{E}{\gamma}{F})$. In case $\gamma=\norm{\ }_{\textrm{min}}$, in \cite[page 37]{l} it is shown that this homomorphism extends to an isometric \hm $\varphi_{\textrm{min}}:\tmin{\mathcal{L}(E)}{\mathcal{L}(F)}\to \mathcal{L}(\tmin{E}{F}$ (which in particular proves our last statement). Suppose $T\in \bats{\ml(E)}{\ml(F)}$ is such that $\varphi_{\gamma}(T)=0$. Since $\varphi_\gamma$ and  $\varphi_{\textrm{min}}$ agree on $\bats{E}{F}$, the fact $\varphi_{\gamma}(T)|_{\bats{E}{F}}=0$, implies $\varphi_{\textrm{min}}(T)|_{\bats{E}{F}}=0$ and, since $\bats{E}{F}$ is dense in $\tmin{E}{F}$, this entails $\varphi_{\textrm{min}}(T)=0$. Since $\varphi_{\textrm{min}}$ is injective, we conclude that $T=0$. Consequently $\varphi_\gamma$ is injective, and therefore, identifying $\bats{\ml(E)}{\ml(F)}$ with the *-subalgebra $\varphi_\gamma(\bats{\ml(E)}{\ml(F)})$, it is enough (and necessary) to take $\overline{\gamma}$ as the restriction of the norm on $\ml(\bts{E}{\gamma}{F})$ to $\bats{\ml(E)}{\ml(F)}$.  
\end{proof}

\subsection{Positive $*$-algebraic bundles and
  Fell bundles}\label{subsec:algebraic}  
\begin{df}\label{df:pre}
Let $G$ be a discrete group, and suppose that $\mc=(C_t)_{t\in G}$ is
a family of complex vector spaces. We identify \mc with the disjoint
union of the spaces $C_t$. We then say that $\mc$ is a $*$-algebraic  
bundle over $G$, with product $\cdot :\mc\times\mc\to \mc$ and
involution $*:\mc\to \mc$ if, $\forall a,b\in\mc$, $t,s\in G$, the
following holds:
\bc
\begin{tabular}{lll}
1) $C_sC_t\subseteq C_{st}$&\hspace*{.3cm}& 
5) $(C_t)^*\subseteq C_{t^{-1}}$\\ 
2) The product $\cdot$ is bilinear on $C_s\times C_t\to  C_{st}$&&
6) $(ab)^*=b^*a^*$.\\ 
3) The product on $\mc$ is associative.&&
7) $a^{**}=a$.\\
4) $*$ is conjugate linear from $C_t$ into $C_{t^{-1}}$. && 
\end{tabular}
\ec
 The vector spaces $C_t$ are called the fibers of the bundle. Note that each $C_t$ is a $*$-tring with the product $(a,b,c):=ab^*c$, and in particular $C_e$ is a $*$-algebra (here and in the rest of the paper $e$ will denote the unit of a
group).    
\end{df}
\par Suppose that $\mc=(C_t)$ is a $*$-algebraic bundle over $G$, and
that $\mi=(I_t)$ is a subset of \mc such that \mi is also a
$*$-algebraic bundle with the operations inherited from \mc, which
moreover satisfies $\mc\mi\subseteq \mi$ and $\mi\mc\subseteq
\mi$. Then we say that $\mi$ is a (two-sided) ideal of $\mc$. It is
easy to see $\mc/\mi:=(C_t/I_t)$ is also a $*$-algebraic bundle with
the obviuos operations naturally induced on the quotients by the
operations on $\mc$.      

\begin{df}\label{df:posfb}
Let $\mc=(C_t)_{t\in G}$ be a $*$-algebraic bundle, and $\alpha$ a
$C^*$-seminorm on $C_e$. We say that $\mc$ is an $\alpha$-positive
$*$-algebraic bundle if for each $c\in\mc$ the element $c^*c$ is
positive in the Hausdorff completion $(C_e)_\alpha$ of $C_e$. We say
that $\mc$ is a 
positive $*$-algebraic bundle if it is $\alpha$-positive $\forall
\alpha\in\ms\mn(C_e)$.
\par In other words, $\mc$ is positive if $c^*c\in C_e^+$ in the meaning of $C_e^+$ according to Definition~\ref{df:positiveelements}.  
\end{df}

\begin{df}\label{df:normonfb}
Let $\mc=(C_t)_{t\in G}$ be a $*$-algebraic bundle. Let
$\norm{\cdot}:\mc\to \R$ be such that:  
\begin{enumerate}
\item[8)] $(C_t,\norm{\cdot})$ is a seminormed space, $\forall t\in
  G$. 
\item[9)] $\norm{c_1c_2}\leq\norm{c_1}\,\norm{c_2}$ $\forall
  c_1,c_2\in \mc$. 
\item[10)] $\norm{c^*c}=\norm{c}^2$.
\end{enumerate}
We then say that $\norm{\ }$ is a $C^*$-seminorm on $\mc$, and that it
is a $C^*$-norm if each $(C_t,\norm{\cdot})$ is a normed space. We
represent respectively by $\ms\mn(\mc)$ and $\mn(\mc)$ the sets of of
$C^*$-seminorms and $C^*$-norms on $\mc$.  
\par If, moreover, 
\begin{enumerate}
\item[11)] $\mc$ is $\alpha$-positive, where $\alpha$ is the
  restriction of $\norm{\ }$ to the unit fiber $C_e$,
\ee
we will say that $(\mc,\cdot,*,\norm{\cdot})$ is a semi-pre-\fb 
over the discrete group $G$, and that it is a pre-Fell bundle if
$\norm{\ }$ is a $C^*$-norm. A pre-Fell bundle $\mc$     
is called a Fell bundle if each $(C_t,\norm{\cdot})$ is complete for
all $t\in G$. 
\end{df}
Note that 9) and 10) imply that $\norm{c^*}=\norm{c}$, $\forall c\in \mc$, and also that $\norm{\ }$ is a C*-norm on the $*$-tring $C_t$. 
\par The proof of the following result is routine, and it is left to
the reader.
\begin{prop}
Let $\mc^0=(C_t^0)_{t\in G}$ be a pre-\fb over the discrete group   
$G$, with $C^*$-norm $\norm{\cdot}$. For $t\in G$, let $C_t$ be the  
completion of $C_t^0$, and consider the family of Banach spaces 
$(C_t)_{t\in G}$ with the extended norm. Then the product and
involution on $\mc^0$ extend by continuity to $\mc$, and with the extended operations and norm $\mc$ is a Fell bundle  
over $G$. We say that $\mc$ is a completion of the pre-\fb $\mc^0$.
\end{prop}

\par Given a semi-pre-Fell bundle $\mc=(C_t)_{t\in G}$, let
$\mi:=\{x\in \mc: \norm{x}=0\}$. Then \mi can be identified with the
$*$-algebraic bundle $\mi=(I_t)_{t\in G}$, where $I_t:=\mi\cap C_t$,
$\forall t\in G$. Note that $\mi$ is also an ideal of $\mc$, for
property 9) above implies $\mc\mi\subseteq \mi$ and $\mi\mc\subseteq
\mi$. It is easy to check that $\mc/\mi:=(C_t/I_t)_{t\in G}$ is a
pre-\fb with the norm induced by the seminorm on \mc: if $c\in C_t$,
then $\norm{c_t+I_t}:=\norm{c_t}$. We will say 
that the Fell bundle $\mc_{\norm{\,}}$ obtained by completing this pre-\fb $\mc/\mi$ is
the \textit{Hausdorff completion of $\mc$}.    

\begin{df}\label{df:farrows}
Let $\ma=(A_t)_{t\in G}$ and $\mb=(B_t)_{t\in G}$ be *-algebraic
bundles over the discrete group $G$. A \hm  $\phi:\ma\to\mb$ is a map 
such that $\phi (A_t)\subseteq B_t$, $\forall t\in G$, and, $\forall
a,b\in A$, $t\in G$: 
1) $\phi\r{A_t}:A_t\to B_t$ is linear; 2) $\phi (ab)=\phi (a)\phi
(b)$; 3) $\phi (a^*)=\phi (a)^*$.
If $\ma$ and $\mb$ are semi-pre-\fbs we also require that $\phi$ is
continuous on each fiber $A_t$.
\end{df}
\par Note that if $\ma$ and $\mb$ are semi-pre-\fbs 
and $\phi :\ma\to \mb$ is a \hm of *-algebraic bundles, then $\phi$   
is continuous if and only if 
$\phi :A_e\to  B_e$ is continuous, because if $x\in\ma$, then 
\[\norm{\phi (x)}^2=\norm{\phi (x)^*\phi (x)}
                   =\norm{\phi (x^*x)}\leq\norm{\phi\r{A_e}}\,\norm{x^*x}
                   =\norm{\phi\r{A_e}}\,\norm{x}^2.\] 
\par In particular every \hm of *-algebraic bundles between \fbs over
discrete groups is continuous. Observe also that, with the notion of
\hm just introduced, any two Hausdorff completions of a given semi-pre-\fb
are necessarily isomorphic, and therefore the Hausdorff completion of
a semi-pre-\fb is essentially unique. 
\par If $\beta$ is a $C^*$-seminorm on the $*$-algebraic bundle
$\mc=(C_t)$, and $\alpha:=\beta|_{C_e}$, it is clear that $\alpha\in
\ms\mn(C_e)$. Besides, by properties 9) and 10) we have:  
\begin{gather}\label{eqn:extensible}
\alpha(c^*c)=\beta(c^*c)=\beta(c)^2=\beta(cc^*)=\alpha(cc^*), \quad\forall c\in
\mc.\end{gather} 
A natural question that arises is whether a $C^*$-seminorm on
$C_e$ can be extended to a $C^*$-seminorm on $\mc$. It follows that if $\mc$ is
an $\alpha$-positive $*$-algebraic bundle, the
necessary condition \eqref{eqn:extensible} is also sufficient for this
to be true:
\begin{prop}\label{prop:tpextnorms}
 Let $\mc=(C_t)_{t\in G}$ be a $*$-algebraic bundle over
 the discrete group $G$, and $\alpha\in \ms\mn(C_e)$ such that $\mc$ is $\alpha$-positive. Then $\alpha$ can be extended to a C*-seminorm on $\mc$ if and only if $\alpha$ satisfies the relation \eqref{eqn:extensible} above. In this case its extension is given by $\tilde{\alpha}:\mc\to [0,\infty)$ such that 
 $\tilde{\alpha}(c):=\sqrt{\alpha(c^*c)}$, $\forall c\in \mc$. Moreover
  $\tilde{\alpha}\in\mn(\mc)\iff\alpha\in \mn(C_e)$.  
\end{prop}       
\begin{proof}
Each fiber $C_t$ is a right module over $C_e$, and $\pr{\
}{\,}_r^t:C_t\times C_t\to C_e$ such that $\pr{c}{d}_r^t:=c^*d$ is a
right semi-inner product on $C_t$ because $\mc$ is
$\alpha$-positive. Then $\tilde{\alpha}|_{C_t}$ is a seminorm on 
$C_t$ (see \cite[page 3]{l} or \cite[Proposition~3.30]{trings}).
Therefore we have that $\alpha(\pr{c}{d}_r^t)\leq \tilde{\alpha}(c)\tilde{\alpha}(d)$ and $\tilde{\alpha}(ca)\leq \tilde{\alpha}(c)\alpha(a)$ $\forall c,d\in C_t$, $a\in C_e$ (see for instance \cite[Proposition~3.30]{trings}). 
Similarly, $C_t$ is a left module over $C_e$, and $\pr{\ }{\,}_l^t:C_t\times C_t\to C_e$ such that $\pr{c}{d}_l^t:=cd^*$ is a left semi-inner product on $C_t$, which induces the seminorm $\Tilde{\tilde{\alpha}}$ such that $\Tilde{\tilde{\alpha}}(c):=\alpha(cc^*)$ $\forall c\in C_t$, and we have $\alpha(\pr{c}{d}_l^t)\leq \Tilde{\tilde{\alpha}}(c)\Tilde{\tilde{\alpha}}(d)$ and $\Tilde{\tilde{\alpha}}(ac)\leq \alpha(a)\Tilde{\tilde{\alpha}}(c)$ $\forall c,d\in C_t$, $a\in C_e$. 
Now suppose that \eqref{eqn:extensible} holds for $\alpha$, that is $\tilde{\alpha}=\Tilde{\tilde{\alpha}}$. Then, if $c\in C_s,d\in C_t$, recalling the above inequalities and observing that $c^*c\in C_e$, we have:
\begin{gather*}
  \tilde{\alpha}(cd)^2=\alpha(d^*c^*cd)=\alpha(\pr{d}{c^*cd}_r^t)
  \stackrel{}{\leq} \tilde{\alpha}(d)\tilde{\alpha}(c^*cd)
  =\tilde{\alpha}(d)\Tilde{\tilde{\alpha}}(c^*cd)\\
  \leq \tilde{\alpha}(d)\alpha(c^*c)\Tilde{\tilde{\alpha}}(d)
  = \tilde{\alpha}(c)^2\tilde{\alpha}(d)\Tilde{\tilde{\alpha}}(d) 
  =\tilde{\alpha}(c)^2 \tilde{\alpha}(d)^2.
\end{gather*}
On the other hand:
$\tilde{\alpha}(c^*c)=\sqrt{\alpha(c^*cc^*c)}=\sqrt{\alpha(c^*c)^2}=\tilde{\alpha}(c)^2$. We 
conclude that $\tilde{\alpha}$ satisfies properties 8)--10) of
Definition~\ref{df:normonfb}, so it is a $C^*$-seminorm on
$\mc$. The converse has already been observed, and the last statement is clear.     
\end{proof}
\par Since $C_t$ can be considered as both a right and a left $C_e$-module, condition \eqref{eqn:extensible} expresses the fact that the C*-seminorms induced on $C_t$ in both cases by the C*-norms $\alpha$ agree.   
\par As in the case of $*$-algebras, the sets $\ms\mn(\mc)$ and
$\mn(\mc)$ of $C^*$-seminorms and $C^*$-norms on a $*$-algebraic
bundle $\mc$ are partially ordered sets. Moreover, the considerations above lead to consider also the (partially ordered) sets:  
\begin{gather*}
  \ms\mn_{\mc}(C_e):=\{\alpha\in\ms\mn(C_e):\,\alpha (c^*c)=\alpha(cc^*)\,\forall c\in \mc\}\\
  \mn_{\mc}(C_e):=\ms\mn_{\mc}(C_e)\cap\mn(C_e)
\end{gather*}    
   
\begin{thm}\label{thm:tpextnorms}
Let $\mc=(C_t)_{t\in G}$ be a positive $*$-algebraic bundle over the
discrete group $G$. Then the map $\Phi:\ms\mn(\mc)\to
\ms\mn_{\mc}(C_e)$ given by $\Phi(\beta):=\beta|_{C_e}$ is an
isomorphism of partially ordered sets, whose inverse $\Psi$ is given
by $\Psi(\alpha)=\tilde{\alpha}$, where
$\tilde{\alpha}(c):=\sqrt{\alpha(c^*c)}$, $\forall c\in \mc$. Besides:
$\Phi(\mn(\mc))=\mn_{\mc}(C_e)$.
\end{thm}
\begin{proof}
It is clear that both $\Phi$ and $\Psi$ are order preserving, and
Proposition~\ref{prop:tpextnorms} shows that $\Phi\circ\Psi
=Id_{\ms\mn(C_e)}$. The fact that $\Psi\circ\Phi=Id_{\ms\mn(\mc)}$
follows directly from the definition of $\tilde{\alpha}$ and property
10) in Definition~\ref{df:normonfb}.
\end{proof}
\par Again as in the case of $*$-algebras, note that if $\alpha\geq\beta$ are C*-seminorms on the $*$-algebraic bundle $\mc$, then every $\alpha$-positive element of $\mc$ is $\beta$-positive as well. Moreover, the indentity on $\mc$ induces a surjective homomorphism of Fell bundles $\sigma^\alpha_\beta:\mc_\alpha\to\mc_\beta$. 
\par We end the section with the definition of general Fell bundles and related concepts.  
\par A \textit{\fb} (or \textit{$C^*$-algebraic bundle}) $\mb=(B_t)_{t\in G}$ over the
locally compact group $G$ is a Banach bundle $\mb$ over $G$, with
fiber $B_t$ over $t\in G$, and such that there exist continuous
product and involution defined on \mb and satisfying conditions
1)--11) of Definition \ref{df:pre}. Recall that a \textit{Banach
bundle} (\cite[II-13.4]{fd}) over a Hausdorff space $X$, called base
space, is a pair $(\mb,p)$ formed by a Hausdorff space $\mb$, called
total space, and a continuous open surjection $p:\mb\to X$, together
with continuous maps $\norm{\ }:\mb\to\R$, $+:\{ (b,b')\in \mb\times
\mb:\, p(b)=p(b')\}\to \mb$ and $\C\times \mb\to \mb$ such that each
fiber $B_x:=p^{-1}(\{ x\})$ becomes a complex Banach space with  
the restrictions of these maps, and such that it satisfies the 
additional property:
\textit{ if $x\in X$ and $(b_i)\subseteq \mb$ is a net such that
      $p(b_i)\to x$ and $\norm{b_i}\to 0$, then $b_i\to 0_x\in \mb$},
where $0_x$ is the zero element of $B_x$.    
A \textit{\hm of Banach bundles} $\phi:\ma\to\mb$ over $X$ is a 
continuous map such that $\phi_x:=\phi\r{A_x}:A_x\to B_x$ is a bounded 
linear operator, $\forall x\in X$, and $\norm{\phi}:=\sup_{x\in
X}\norm{\phi_x}<\infty$.  
\par A section of \mb is a function $\xi:X\to\mb$ such that 
$p(\xi(x))=x$, $\forall x\in X$. If $\xi,\eta$ are continuous sections 
of \mb, and $\al\in\C$, then $t\mapsto \al\xi(t)+\eta(t)$ is again a
continuous section. We will denote by $C_c(\mb)$ the vector space of
continuous sections of compact support of the Banach bundle \mb. If
$K\subseteq X$ is a compact subset, we denote by $C_K(\mb)$ the
subspace of $C_c(\mb)$ whose elements are those with support contained
in $K$. The map $\norm{\ }_K:C_K(\mb)\to\R$ given by  
$\norm{\xi}_K=\max_{x\in X}{\xi(x)}$ is a norm and
$(C_K(\mb),\norm{\ }_K)$ is a Banach space. We endow $C_c(\mb)$ with
the locally convex inductive limit topology induced by the family
$\{(C_K(\mb),\iota_K)\}_{K}$, where $K$ runs over the family of
compact subsets of $X$, and $\iota_K:C_K(\mb)\inc C_c(\mb)$ is the
natural inclusion. We refer the reader to \cite{fd} for further
information on Banach bundles.    
\par If $X$ is a topological space, $X_d$ will denote the set $X$ with
the discrete topology and, if \mb is a Banach bundle over $X$, we will
denote by $\mb_d$ the Banach bundle over $X_d$ whose fiber over $x\in
X$ is the corresponding fiber of $\mb$. That is,  
$\mb_d$ is the disjoint union of the fibers $B_x$, $x\in X$. Since \mb
is a topological space the notation just introduced is 
ambiguous. Thus, in order to avoid any confusion we will use 
calligraphic letters only to represent Banach bundles. Note that
if \ma is a \fb over $G$, then $\ma_d$ is a \fb over $G_d$.   

\begin{df}\label{df:farrows2}
Let \fela  and $\mb =(B_t)_{t\in G}$ be \fbs over the locally compact
group $G$. We say that a \hm of Banach bundles $\phi :\ma\to\mb$ is a
\hm of \fbs if $\phi :\ma_d\to\mb_d$ is a \hm of \fbs over $G_d$  
(see Definition \ref{df:farrows}). 
\end{df}

\par Along this work we will use repeatedly the following two 
results. The first one is \textit{Cohen-Hewitt theorem: if $B$ is a
Banach algebra with approximate unit and if $E$ is a non-degenerate
Banach $B$-module (i.e. $\cgen{EB}=E$), then for each $x\in E$ there exist 
$y\in E$, $b\in B$, such that $x=yb$}. Although the use of this
theorem is not strictly necessary for our purposes, it 
facilitates the exposition and allows us to avoid the repetition of 
similar approximation arguments. A proof of this theorem may 
be found in \cite{fd} (there is a nice proof for Hilbert modules in
\cite{rw}). 
\par The second of the mentioned results is the theorem of
\textit{Douady-dal Soglio H\'erault}, which is fundamental in the
theory of Banach bundles: \textit{let $X$ be a Hausdorff 
space, and $(\mb,p)$ a Banach bundle over $X$; if $X$ is paracompact
or locally compact, then for each $b\in\mb$ there exists a continuous
section of compact support $\xi$ of $\mb$ such that
$\xi\big(p(b)\big)=b$}. The reader is referred to \cite[Apendix C]{fd}
for a proof.

\section{Tensor Products of Fell Bundles}\label{sec:tensfell}
\par Our aim in what follows is to introduce tensor products of Fell
bundles. A tensor product of the Fell bundles 
$\ma=(A_t)_{t\in G}$ and $\mb=(B_s)_{s\in H}$ over the groups 
$G$ and $H$ will be a Fell bundle $\mc=(C_r)_{r\in G\times H}$ over 
$G\times H$, and we will have that $C_e$ is a tensor product of $A_e$
and $B_e$ (recall that $e$ denotes the unit of the group). We will show that there exist, up to isomorphisms, unique tensor products $\mc_{\max}$  
and $\mc_{\min}$ of $\ma$ and $\mb$, such that $(\mc_{\max})_e= 
\tmax{A_e}{B_e}$ and $(\mc_{\min})_e=\tmin{A_e}{B_e}$. 
\par In the first part of the section we consider the case of bundles over
discrete groups. The treatment of the general case is postponed to 
the the second part of the section. Finally, the end of the section is devoted to study the \reps of tensor products.   

\subsection{Tensor products of Fell bundles over discrete
groups}\label{subsec:dtens}  
\par Let $\ma =(A_t)_{t\in G}$ and $\mb =(B_s)_{s\in H}$ be Fell
bundles over the groups $G$ and $H$ respectively. Consider, for $t\in
G$, $s\in H$, the algebraic tensor product \bats{A_t}{B_s}. When we
let $t$, $s$ run in $G$ and $H$, we obtain a family
$\{\bats{A_t}{B_s}\}_{(t,s)\in G\times H}$ of vector
spaces. Let denote by $\bats{\ma}{\mb}$ the disjoint union of these
vector spaces. For $(t,s),(t',s')\in G\times H$, we have unique linear
maps $\bigl(\bats{A_t}{B_s}\bigr)\times
\bigl(\bats{A_{t'}}{B_{s'}}\bigr)\to\bats{A_{tt'}}{B_{ss'}}$ such that
$(a_t\odot b_s,a_{t'}\odot b_{s'})\mapsto\, a_ta_{t'}\odot b_sb_{s'}$,
and unique conjugate linear maps
$\bats{A_t}{B_s}\to\bats{A_{t^{-1}}}{B_{s^{-1}}}$  
such that $a_t\odot b_s\mapsto\, a_t^*\odot b_s^*$. Put together,
these families of maps define a product $\cdot
:\big(\bats{\ma}{\mb}\big)\times\big(\bats{\ma}{\mb}\big)\to
\big(\bats{\ma}{\mb}\big)$ and an involution  
$*:\big(\bats{\ma}{\mb}\big)\to 
\big(\bats{\ma}{\mb}\big)$ such that the product is associative,
bilinear on every
$\big(\bats{A_t}{B_s}\big)\times\big(\bats{A_{t'}}{B_{s'}}\big)\to  
\bats{A_{tt'}}{B_{ss'}}$, * is conjugate linear when restricted to 
$\bats{A_t}{B_s} \to \bats{A_{t^{-1}}}{B_{s^{-1}}}$ and $(x\cdot
y)^*=y^*\cdot x^*$, $\forall x,y\in\bats{\ma}{\mb}$. In other words,
\bats{\ma}{\mb} is a *-algebraic bundle, in the sense of
Definition~\ref{df:pre}.  We will say that \bats{\ma}{\mb} is the
\textit{algebraic tensor product of $\ma$ and $\mb$}. 

\begin{prop}\label{prop:positive}
The algebraic tensor product of Fell bundles is a positive
$*$-algebraic bundle (Definition~\ref{df:posfb}). 
\end{prop}
\begin{proof}
Let \ma and \mb be Fell bundles over $G$ and $H$ respectively. We have
to show that for any $s\in G$, $t\in H$, and elements
$a_1,\ldots,a_n, a_1',\ldots,a_n'\in A_t$,
$b_1,\ldots,b_n,b_1',\ldots,b_n'\in B_t$, the element
$(a_1^*a_1'+\cdots +a_n^*a_n')\odot (b_1^*b_1'+\cdots +b_n^*b_n')$ is a
positive element of $A_e\bigodot B_e$. Since $A_s$ and $B_t$ are
positive $C^*$-trings, and Hilbert modules over $A_e$ and $B_e$, this fact follows from Proposition~\ref{prop:postimesposispos}.    
\end{proof}
\begin{df}\label{df:ftp}
Let $\ma =(A_t)_{t\in G}$ and $\mb =(B_s)_{s\in H}$ be \fbs over the
discrete groups $G$ and $H$, and consider their algebraic tensor product 
\bats{\ma}{\mb}. If $\alpha$ is a $C^*$-norm on \bats{\ma}{\mb}, 
we will call the completion \bts{\ma}{\alpha}{\mb} of  
$(\bats{\ma}{\mb},\alpha)$ a tensor product of $\ma$ and $\mb$.
\end{df}

\par If $\bts{\ma}{\alpha}{\mb}$ is a tensor product of $\ma$ and
$\mb$, then the unit fiber $\big(\bts{\ma}{\alpha}{\mb}\big)_e$ is a tensor product of   
$A_e$ and $B_e$. In fact, if we know the $C^*$-norm determined by  
$\big(\bts{\ma}{\alpha}{\mb}\big)_e$ on \bats{A_e}{B_e}, 
then we know the norm of every $x\in \bts{\ma}{\alpha}{\mb}$, because
it is equal to $\sqrt{\alpha (x^*x)}$. Therefore, two tensor products
will be isomorphic if and only if their fibers on the identity element 
are the same tensor product of $A_e$ and $B_e$. 
This raises the question of whether or not a given tensor product of   
$A_e$ and $B_e$ determines a tensor product of the \fbs 
$\ma$ and $\mb$. According to Proposition~\ref{prop:tpextnorms}, if
$\alpha$ is a $C^*$-norm on \bats{A_e}{B_e}, then $\alpha$ can be 
extended to a $C^*$-norm on \bats{\ma}{\mb} 
if and only if $\alpha (x^*x)=\alpha (xx^*)$, 
$\forall x\in \bats{\ma}{\mb}$, and in this case the extension is
unique. Writing $x=\sum_{i=1}^na_i\odot b_i\in A_r\odot B_s$, this condition 
is $\alpha(\sum_{i,j=1}^nx_i^*x_i\odot y_i^*y_i)=\alpha(\sum_{i,j=1}^nx_ix_i^*\odot y_iy_i^*)$. Although we will not go deeper into this problem, 
we will see that this is in fact the case for the maximal and minimal   
tensor products (see Proposition~\ref{prop:extmm} below).    
We begin with a result certainly well-known; for lack of reference we
include a proof of it.     

\begin{lem}\label{lem:maxideals}
Let $I$ and $J$ be ideals of the \css  $A$ and $B$ 
respectively. Then \tmax{I}{J} is the closure of \bats{I}{J} in
\tmax{A}{B}.  
\end{lem}
\begin{proof} Let $\pi :\tmax{I}{J}\to B(H)$ be a faithful and
non-degenerate \rep of \tmax{I}{J}. Then there are faithful and
non-degenerate \reps $\pi_I:I\to  B(H)$ and $\pi_J:J\to  B(H)$,
such that $\pi_I(x)\pi_J(y)=\pi (x\otimes y)=\pi_J(y)\pi_I(x)$,
$\forall x\in I$, $y\in J$ (\cite[T.6.4]{w}). Since $\pi_I$ and $\pi_J$
are non-degenerate they have unique extensions $\pi_A:A\to B(H)$ and  
$\pi_B:B\to  B(H)$ to \reps of $A$ and $B$ respectively  
(\cite[VI-19.11]{fd}). If $a\in A$, $x\in I$, $b\in B$ and $y\in J$,
then $\pi_A(ax)\pi_B(by)=\pi_B(by)\pi_A(ax)$, because $ax\in I$ and
$by\in J$. Since $\pi_I$ and $\pi_J$ are non-degenerate, we conclude
that $\pi_A(a)\pi_B(b)=\pi_B(b)\pi_A(a)$, $\forall a\in A$, $b\in 
B$. Hence there exists a \rep $\tilde{\pi}:\tmax{A}{B}\to  B(H)$ such
that $\tilde{\pi}(a\otimes b)=\pi_A(a)\pi_B(b)$, $\forall a\in A,\,
b\in B$. Thus $\tilde{\pi}$ is an extension of
$\pi\r{\bats{I}{J}}$. Since $\tilde{\pi}$ is contractive, we conclude
that if $x\in \bats{I}{J}$, its norm as an element of \tmax{A}{B} is
greater or equal to its norm in \tmax{I}{J}, and therefore they
agree. 
\end{proof}
 
\begin{prop}\label{prop:extmm}
Let $\ma=(A_t)_{t\in G}$ and $\mb=(B_s)_{s\in H}$ be \fbs over 
the discrete groups $G$ and $H$. Then the norms 
$\norm{\cdot}_{\min}$ and $\norm{\cdot}_{\max}$
on \bats{A_e}{B_e} can be extended to 
$C^*$-norms on \bats{\ma}{\mb}. 
\end{prop}
\begin{proof} Let $A_t^*A_t:=\cgen\{a_t^*a_t:\ a_t\in
A_t\}\subseteq A_e$ and  $B_s^*B_s:=\cgen\{b_s^*b_s:\ b_s\in
B_s\}\subseteq B_e$. Then 
$A_t^*A_t$ and $B_s^*B_s$ are ideals in $A_e$ and $B_e$
respectively, and $A_t$ may be seen as a positive \ct with 
$A_t^r=A_t^*A_t$ and $A_t^l=A_tA_t^*$, and  
similarly $B_s$. Recall that there exists a maximum $C^*$-norm 
$\mu$ on $\bats{A_t}{B_s}$. By \cite[Corollary~5.13]{trings}, we must have
$(\tmax{A_t}{B_s})^r=\bts{A_t^*A_t}{\mu^r}{B_s^*B_s}$ and  
$(\tmax{A_t}{B_s})^l=\bts{A_tA_t^*}{\mu^l}{B_sB_s^*}$, where $\mu^r$
denotes the maximum norm on \bats{A_t^*A_t}{B_s^*B_s} and $\mu^l$
denotes the maximum norm on \bats{A_tA_t^*}{B_sB_s^*}. Now, Lemma
\ref{lem:maxideals} implies that $\mu^r$ and $\mu^l$ are 
restrictions of the maximum norm of \bats{A_e}{B_e}. Since 
\tmax{A_t}{B_s} is a Hilbert
\big(\bts{A_tA_t^*}{\mu^l}{B_sB_s^*}-\bts{A_t^*A_t}{\mu^r}{B_s^*B_s}\big)
-bimodule we have, for $x\in \bats{A_t}{B_s}$
\[\norm{xx^*}_{\max}=\norm{xx^*}_{\mu^l}= 
  \norm{x}_{\mu}^2=\norm{x^*x}_{\mu^r}=\norm{x^*x}_{\max}.\] 
Thus $\norm{\ }_{\max}$ may be extended to all of \bats{\ma}{\mb} by
Proposition \ref{prop:tpextnorms}.
\par On the other hand, it is well-known that if $C$ and $D$ are \scss
of the \css $A$ and $B$ respectively, then the restriction of the
spatial norm on \bats{A}{B} to \bats{C}{D} is precisely the spatial
norm on \bats{C}{D} (see for instance \cite[Corollary~B.14]{trings}, or simply Proposition~\ref{prop:biminimal}). Therefore the same arguments given above for 
$\norm{\ }_{\max}$ also apply to the spatial norm on \bats{A_e}{B_e}
and hence $\norm{\ }_{\min}$ can also be extended to \bats{A_e}{B_e}. 
\end{proof}
\subsection{Tensor products of Fell bundles over locally compact groups}\label{subsec:ctens}  
\par We will extend next the construction done in the previous section
to the case of \fbs over arbitrary locally compact groups. 

\par Suppose now that $\ma=(A_t)_{t\in G}$ and $\mb=(B_s)_{s\in H}$
are \fbs over the locally compact groups $G$ and $H$, 
and let \bts{\ma_d}{\al}{\mb_d} be a tensor product of  
$\ma_d$ and $\mb_d$ as in the previous section. We will endow 
\bts{\ma_d}{\al}{\mb_d} with a topology such that
\bts{\ma_d}{\al}{\mb_d} will be a \fb over $G\times H$. 
\par For $f\in C_c(\ma )$, $g\in C_c(\mb)$, let  
$f\oslash g:G\times H\to\bts{\ma_d}{\al}{\mb_d}$ be such that  
$(f\oslash g)(t,s)=f(t)\otimes g(s)$, $\forall t\in G$, $s\in H$. 
Every $f\oslash g$ is a section of \bts{\ma_d}{\al}{\mb_d}. We consider
the vector space \[L:=\gen\{f\oslash g:\, f\in C_c(\ma ), g\in
C_c(\mb )\},\] which is a vector subspace of the space of sections of
\bts{\ma_d}{\al}{\mb_d}. 
The topology we want to define on
\bts{\ma_d}{\al}{\mb_d} is determined by the requirement that every
element of $L$ is a continuous section:   

\begin{prop}\label{prop:bb}
With the notation above we have: 
\be
\item For each $l\in L$, the map $G\times H\to  \mathbb{R}$ 
  such that $(t,s)\mapsto\, \alpha\big(l (t,s)\big)$ is continuous.
\item For each $(t,s)\in G\times H$, the set $L (t,s):=\{l (t,s):\,
  l\in L\}$ is dense in \bts{A_t}{\al}{B_s}.
\item There exists a unique topology on \bts{\ma_d}{\al}{\mb_d} for
  which \bts{\ma_d}{\al}{\mb_d} is a Banach bundle over 
  $(G\times H)$ and such that $L$ is contained in the space of
  continuous sections of the bundle \bts{\ma_d}{\al}{\mb_d} with this 
  topology.  
\ee
The Banach bundle over $G\times H$ thus obtained will be denoted by 
$\bts{\ma}{\al}{\mb}$.
\end{prop}
\begin{proof}
Since 3) is a consequence of 1) and 2) (\cite[II-13.18]{fd}) it is
enough to prove the first two assertions. We begin by 2). If 
$x=\sum_{i=1}^na_i\otimes b_i\in\bts{A_t}{\al}{B_s}$, there exist
continuous sections $f_i$, $ g_i$ of 
\ma and \mb respectively such that $f_i(t)=a_i$, $ g_i(s)=b_i$,
$\forall i=1\ldots n$ (\cite[C.17]{fd}). If $l
=\sum_{i=1}^nf_i\oslash g_i$, then $l\in L$, and 
$ l (t,s)=\sum_{i=1}^nf_i(t)\otimes g_i(s)=\sum_{i=1}^na_i\otimes
b_i=x.$ Hence 2) follows, because \bats{A_t}{B_s} is dense in
\bts{A_t}{\al}{B_s}.  
\par To prove 1), fix $l =\sum_{i=1}^nf_i\oslash g_i\in L$, 
and let $(t,s)\to  (t_0,s_0)$. Then: 
\[ \alpha\big(l (t,s)\big)^2=\alpha\big(\sum_{i=1}^nf_i(t)\otimes g_i(s)\big)^2=
\alpha\big(\sum_{i,j=1}^nf_i(t)^*f_j(t)\otimes g_i(s)^* g_j(s)\big)\]
Note that $t\mapsto\,f_i(t)^*f_j(t)$ and $s\mapsto\, g_i(s)^* g_j(s)$
are continuous maps, because the $f_i$ and $ g_i$ are continuous
sections, and the involutions of \ma and \mb are continuous as
well. Now the ``cross-norm'' property (i.e.: $\alpha(a\otimes b)= 
\alpha(a)\alpha(b)$) of the $C^*$-norms on tensor products implies   
that $a\otimes b\to  a_0\otimes b_0$ when $a\to  a_0$ and $b\to
b_0$. Therefore,   
\[\sum_{i,j=1}^nf_i(t)^*f_j(t)\otimes g_i(s)^* g_j(s)\to  
\sum_{i,j=1}^nf_i(t_0)^*f_j(t_0)\otimes g_i(s_0)^* g_j(s_0)\] 
when $(t,s)\to  (t_0,s_0)$. Thus $\alpha\big(l (t,s)\big)\to  \alpha\big(l
 (t_0,s_0)\big)$ if $(t,s)\to  (t_0,s_0)$.  
\end{proof}

\par If $G$ is a group, * is an involution on a set $X$ and $l: G\to
X$ is a map, we define a new map $\tilde{l}:G\to X$ as 
$\tilde{l}(t)=l(t^{-1})^*$. In particular, if $l$ is a continuous 
section of compact support of a \fb \mb, then $\tilde{l}$ also is. 

\begin{lem}\label{lem:invcont}
The involution  $*:\bts{\ma}{\al}{\mb}\to\bts{\ma}{\al}{\mb}$ is 
continuous.
\end{lem}
\begin{proof}
\par We know from 
\cite[II-13.18]{fd} that a base for the topology defined in
Proposition \ref{prop:bb} is given by the sets 
\[\mathcal{W}(l ,U,\epsilon )=\{w\in\bts{\ma}{\al}{\mb}:\, 
  p(w)\in U, \text{ and } \alpha\big(l\big(p(w)\big)-w\big)<\epsilon\},\] 
where $p:\bts{\ma}{\al}{\mb}\to G\times H$ is the projection,
$U\subseteq G\times H$ is an open subset, $l =\sum_if_i\oslash 
g_i$, with $f_i\in C_c(\ma)$, $ g_i\in C_c(\mb)$, and $\epsilon >0$. 
In other words, $\mathcal{W}(l ,U,\epsilon )=\bigcup_{t\in U}B(l(t),\epsilon)$, where $B(l(t),\epsilon)\subseteq (\bts{\ma}{\al}{\mb})_t$ is the open $\epsilon$-ball with center $l(t)$.  
Then we have:
\begin{align*}
\mathcal{W}(l ,U,\epsilon )^*
 &=\{w^*\in\bts{\ma}{\al}{\mb}:\, p(w)\in U, \text{ and } 
                          \alpha\big(l\big(p(w)\big)-w\big)<\epsilon\}\\
 &=\{w^*\in\bts{\ma}{\al}{\mb}:\, p(w^*)\in U^{-1}, \text{ and } 
                          \alpha\big(l\big(p(w^*)^{-1}\big)^*-w^*\big)<\epsilon\}\\
 &=\{z\in\bts{\ma}{\al}{\mb}:\, p(z)\in U^{-1}, \text{ and } 
                          \alpha\big(\tilde{l}\big(p(z)\big)-z\big)<\epsilon\}\\
 &=\mathcal{W}(\tilde{l} ,U^{-1},\epsilon ),
\end{align*}
Thus $*$ is continuous.  
\end{proof}

\begin{prop}\label{prop:prodcont}
The product $(\bts{\ma}{\al}{\mb})\times (\bts{\ma}{\al}{\mb})\to
\bts{\ma}{\al}{\mb}$ is continuous. 
\end{prop}
\begin{proof}
We claim that if $a\to a_0$ in $\ma$ and $b\to b_0$ in $\mb$, then  
$a\otimes b\to a_0\otimes b_0$ in $\bts{\ma}{\al}{\mb}$. Let
$W\subseteq \bts{\ma}{\al}{\mb}$ be an open set such that $a_0\otimes 
b_0\in W$, and let $f\in C_c(\ma)$, $g\in C_c(\mb)$ be such that
$f(t_0)=a_0$ and $g(s_0)=b_0$.  
Then $(f\oslash g)(t_0,s_0)=a_0\otimes b_0$. Since  
$f\oslash g\in C_c(\bts{\ma}{\al}{\mb})$, and since the norm $\alpha$ is
continuous, there exist $\epsilon >0$ and open sets $U\subseteq G$ and
$V\subseteq H$ such that $(t_0,s_0)\in U\times V$ and
$W\cap(\bts{\ma}{\al}{\mb})_{(t,s)}\supseteq   
B\big((f\oslash g)(t,s),\epsilon \big)$, $\forall (t,s)\in U\times V$.
Consider now the open subsets $W(f,U,\epsilon^{1/2})$ and
$  W(g,V,\epsilon^{1/2})$ of $\ma$ and $\mb$ containing $a_0$ and $b_0$
respectively. We have $
W(f,U,\epsilon^{1/2})\otimes W(g,V,\epsilon^{1/2})
=\{ a_t\otimes b_s\in \bts{\ma}{\al}{\mb}:\, (t,s)\in U\times V,$
and $ \norm{f(t)-a_t}<\epsilon^{1/2},\  \norm{g(s)-b_s}<\epsilon^{1/2}\}
\subseteq\{ x_{(t,s)}\in \bts{\ma}{\al}{\mb}:\, (t,s)\in U\times V,$
and $\alpha\big((f\oslash g)(t,s)-x_{(t,s)}\big)<\epsilon\}
=W(f\oslash g ,U\times V,\epsilon )
=\bigcup_{(t,s)\in U\times V}B\big((f\oslash g)(t,s),\epsilon\big)
\subseteq W,$ so it follows that $a\otimes b\to a_0\otimes b_0$ when $a\to
a_0$, $b\to b_0$. 
\par Note that $\forall f,f'\in C_c(\ma)$, $g,g'\in C_c(\mb)$,
the map $\mu: (G\times H)\times (G\times H)\to
\bts{\ma}{\al}{\mb}$ given by $\big((t,s),(t',s')\big)\mapsto 
(f\oslash g)(t,s)\otimes (f'\oslash g')(t',s')$ is continuous. Indeed  
the products on $\ma$ and $\mb$ are continuous, 
$f,f',g,g'$ are continuous as well, and since $\mu (t,s,t',s')=
f(t)f'(t')\otimes g(s)g'(s')$, the continuity of $\mu$ follows from the
claim at the beginning of the proof.
\par Now pick elements $x_0\in (\bts{\ma}{\al}{\mb})_{(t_0,s_0)}$ and 
 $x_0'\in (\bts{\ma}{\al}{\mb})_{(t_0',s_0')}$, and let $m\in L$, 
$1>\epsilon > 0$ such that $x_0x_0'\in W(m,Z,\epsilon )$, where $Z$ is  
some open subset of $G\times H$ containing $(t_0t_0',s_0s_0')$.  
Let $M>\epsilon+1+\alpha(x_0)+\alpha(x_0')$ and $l,l'\in L$ such that  
$\alpha(l(t_0,s_0)-x_0)<\epsilon/2M$, $\alpha(l'(t_0',s_0')-x_0')<\epsilon/2M$. 
Then we have $\alpha(l(t_0,s_0)l'(t_0',s_0')-x_0x_0')=:d<\epsilon$. Let 
$d<\epsilon' <\epsilon$. As seen above, the map 
$ll':(G\times H)\times (G\times H)\to\bts{\ma}{\al}{\mb}$ such that  
$\big((t,s),(t',s')\big)\mapsto l(t,s)l'(t',s')$ is continuous, so
there exist open neighborhoods $U$, $V$, $U'$ and $V'$ of $t_0$, $s_0$, 
$t_0'$ and $s_0'$ respectively such that 
$ll'\big((U\times V)\times (U'\times V')\big)\subseteq W(m,Z,\epsilon')$. 
Let now $N>1+\norm{l}_{\infty}+\norm{l'}_{\infty}$, 
$0<\del<(\epsilon-\epsilon')/4N$. We have  
$W(l,U\times V,\del)W(l',U'\times V',\del))\subseteq
W(m,Z,\epsilon)$. In fact, if $x_{(t,s)}\in W(l,U\times V,\del)$, 
$x_{(t',s')}'\in W(l,U'\times V',\del)$
\begin{align*}
\alpha(x_{(t,s)}x_{(t',s')}'-m(tt',ss'))
&\leq\alpha(x_{(t,s)}x_{(t',s')}'-l(t,s)l'(t',s'))\\ 
&\hspace*{2.5cm} +\alpha(l(t,s)l'(t',s')-m(tt',ss'))\\
&\leq\epsilon'+\alpha(x_{(t,s)}\big(x_{(t',s')}'-l'(t',s')\big))\\
&\hspace*{2.5cm} +\alpha(\big(x_{(t,s)}-l(t,s)\big)l'(t',s'))\\
&<\epsilon'+\frac{\epsilon-\epsilon'}{4}
\bigg[\frac{1}{N}(\alpha(x_{(t,s)}-l(t,s))\, \alpha(l(t,s)))+1\bigg]\\
&<\epsilon
\end{align*}  
\end{proof}

\begin{df}\label{df:conttensprods}
Let \ma and \mb be \fbs over the locally compact groups $G$ and $H$,  
and let $\al$ be a $C^*$-norm on $\bats{\ma}{\mb}$. 
The tensor product \bts{\ma}{\al}{\mb} of $\ma$ and $\mb$ with respect
to $\al$ is the \fb obtained by completing the algebraic tensor product
\bats{\ma}{\mb} with respect to the $C^*$-norm $\al$,  
furnished with the topology provided by Proposition \ref{prop:bb}.
\end{df}

\begin{prop}\label{prop:extmm2}
Let \ma and \mb be \fbs over the locally compact groups $G$ and $H$.  
If $\al\geq\beta$ are $C^*$-norms on $\bats{\ma}{\mb}$, then there
exists a unique \hm of \fbs 
$\sigma^{\al}_{\beta}:\bts{\ma}{\al}{\mb}\to\bts{\ma}{\beta}{\mb}$  
such that $\sigma^{\al}_{\beta}(a\otimes b)= a\otimes b$, 
$\forall a\in \ma$, $b\in\mb$. This \hm is onto. Moreover, if $\al\geq\beta\geq\gamma$ are $C^*$-norms on $\bats{\ma}{\mb}$, we have $\sigma^{\al}_{\gamma}=\sigma^{\al}_{\beta}\sigma^{\beta}_\gamma$.  
\end{prop}
\begin{proof}
  Since $\alpha\geq\beta$ for each $(r,s)\in G\times H$ the identity map on \bats{A_r}{B_s} has a (unique) continuous extension to a map $\bts{A_r}{\alpha}{B_s}\to \bts{A_r}{\beta}{B_s}$, which is surjective because its image is both dense and closed. The collection of all these maps is clearly a homomorphism $\sigma^{\al}_{\beta}$ from $\bts{\ma_d}{\al}{\mb_d}$ into $\bts{\ma_d}{\beta}{\mb_d}$. It is also continuous from $\bts{\ma}{\al}{\mb}$ into $\bts{\ma}{\beta}{\mb}$, because the vector space $L$ of sections used to define the involved topologies is exactly the same, and the map $\sigma^{\al}_{\beta}$ is the identity on the set of such sections. Thus $\sigma^{\al}_{\beta}$ is continuous by \cite[II-13.16]{fd}. The last assertion follows from the uniqueness of the maps~$\sigma^\alpha_\beta$. 
\end{proof}

\par Summarizing the constructions and results obtained up to now we have:

\begin{thm}\label{thm:minimax}
  Let $\ma=(A_t)_{t\in G}$ and $\mb=(B_s)_{s\in H}$ be \fbs over the locally compact groups $G$ and $H$.  Then $\ms\mn_{\bats{\ma}{\mb}}(\bats{A_e}{B_e})\cong \ms\mn(\bats{\ma}{\mb})$ and $\mn_{\bats{\ma}{\mb}}(\bats{A_e}{B_e})\cong \mn(\bats{\ma}{\mb})$ as a posets.  Moreover $\mn(\bats{\ma}{\mb})$ has a minimum and a maximum elements, namely the unique extensions of $\norm{\cdot}_{\min}$ and $\norm{\cdot}_{\max}$ on \bats{A_e}{B_e} to C*-norms on all of \bats{\ma}{\mb}.  
\end{thm}

\par As a consequence we can extend Propositions \ref{prop:biminimal} and \ref{prop:bimaximal} to the context of Fell bundles:

\begin{prop}\label{prop:biminimalfell}
  Let $\pi_1:\ma_1\to \mb_1$ and $\pi_2:\ma_2\to \mb_2$ be homomorphisms of Fell bundles. Then $\pi_1\odot\pi_2:\bats{\ma_1}{\ma_2}\to \bats{\mb_1}{\mb_2}$ is $\norm{\ }_{\textrm{min}}$-continuous, so it extends to a homomorphism $\tmin{\pi_1}{\pi_2}:\tmin{\ma_1}{\ma_2}\to \tmin{\mb_1}{\mb_2}$. 
\end{prop}
  
\begin{prop}\label{prop:bimaximalfell}
  Let $\pi_1:\ma_1\to \mb_1$ and $\pi_2:\ma_2\to \mb_2$ be homomorphisms of Fell bundles. Then $\pi_1\odot\pi_2:\bats{\ma_1}{\ma_2}\to \bats{\mb_1}{\mb_2}$ is $\norm{\ }_{\textrm{max}}$-continuous, so it extends to a homomorphism $\tmax{\pi_1}{\pi_2}:\tmax{\ma_1}{\ma_2}\to \tmax{\mb_1}{\mb_2}$. 
\end{prop}

Consequently, as in the case of C*-algebras and of C*-trings, we see that the minimal and maximal tensor products of Fell bundles is a bifunctor $\mathsf{F}\times\mathsf{F}\to \mathsf{F}$, where $\mathsf{F}$ is the category of Fell bundles.  

\subsection{Representations of tensor products}\label{subsec:tensreps} 
\par We will study now the \reps of tensor products of \fbs on
Hilbert modules. The results obtained, similar to the case of \css,
will be useful in the next section. The first of them tells us how to
obtain a \rep of $\tmin{\ma}{\mb}$ starting with \reps of $\ma$ and
$\mb$. The second one shows that there exists a bijective
correspondence between non-degenerate \reps of $\tmax{\ma}{\mb}$ and
non-degenerate commuting \reps of $\ma$ and $\mb$. 

\begin{df}\label{df:freps}
Let \ma be a *-algebraic bundle over the discrete group $G$, and \mh a
Hilbert module. A map $\pi :\ma\to \adj{\mh}$ is called a \rep of
$\ma$ on $\mathcal{H}$ if $\pi(ab)=\pi(a)\pi(b)$, $\pi(a^*)=\pi(a)^*$
and $\pi\r{A_t}$ is linear, $\forall a,b\in\ma, t\in G$. The \rep
$\pi$ is said to be non-degenerate if
$\cgen\pi(\ma)\mathcal{H}=\mh$. This is equivalent to the restriction
$\pi\r{A_e}$ to be non-degenerate. 
\end{df}

\begin{df}\label{df:frep} 
Let \ma be a \fb over the locally compact group~$G$. A \rep of $\ma$
on the Hilbert module $\mathcal{H}$ is a \rep
$\phi:\ma_d\to\adj{\mathcal{H}}$ which is strongly continuous, that
is, $\forall h\in\mathcal{H}$ the map $\ma\to\mathcal{H}$ given by 
$a\mapsto\pi(a)h$ is continuous. 
\end{df}

\par Note that for $G$ discrete every \rep of the \fb \ma  
is automatically continuous, because $\norm{\pi(a)}\leq\norm{a}$,
$\forall a\in\ma$, as is easy to check. 
\par If \ma is a \fb (or just an
*-algebraic bundle), and $\mathcal{H}$ is a Hilbert module, we will
denote by $R(\ma ,\mh )$ the family of non-degenerate \reps of $\ma$
on $\mh$. If \mb is 
another \fb (or *-algebraic bundle), we set: 
$R(\ma,\mb,\mh):=\{ (\pi_1,\pi_2)\in R(\ma,\mh)\times R(\mb,\mh):\ 
\pi_1(a)\pi_2(b)=\pi_2(b)\pi_1(a),\, \forall a\in\ma,b\in\mb\}$. 
If $A$ and $B$ are *-algebras, we will also use the notations
$R(A,\mh)$, $R(A,B,\mh)$, with the same meaning.
\medskip
\par In what follows, given right Hilbert modules $\mh$ and $\mk$, over the C*-algebras $C$ and $D$ respectively, we will consider their \textit{exterior tensor product} $\tmin{\mh}{\mk}$, which is a right Hilbert module over $\tmin{C}{D}$. The reader is referred to Subsection~\ref{subsec:tenstrings}, as well as \cite{l} or \cite{trings}.

\begin{prop}\label{prop:minreps}
Let \ma and \mb be \fbs over the locally compact groups $G$ and $H$
respectively, and let $\pi_{\ma}\in R(\ma,\mh_{\ma})$, $\pi_{\mb}\in 
R(\mb,\mh_{\mb})$. Then there exists a unique \rep   
$\pi_{\ma}\otimes\pi_{\mb}\in R(\tmin{\ma}{\mb},\tmin{\mh_{\ma}}{\mh_{\mb}})$ 
such that $(\pi_{\ma}\otimes\pi_{\mb})(a\otimes b)=
\pi_{\ma}(a)\otimes\pi_{\mb}(b)$, 
$\forall a\in\ma$ and $\forall b\in\mb$. If $\pi_{\ma}\r{A_e}$ and  
$\pi_{\mb}\r{B_e}$ are faithful, then 
$(\pi_{\ma}\otimes\pi_{\mb})\r{(\tmin{\ma}{\mb})_e}$ also is
faithful. 
\end{prop}
\begin{proof}
According to \cite[pages 36 and 37]{l} (see also Corollary~\ref{cor:alfita}), we have an isometric embedding    
$\tmin{\adj{\mathcal{H}_{\ma}}}{\adj{\mathcal{H}_{\mb}}}\inc  
\adj{\tmin{\mh_{\ma}}{\mh_{\mb}}},$ such that, $\forall
T\in\adj{\mh_{\ma}}, S\in\adj{\mh_{\mb}}$, $h_{\ma}\in\mh_{\ma}$,
$h_{\mb}\in \mh_{\mb}$:  $(T\otimes S)(h_{\ma}\otimes h_{\mb})=
T(h_{\ma})\otimes S(h_{\mb})$. Thus we may consider, for each  
$(t,s)\in G\times H$, the map $A_t\times B_s\to 
\adj{\tmin{\mh_{\ma}}{\mh_{\mb}}}$ such that 
$(a_t,b_s)\mapsto\pi_{\ma}(a_t)\otimes\pi_{\mb}(b_s)$. This map is  
bilinear, so there exists a unique linear map 
$(\pi_{\ma}\otimes\pi_{\mb})_{(t,s)}:\bats{A_t}{B_s}\to 
\adj{\tmin{\mh_{\ma}}{\mh_{\mb}}}$ 
such that $(\pi_{\ma}\otimes\pi_{\mb})_{(t,s)}(a_t\odot b_s) 
=\pi_{\ma}(a_t)\otimes\pi_{\mb}(b_s)$, $\forall a_t\in A_t$, $b_s\in
B_s$. The collection of these linear maps is a \rep  
$\pi_{\ma}\otimes\pi_{\mb}$ of the pre-\fb \bats{\ma}{\mb}. 
Restricted to \bats{A_e}{B_e} this map coincides with  
$\pi_{\ma}\r{A_e}\otimes\pi_{\mb}\r{B_e}:\bats{A_e}{B_e}\to
\tmin{\adj{\mathcal{H}_{\ma}}}{\adj{\mathcal{H}_{\mb}}}$ 
which is contractive with respect to $\norm{\ }_{\min}$ on
\bats{A_e}{B_e} (\cite[T.5.19]{w}). It follows that 
$\pi_{\ma}\otimes\pi_{\mb}$ extends to a \rep of the \fb 
\tmin{\ma_d}{\mb_d}. Moreover this \rep is continuous in the topology
of $\tmin{\ma}{\mb}$ (recall Definition \ref{df:conttensprods} and
Proposition \ref{prop:bb}). Indeed, for $h_{\ma}\in \mh_{\ma}$ and 
$h_{\mb}\in \mh_{\mb}$ consider the Banach bundle $G\times H\times
(\tmin{\mh_{\ma}}{\mh_{\mb}})$ over $G\times H$ (with the product
topology and the natural projection), and the map 
$\Phi: \tmin{\ma}{\mb}\to G\times H\times (\tmin{\mh_{\ma}}{\mh_{\mb}})$ 
given by $\Phi (c_{t,s})=\left(t,s,(\pi_{\ma}\otimes\pi_{\mb})(c_{t,s})
(h_{\ma}\otimes h_{\mb})\right),\ \forall c_{t,s}\in \tmin{A_t}{B_s}.$ 
Let $L$ be as in Proposition \ref{prop:bb}. To see that $\Phi$ is  
a continuous \hm of Banach bundles it is enough to show, according to
\cite[II-13.16]{fd}, that for all $l\in L$ the map $\Phi l$ is a 
continuous section of the bundle $G\times
H\times(\tmin{\mh_{\ma}}{\mh_{\mb}})$. Clearly it is sufficient to 
check this for sections of the form $f\oslash g$, with 
$f\in C_c(\ma)$, $ g\in C_c(\mb)$. Thus assume that 
$(t,s)\to(t_0,s_0)$ in $G\times H$. We have to show that 
$\Phi\left(f (t)\otimes g (s)\right)\to
\Phi\left(f (t_0)\otimes g (s_0)\right)$, which is equivalent to showing 
that $\pi_{\ma}(f(t))h_{\ma}\otimes\pi_{\mb}( g(s))h_{\mb}$  
converges to $\pi_{\ma}(f(t_0))h_{\ma}\otimes\pi_{\mb}(
g(s_0))h_{\mb}$. Now, if $\varepsilon (t,s)=
\norm{\pi_{\ma}(f(t))h_{\ma}\otimes\pi_{\mb}( g(s))h_{\mb}-
      \pi_{\ma}(f(t_0))h_{\ma}\otimes\pi_{\mb}( g(s_0))h_{\mb}}$, we
have: 
\begin{align*}
 \begin{split}
\varepsilon (t,s)
     &\leq\norm{\pi_{\ma}(f(t))h_{\ma}\otimes\pi_{\mb}( g(s))h_{\mb}-
          \pi_{\ma}(f(t))h_{\ma}\otimes\pi_{\mb}( g(s_0))h_{\mb}}\\
     &\hspace*{1.3cm} +\norm{\pi_{\ma}(f(t))h_{\ma}\otimes
      \pi_{\mb}( g(s_0))h_{\mb}-\pi_{\ma}(f(t_0))h_{\ma}\otimes
      \pi_{\mb}( g(s_0))h_{\mb}}\\
     &\leq\norm{\pi_{\ma}(f(t))}\,\norm{h_{\ma}}\,
          \norm{\pi_{\mb}( g(s))h_{\mb}-\pi_{\mb}( g(s_0))h_{\mb}}\\
     &\hspace*{1.3cm} +\norm{\pi_{\ma}(f(t))h_{\ma}-\pi_{\ma}(f(t_0))h_{\ma}}\, 
          \norm{\pi_{\mb}( g(s_0))h_{\mb}}\\
     &\leq\norm{f}_{\infty}\,\norm{h_{\ma}}\, 
          \norm{\pi_{\mb}( g(s))h_{\mb}-\pi_{\mb}( g(s_0))h_{\mb}}\\
     &\hspace*{1.3cm} +\norm{\pi_{\ma}(f(t))h_{\ma}-\pi_{\ma}(f(t_0))h_{\ma}}\,
          \norm{ g}_{\infty}\,\norm{h_{\mb}},
 \end{split}
\end{align*} 
which converges to zero because $\pi_{\ma}$ and $\pi_{\mb}$ are
continuous \reps.   
\par The fact that $\Phi$ is continuous implies that 
$\forall h_{\ma}\in \mh_{\ma}$, $h_{\mb}\in \mh_{\mb}$, the map 
$\tmin{\ma}{\mb}\to \tmin{\mh_{\ma}}{\mh_{\mb}}$ such that 
$c\mapsto (\pi_{\ma}\otimes\pi_{\mb})(c)
(h_{\ma}\otimes h_{\mb})$ is continuous. Since 
$\norm{(\pi_{\ma}\otimes\pi_{\mb})(c)}\leq \norm{c}$, 
$\forall c\in\tmin{\ma}{\mb}$, we also have that  
$c\mapsto (\pi_{\ma}\otimes\pi_{\mb})(c)(h)$ is continuous, 
$\forall h\in\tmin{\mh_{\ma}}{\mh_{\mb}}$. It follows that 
$\pi_{\ma}\otimes\pi_{\mb}$ is a \rep.
\par If $\pi_{\ma}$, $\pi_{\mb}$ are non-degenerate, then so are 
${\pi_{\ma}}\r{A_e}$ and $\pi_{\mb}\r{B_e}$. 
By Cohen-Hewitt,
given $h_{\ma}\in\mh_{\ma}$, $h_{\mb}\in\mh_{\mb}$, there exist $a\in
A_e$, $b\in B_e$, $h_{\ma}'\in\mh_{\ma}$, $h_{\mb}'\in\mh_{\mb}$ such
that $\pi_{\ma}(a)h_{\ma}'=h_{\ma}$ and  
$\pi_{\mb}(b)h_{\mb}'=h_{\mb}$. Therefore   
$(\pi_{\ma}\otimes\pi_{\mb})(a\otimes b)(h_{\ma}'\otimes h_{\mb}')= 
h_{\ma}\otimes h_{\mb}$. Consequently 
$(\pi_{\ma}\otimes\pi_{\mb}(\tmin{\ma}{\mb}))
(\tmin{\mh_{\ma}}{\mh_{\mb}})$ is dense in
$\tmin{\mh_{\ma}}{\mh_{\mb}}$.  
\par Finally, $(\pi_{\ma}\otimes\pi_{\mb})\r{(\tmin{\ma}{\mb})_e}
=\pi_{\ma}\r{A_e}\otimes\pi_{\mb}\r{B_e}$, and this one is injective
if and only if $\pi_{\ma}\r{A_e}$ and $\pi_{\mb}\r{B_e}$ are
injective (\cite[T.5.19]{w}). 
\end{proof}

\par Suppose that \fela is a \fb and $L,R:\ma\to\ma$
are continuous maps such that there exists $t\in G$ for which
$L(A_s)\subseteq 
A_{ts}$, $R(B_s)\subseteq B_{st}$, $\forall s\in G$, $L\r{A_s}:A_s\to
A_{ts}$, $R\r{A_s}:A_s\to A_{st}$ are linear and bounded, and
$\norm{L}:=\sup_s\norm{L\r{A_s}}<\infty$,  
$\norm{R}:=\sup_s\norm{R\r{A_s}}<\infty$. Then $(L,R)$ is called a  
multiplier of order $t$ of $\ma$ (\cite{fd}) if $\forall
a_1,a_2\in\ma$ the following holds:  
\[
a_1L(a_2)=R(a_1)a_2\qquad
L(a_1a_2)=L(a_1)a_2\qquad
R(a_1a_2)=a_1R(a_2)
\]
The set of multipliers of $\ma$ of order $t$ is denoted by  
$M_t(\ma)$, and $M(\ma)=\bigcup_{t\in G}M_t(\ma)$ denotes the set of
all multipliers of \ma (the notation differs from the one used in
\cite{fd}). Every $M_t(\ma)$ is a Banach space with the obvious
operations and the norm: $\norm{(L,R)}_0=\max\{\norm{L},\norm{R}\}.$
In fact we have $\norm{L}=\norm{R}$. Moreover we have a product
and an involution on $M(\ma)$:
\[
(L_1,R_1)(L_2,R_2)=(L_1L_2,R_2R_1)\qquad
(L,R)^*=(R^*,L^*)
\]
where $L^*(a)=L(a^*)^*$ and $R^*(a)=R(a^*)^*$. With these operations
and norm $M(\ma)$ is a \fb over $G_d$. In addition $M(\ma)$ has a 
topology, in which $u_i=(L_i,R_i)$ converges to $u=(L,R)$ if $\forall
a\in \ma$ we have that $L_i(a)\to L(a)$ and $R_i(a)\to R(a)$. By
analogy to the case of \css, we call this topology \textit{strict} (in 
\cite[VIII-15.2]{fd} this topology is called \textit{strong}). If 
$u=(L,R)\in M(\ma )$, we write $ua$ and $au$ instead of $L(a)$ and
$R(a)$ respectively.
\par There is an isometric and continuous inclusion $\ma\inc M(\ma)$,
given by $a\mapsto (L_a,R_a)$, where $L_a$ is multiplication
by $a$ on the left, and $R_a$ is multiplication by $a$ on the
right. In particular, the topology of $\ma$ is stronger than the
topology inherited from the strict topology of $M(\ma)$. If $A_e$ is
unital, these topologies agree. 
There is also an isomorphism $M(A_e)\cong 
M_e(\ma)$: since $A_t$ is a Hilbert $A_e$-bimodule, then it is also a
Hilbert $M(A_e)$-bimodule, and it can be shown that the actions of
left and right multiplications by elements of $M(A_e)$ on $\ma$ define
multipliers of order $e$ (see \cite[VIII-3.8]{fd}).  
If $\pi :\ma\to\adj{\mh }$ is a non-degenerate \rep of $\ma$,  
then there exists a unique extension (\cite[VIII-15.3]{fd}) of $\pi$
to a \rep $\pi':M (\ma )\to\adj{\mh }$ such that $\forall h\in\mh$, 
the map $u\mapsto\pi'(u)h$ is strictly continuous on cylinders of $\ma$
(the cylinder of radius $r$ of $\ma$ is $C_r:=\{ 
a\in\ma :\, \norm{a}\leq r\}$).

\begin{lem}\label{lem:multis}
The maps $M_t(\ma)\times\ma\to\ma$: $(u,a)\mapsto ua$ and
$(u,a)\mapsto au$ are continuous, $\forall t\in G$.
\end{lem}
\begin{proof}
Recall that for any multiplier $u\in M_t(\ma)$, the maps $\ma\to\ma$:  
$a\mapsto ua$ $a\mapsto au$ are continuous.
Suppose that $(u_i,a_i)\to (u,a)$ in $M_t(\ma)\times\ma$, with
$a_i\in A_{s_i}, a\in A_s$. Since the norm $\norm{\cdot}:\ma\to\R$ is
continuous, there exist $M\geq 0$ and $i_0$ such that $\forall i\geq
i_0$ we have $\norm{a_i}\leq M$. Hence if $i\geq i_0$:
$\max\{\norm{u_ia_i-ua_i},\norm{a_iu_i-a_iu}\}\leq M\, 
\norm{u_i-u}\to 0$, so we have $(u_ia_i-ua_i)\to 0_{ts}$ and
$(a_iu_i-a_iu)\to 0_{st}$  when $i\to\infty$. On the other hand, we 
have that $ua_i\to ua$ and $a_iu\to au$. Thus 
$u_ia_i\to ua$ and $a_iu_i\to au$ if $i\to\infty$. 
\end{proof}

\par The next result is analogous to \cite[Lemma~T.6.1.]{w}.    

\begin{lem}\label{lem:incmult}
Let \fela and \felbh be \fbs and $\bts{\ma}{}{\mb}$ a tensor product  
of $\ma$ and $\mb$. Then there exist unique inclusions 
$\iota_{\ma}: M(\ma )\to M(\bts{\ma}{}{\mb})$ and   
$\iota_{\mb}: M(\mb )\to M(\bts{\ma}{}{\mb})$ such that 
$\iota_{\ma}(u)\iota_{\mb}(v)=\iota_{\mb}(v)\iota_{\ma}(u)$, 
$\forall u\in M(\ma ),\, v\in M(\mb )$, and such that 
$\iota_{\ma}(a)\iota_{\mb}(b)=a\otimes b$, $\forall a\in \ma$, $b\in\mb$.  
These inclusions are isometric and continuous in the strict topologies
when restricted to cylinders.  
\end{lem}
\begin{proof}
Let $u\in M_t(\ma )$. For $r\in G$, $s\in H$, the map 
$A_r\times B_s\to\bts{A_{tr}}{}{B_s}$ such that $(a_r,b_s)\mapsto
(ua_r,b_s)$ is bilinear, and therefore there exists a unique linear
map $L_u:\bats{A_t}{B_s}\to\bts{A_{tr}}{}{B_s}$ such that $a_r\otimes
b_s\mapsto ua_r\otimes b_s$. Similarly, there exists
$R_u:\bats{A_t}{B_s}\to \bts{A_{rt}}{}{B_s}$ such that $R_u(a_r\otimes
b_s)=a_ru\otimes b_s$. The collection of such maps 
define two applications $L_u,R_u:\bats{\ma}{\mb}\to\bts{\ma}{}{\mb}$
such that $\forall x,y\in \bats{\ma}{\mb}\subseteq\bts{\ma}{}{\mb}$
satisfy: $L_u(xy)=L_u(x)y$, $R_u(xy)=xR_u(y)$, $xL_u(y)=R_u(x)y$. If
we prove that $L_u$, $R_u$ are bounded, then they extend by continuity
on each fiber to continuous operators, which still satisfy the above
algebraic relations. In other words, the pair formed by these
extensions will be a multiplier of order $(t,e)$ of
$(\bts{\ma}{}{\mb})_d$.  
\par Let $x=\sum_{i=1}^na_i\otimes b_i\in\bts{A_r}{}{B_s}$. Then:  
$\norm{L_ux}^2=\norm{\sum_{i=1}^nua_i\otimes b_i}^2=
  \norm{\sum_{i,j=1}^na_i^*u^*ua_j\otimes b_i^*b_j}.$ 
Let $\mathfrak{u}=(u_{ij})\in M_n(M(A_e))$, 
$\mathfrak{a}=(a_{ij})\in M_n(A_r)$, 
given by:
$u_{ij}=\begin{cases}\sqrt{u^*u}& \text{ if $i=j$,}\\0& \text{otherwise}
\end{cases}$ and
$a_{ij}=\begin{cases}a_j& \text{ if $i=1$,}\\0& \text{otherwise}
\end{cases}$. Since $A_r$ is a right Hilbert $M(A_e)$-module, then   
$M_n(A_r)$ is a right Hilbert $M_n(M(A_e))$-module. 
Then we have 
$\pr{\mathfrak{u}\mathfrak{a}}{\mathfrak{u}\mathfrak{a}}\leq
\mathfrak{u}^*\mathfrak{u}\pr{\mathfrak{a}}{\mathfrak{a}}\leq
\norm{u}^2\pr{\mathfrak{a}}{\mathfrak{a}}$. Thus  
$\mathfrak{c}:=\norm{u}^2\pr{\mathfrak{a}}{\mathfrak{a}}_r-
\pr{\mathfrak{u}\mathfrak{a}}{\mathfrak{u}\mathfrak{a}}_r\geq 0$.
An easy computation shows that if $\mathfrak{c}=(c_{ij})$, then 
$c_{ij}=a_i^*(\norm{u}^2-u^*u)a_j$. On the other hand, $M_n(B_s)$ is a
Hilbert 
$M_n(M(B_e))$-module. In particular if $\mathfrak{b}=(b_{ij})\in
M_n(B_s)$ is given by $b_{ij}=\begin{cases}b_j& \text{ if $i=1$}\\ 
                                      0& \text{ otherwise
}\end{cases}$, the element $\mathfrak{b}^*\mathfrak{b}=(b_i^*b_j)$ is
positive in $M_n(M(B_e))$. Now, Lemma~\ref{lem:lance} implies that  
$\mathfrak{c}\otimes\mathfrak{b}^*\mathfrak{b}=
\left(a_i^*(\norm{u}^2-u^*u)a_j\otimes b_i^*b_j\right)$ is a positive
element in any $C^*$-completion of $\bats{M_n(M(A_e))}{M_n(M(B_e))}$,  
and $\sum_{i,j=1}^na_i^*(\norm{u}^2-u^*u)a_j\otimes b_i^*b_j$ is
a positive element in any $C^*$-completion of
$\bats{M(A_e)}{M(B_e)}$ (alternatively, the positivity of $\mathfrak{c}\otimes\mathfrak{b}^*\mathfrak{b}$ can be deduced from the proof of \cite[Lemma~4.3]{l}, which does not really use that the norm involved is $\norm{\ }_{\textrm{min}}$). Thus 
$\norm{\sum_{i,j=1}^na_i^*u^*ua_j\otimes b_i^*b_j}\leq\norm{u}^2\, 
\norm{\sum_{i,j=1}^na_i^*a_j\otimes b_i^*b_j}$, for any $C^*$-norm 
on $\bats{A_e}{B_e}$. This shows that $\norm{L_ux}^2\leq\norm{u}^2\, 
\norm{x}^2$, so $L_u$ is bounded. Similarly we see that
$\norm{R_ux}^2\leq \norm{u}^2\,\norm{x}^2$, and therefore $(L_u,R_u)$
extends to a multiplier $\iota_{\ma}(u)$ on $(\bts{\ma}{}{\mb})_d$,
and $\norm{\iota_{\ma}(u)} 
\leq\norm{u}$. In fact $\norm{\iota_{\ma}(u)}=\norm{u}$: if $a\in\ma$,  
$b\in\mb$ are such that $\norm{a}, \norm{b}\leq 1$, then 
$\norm{\iota_{\ma}(u)}\geq\norm{\iota_{\ma}(u)(a\otimes b)}=
\norm{ua}\,\norm{b}=\norm{ua}$, 
and therefore $\norm{\iota_{\ma}(u)}\geq\norm{u}$. Then  
$\norm{\iota_{\ma}(u)}=\norm{u}$, so $\iota_{\ma}$ is an isometry. 
\par To see that $\iota_{\ma}(u)\in M(\bts{\ma}{}{\mb})$, it remains
to prove that it is continuous. To this end consider $f\in
C_c(\ma)$, $g\in C_c(\mb)$. Then the maps $G\times
H\to\bts{\ma}{}{\mb}$ such that $(t,s)\mapsto uf(t)\otimes g(s)$ and 
$(t,s)\mapsto f(t)u\otimes g(s)$ are continuous. 
Suppose that $x_i\to x$ in $\bts{\ma}{}{\mb}$, and let 
$l=\sum_if_i\oslash g_i$ be such that $\norm{l(t,s)-x}<\epsilon$, 
where $x\in\bts{A_t}{}{B_s}$, $x_i\in\bts{A_{t_i}}{}{B_{s_i}}$. Since 
$x_i\to x$ and $l(t_i,s_i)\to l(t,s)$,  
there exists $i_0$ such that   
$\forall i\geq i_0$ we have $\norm{l(t_i,s_i)-x_i}<\epsilon$. Now 
$\norm{L_u l(t_i,s_i)-L_ux_i}\leq\epsilon\norm{u}$, and  
$\norm{L_u l(t,s)-L_ux}\leq\epsilon\norm{u}$, and since 
$L_ul(t_i,s_i)\to L_u(t,s)$, we conclude that $L_ux_i\to L_ux$.  
\par Let see now that $\iota_{\ma}$ is strictly continuous on cylinders.
If $a\in\ma$, $b\in\mb$, and $(u_i)\subseteq\ma$ is a net strictly
convergent to $u\in\ma$, with $\norm{u_i},\norm{u}\leq C$, then: 
$\iota_{\ma}(u_i)(a\otimes b)=u_ia\otimes b\to ua\otimes b=
                                              \iota_{\ma}(u)(a\otimes
                                              b)$ 
and $(a\otimes b)\iota_{\ma}(u_i)=au_i\otimes b\to au\otimes b=
                                              (a\otimes
                                              b)\iota_{\ma}(u)$. 
Then $\iota_{\ma}(u_i)x\to\iota_{\ma}(u)x$ and $x\iota_{\ma}(u_i)\to
x\iota_{\ma}(u)$, $\forall x\in\bats{\ma}{\mb}$. Since  
$\norm{\iota_{\ma}(u_i)},\norm{\iota_{\ma}(u)}\leq C$, we conclude that  
$\iota_{\ma}(u_i)x\to\iota_{\ma}(u)x$ and $x\iota_{\ma}(u_i)\to
x\iota_{\ma}(u)$, $\forall x\in\bts{\ma}{}{\mb}$. Thus 
$\iota_{\ma}(u_i)$ converges strictly to $\iota_{\ma}(u)$.  
\par Similarly we construct $\iota_{\mb}:M(\mb )\to M(\bts{\ma}{}{\mb})$: 
if $v\in M(\mb )$ and $a\in\ma$, $b\in\mb$, then
$\iota_{\mb}(v)(a\otimes b)=a\otimes vb$, and $(a\otimes
b)\iota_{\mb}(v)=a\otimes bv$. It is clear that 
$\iota_{\ma}(u)\iota_{\mb}(v)=\iota_{\mb}(v)\iota_{\ma}(u)$, $\forall u
\in M(\ma )$, $v\in M(\mb )$, and also that $\iota_{\ma}(a)\iota_{\mb}(b)=
a\otimes b$, $\forall a\in \ma$, $b\in\mb$.
\par This way we obtain a map $M(\ma )\times M(\mb )\to 
M(\bts{\ma}{}{\mb})$, given by
$(u,v)\mapsto\iota_{\ma}(u)\iota_{\mb}(v)$, which is 
bilinear on each $M_t(\ma )\times M_s(\mb )$, and therefore we get 
a map $\bats{M(\ma )}{M(\mb )}\to M(\bts{\ma}{}{\mb})$, 
which is linear on each $M_t(\ma )\otimes M_s(\mb )$, and which is
a \hm of \fbs because $\iota_{\ma}(u)$ and $\iota_{\mb}(v)$ commute,
$\forall u\in M(\ma )$, $v\in M(\mb )$. 
\end{proof}

\begin{prop}\label{prop:maxreps}
Let \ma and \mb be \fbs over the locally compact groups $G$ and $H$
respectively, and let $\mh$ be a Hilbert module. Then for each  
$(\pi_1,\pi_2)\in R(\ma,\mb,\mh)$, there exists a unique 
$\pi\in R(\tmax{\ma}{\mb},\mh )$ such that 
$\pi (a\otimes b)=\pi_1(a)\pi_2(b)$, $\forall a\in \ma$, $b\in \mb$,
and the map $(\pi_1,\pi_2)\mapsto \pi$ thus defined is a bijection
between $R(\ma,\mb,\mh)$ and $R(\tmax{\ma}{\mb},\mh)$.  
\end{prop}
\begin{proof}
Let $(\pi_1,\pi_2)\in R(\ma ,\mb ,\mh )$. The map $\ma\times 
\mb\to\adj{\mh}$ such that $(a_t,b_s)\mapsto\pi_1(a_t)\pi_2(b_s)$ is 
bilinear on each $A_t\times B_s$, and therefore there exists a unique
$\pi: \bats{\ma}{\mb}\to\adj{\mh}$ such that 
$\pi(a_t\otimes b_s)=\pi_1(a_t)\pi_2(b_s)$. Since $\pi_1(a)$ and
$\pi_2(b)$ commute, $\forall a\in\ma$, $b\in\mb$, we have that 
$\pi:\bats{\ma}{\mb}\to 
\adj{\mh}$ is a \rep of the pre-\fb $\bats{\ma}{\mb}$. Note that 
$x\mapsto\norm{x}':=\max\{\norm{x}_{\max},\norm{\pi
(x)}\}$ is a $C^*$-norm on $\bats{A_t}{B_s}$, so $\norm{\cdot}'= 
\norm{\cdot}_{\max}$, and then $\norm{\pi (x)}\leq 
\norm{x}_{\max}$, $\forall x\in \bats{A_t}{B_s}$. Thus $\pi$
has a unique extension to a \rep
$\pi:(\tmax{\ma}{\mb})_d\to\adj{\mh}$. It is easy to see that this 
\rep is non-degenerate: since $\pi_1$ is non-degenerate, for every
$h\in\mh$ there exist $h'\in\mh$ and $a\in A_e$ such that 
$\pi_2(a)h'=h$. Since $\pi_2$ is non degenerate, there exist
$h''\in\mh$ and $b\in B_e$ such that $\pi_1(b)h''=h'$. Therefore $\pi
(a\otimes b)h''=\pi_1(a)\pi_2(b)h''=\pi_1(a)h'=h$. 
\par To see that $\pi$ is continuous is sufficient, by
\cite[II-13.16]{fd}, to prove that $\forall f\in C_c(\ma )$,  
$ g\in C_c(\mb )$, and $h\in\mh$, the map $F:G\times H\to\mh$ such
that $F(t,s)=\pi\left((f\oslash g)(t,s)\right)h$ is continuous. Now 
\begin{align*}
\begin{split}
\norm{F(t,s)-F(t_0,s_0)}
&=\norm{\pi_1\left(f(t)\right)\pi_2\left( g(s)\right)h-
        \pi_1\left(f(t_0)\right)\pi_2\left( g(s_0)\right)h}\\
&\leq\norm{\pi_1(f (t))\left(\pi_2( g(s))-\pi_2( g(s_0))\right)h}\\
 & \quad +\norm{\left(\pi_1(f(t))-\pi_1(f(t_0))\right)\pi_2( g(s_0))h}\\
&\leq\norm{f}_{\infty}\, \norm{(\pi_2( g(s))h-\pi_2( g(s_0))h}\\
 & \quad + \norm{\pi_1(f(t))\pi_2( g(s_0))h-\pi_1(f(t_0))
   \pi_2((s_0))h}\\
\end{split}
\end{align*} 
Then $F(t,s)\to F(t_0,s_0)$ if $(t,s)\to
(t_0,s_0)$. Therefore $\pi\in R(\tmax{\ma}{\mb},\mh )$. 
\par Conversely suppose that $\pi\in R(\tmax{\ma}{\mb},\mh)$. By 
\cite[VIII-15.3]{fd}, $\pi$ can be uniquely extended to a \rep $\pi'$
of $M(\tmax{\ma}{\mb})$ such that $x\mapsto\pi'(x)h$ is strictly
continuous on cylinders, $\forall h\in\mh$. Let  
$\pi_1=\pi'\iota_{\ma}\r{\ma}:\ma\to\adj{\mh}$, and  
$\pi_2=\pi'\iota_{\mb}\r{\mb}:\mb\to\adj{\mh}$, where $\iota_{\ma}$
and $\iota_{\mb}$ are the inclusions provided by Lemma
\ref{lem:incmult}. Since $\pi'$, $\iota_{\ma}$ and $\iota_{\mb}$ are
continuous on cylinders, immediately follows that $\forall h\in\mh$ the
maps $\ma\to\mh$ and $\mb\to\mh$ given by $a\mapsto\pi_1(a)h$ and
$b\mapsto\pi_2(b)h$ respectively are continuous, from where it follows
that $\pi_1$ and $\pi_2$ are continuous \reps. Since $\iota_{\ma}(a)$
and $\iota_{\mb}(b)$ commute, $\forall a\in \ma$ and $b\in\mb$, then 
$\pi_1(a)$ and $\pi_2(b)$ commute as well. Finally, the \reps $\pi_1$  
and $\pi_2$ are non-degenerate. To see this is enough to show that
$\pi (a\otimes b)h$ is in the image of 
$\pi_1$ and the image of $\pi_2$, $\forall a\in\ma$, $b\in\mb$ and
$h\in\mh$. By Cohen-Hewitt $a$ and $b$ can be factorized as
$a=a_1a_2$, $b=b_1b_2$, and therefore $\pi_1(a_1)\pi (a_2\otimes
b)h=\pi (a_1a_2\otimes b)h=\pi(a\otimes b)h$ and 
$\pi_2(b_1)\pi(a\otimes b_2)h=\pi(a\otimes b_1b_2)h=\pi(a\otimes b)h$.
\par In conclusion we constructed two correspondences
$(\pi_1,\pi_2)\mapsto \pi$ and 
$\pi\mapsto (\pi\iota_{\ma}\r{\ma},\pi\iota_{\mb}\r{\mb})$, which
clearly are mutually inverse. 
\end{proof}

\section{C*-algebras of Tensor Products of Fell
  Bundles}\label{sec:cross} 
\par The first goal of this section is to compare tensor products of the
cross-sectional algebras of the \fbs \ma and \mb with the
cross-sectional algebras of tensor products of \ma and \mb. This is
accomplished in Propositions~\ref{prop:max} and \ref{prop:comp}, and
in Theorem \ref{thm:main}. The second objective is to give some
applications.   
\par Let \mb be a \fb over a locally compact group $G$. Then there are
two important cross-sectional \css associated with \mb: the full cross-sectional
algebra $C^*(\mb)$, and the reduced cross-sectional algebra
$C^*_r(\mb)$. We recall next their definitions. 
\par Suppose that $G$ is a locally compact group with Haar measure
$\la$ and modular function $\Del$.  Let \mb be a \fb over $G$ and let 
$L^1(\mb):=\{f:G\to\ma:\ f(t)\in B_t,\forall t\in G, \text{ and }
(t\mapsto\norm{f(t)})\in L^1(G,\la)\}$. Then $C_c(\mb)$ and $L^1(\mb)$
are *-algebras with the operations: $f*g(t)=\int_{G}f(r)g(r^{-1}t)$,  
$f^*(t)=\Delta(t)^{-1}f(t^{-1})^*$. Moreover, $L^1(\mb)$ is a Banach
*-algebra with the norm: 
$\norm{f}_1=\int_{G}\norm{f(t)}$. The enveloping \cs of
$L^1(\mb)$ is called the cross-sectional algebra of $\mb$, and it is 
denoted by $C^*(\mb)$.
\par Suppose that $\phi :\ma\to\mb$ is a \hm of \fbs.
If $f\in L^1(\ma)$ we have that $\phi^1(f):G\to\mb$, given by 
$\phi^1(f)(t)=\phi\big(f(t)\big)$, belongs to $L^1(\mb)$, and  
$\norm{\phi^1(f)}_1\leq\norm{f}_1$. Moreover $\phi^1$ is a \hm of
*-algebras, so it uniquely extends to a \hm $C^*(\phi ):C^*(\ma)\to
C^*(\mb)$. This way we obtain a functor from the category of \fbs over
$G$ to the category of \css. In fact this functor is the compostion of
the functor $A\mapsto C^*(A)$ from the category of Banach *-algebras
with approximate unit and contractive \hms to the category of \css,
with the functor: $\ma\mapsto L^1(\ma)$,
$\big(\ma\stackrel{\phi}{\to}\mb\big)\mapsto 
\big(L^1(\ma)\stackrel{\phi^1}{\to}L^1(\mb)\big)$. 
\par There is a bijection between non-degenerate \reps of
the \fb \mb and non-degenerate \reps of the \cs $C^*(\mb)$. In one
direction this correspondence consists of passing from a \rep
$\pi:\mb\to\adj{\mh}$ to its integrated form
$\int_G\pi:C^*(\mb)\to\adj{\mh}$, characterized by
$\pr{\int_G\pi(f)\xi}{\eta}=\int_G\pr{\pi(f(t))\xi}{\eta}dt$, $\forall
f\in C_c(\mb)$, $\xi,\eta\in\mh$ (see \cite[VIII-13.2]{fd}).
\par Among the \reps of \mb there is one of particular importance: the
\textit{(left) regular \rep}, which we describe below. Note that
$C_c(\mb)$ is a right $B_e$-module with the action given by pointwise
multiplication. Moreover the map $\pr{\cdot}{\cdot}:C_c(\mb)\times
C_c(\mb)\to B_e$ such that $\pr{\xi}{\eta}=\int_G\xi(s)^*\eta(s)ds$ is
a pre-inner product. Completing $C_c(\mb)$ with respect to the norm
defined by $\pr{\cdot}{\cdot}$ we obtain a full right Hilbert
$B_e$-module, which is denoted by $L^2(\mb)$. Again, it is not dificult to check that $\mb\mapsto L^2(\mb)$ is a functor. There exists a unique
\rep $\Lambda^\mb:\mb\to\adj{L^2(\mb)}$ such that $\La^\mb_{b_t}\xi
(s)=b_t\xi(t^{-1}s)$ $\forall s,t\in G$, $b_t\in B_t$ and $\xi\in
C_c(\mb)$ (if no confusion can arise we write just $\Lambda$ instead of $\Lambda^\mb$). This is called the regular \rep of \mb on $L^2(\mb)$. Its 
integrated form is also called the regular \rep, and it satisfies
$\La_f(\xi)=f*\xi$, $\forall f\in C_c(\mb)\subseteq C^*(\mb)$ and
$\forall \xi\in C_c(\mb)\subseteq L^2(\mb)$, where the
\textit{convolution} $f*\xi$ is defined as:
$f*\xi(t)=\int_Gf(s)\xi(s^{-1}t)ds$. The reduced cross-sectional
algebra is then defined as:
$C^*_r(\mb):=\La^\mb(C^*(\mb))\subseteq\adj{L^2(\mb)}$. 
\par When we look at the regular \rep as a \hm $\La^\mb:C^*(\mb)\to
C^*_r(\mb)$, then it is clear that $\La^\mb$ is onto. In the case that
$\La^\mb$ is also injective, thus an isomorphism, we say that the \fb \mb
is \textit{amenable}. The reader is referred to \cite{examen},
\cite{exeng} and \cite{env} for further information on the reduced
cross-sectional algebra.     
\par It can be shown that $\mb\mapsto C^*_r(\mb)$ also is a functor, and in fact $\Lambda$ is a natural transformation from $C^*$ to $C^*_r$ (\cite[page 277]{fellequiv}).  
\subsection{Cross-sectional algebras}\label{subsec:comparison}  
\begin{lem}\label{lem:bolu2}
Let \fela and \felbh be \fbs and 
suppose that $\bts{\ma}{}{\mb}$ is a tensor product of $\ma$ and
$\mb$. Then there exists a unique \hm of algebras 
$j_c:\bats{C_c(\ma)}{C_c(\mb)}\to C_c(\bts{\ma}{}{\mb})$, such that
$j_c(f\odot g)=f\oslash 
g$, that is: $j_c(f\odot g)(t,s)=f(t)\otimes g(s)$, 
$\forall f\in C_c(\ma)$, $g\in C_c(\mb )$, $t\in G$, and $s\in
H$. Moreover $j_c$ is injective and $j_c(\bats{C_c(\ma)}{C_c(\mb)})$ is
dense in $C_c(\bts{\ma}{}{\mb})$ in the \ilt. 
\end{lem}
\begin{proof}
The existence and uniqueness of the linear map $j_c$ follows from the
universal property of tensor products. It is clear that $j_c$ is a \hm
of *-algebras. To see that it is injective suppose that 
$l=\sum_{i=1}^nf_i\odot g_i\in\ker j_c$.  
Then $0=\pr{j_c(l )}{j_c(l )}=\int_{G\times
H}\sum_{i,j=1}^nf_i(t)^*f_j(t)\otimes g_i(s)^*g_j(s) d(t,s)$. 
On the other hand we have $\int_{G\times
H}\sum_{i,j=1}^nf_i(t)^*f_j(t)\otimes g_i(s)^* g_j(s)
d(t,s)=\big(\int_G\sum_{i,j=1}^nf_i(t)^*f_j(t)dt\big)\otimes 
\big(\int_H\sum_{i,j=1}^ng_i(s)^* g_j(s)ds\big)$. Therefore, if we
think of $l$ as an element of $\bats{L^2(\ma)}{L^2(\mb)}$ we have that   
$\pr{j_c(l)}{j_c(l)}=\pr{l}{l}$, where the latter is the pre-inner
product of $\bats{L^2(\ma)}{L^2(\mb)}$ computed in $l$. Since
$\pr{l}{l}=0$, it follows that $l=0$. 
\par Let see that $j_c(\bats{C_c(\ma)}{C_c(\mb)})$ is dense in
$C_c(\bts{\ma}{}{\mb})$ in the \ilt. It is clear that
$j_c(\bats{C_c(\ma)}{C_c(\mb)}))(t,s)$ is dense in   
$(\bts{\ma}{}{\mb})_{(t,s)}$, $\forall (t,s)\in G\times H$. On the
other hand, if $\Theta =\bats{C_c(G)}{C_c(H)}$, let $\theta\in \Theta$
and $l\in \bats{C_c(\ma)}{C_c(\mb)}$, say $\theta
=\sum_i\phi_i\odot\psi_i$ and $l=\sum_jf_j\odot g_j$, then:  
$\theta j_c(l)(t,s)=\big(\sum_i\phi_i(t)\psi_i(s)\big)
\big(\sum_jf_j(t)\otimes g_j(s)\big)
=\sum_i\sum_j(\phi_if_j)(t)\otimes (\psi_i g_j)(s)=j_c(l')(t,s),$ 
where $l'=\sum_i\sum_j\phi_if_j\odot
\psi_ig_j\in\bats{C_c(\ma)}{C_c(\mb)}$.  
Thus $j_c(\bats{C_c(\ma)}{C_c(\mb)})$ is dense in
$C_c(\bts{\ma}{}{\mb})$ by \cite[Lemma 5.1]{env}.  
\end{proof}

\begin{prop}\label{prop:max}
Let \fela and \felbh be \fbs. Then there exists a unique isomorphism  
$j:\tmax{C^*(\ma)}{C^*(\mb)}\to C^*(\tmax{\ma}{\mb})$, such that  
$j(f\otimes g)(t,s)=f(t)\otimes g(s)$, $\forall f\in C_c(\ma)$, 
$g\in C_c(\mb)$, and $(t,s)\in G\times H$.
\end{prop}
\begin{proof}
\par Recall that if $\mh$ is a Hilbert space and $\mc=(C_t)_{t\in G}$
is a \fb, then there is a bijection between $R(\mc ,\mh)$ and  
$R\big(C^*(\mc ),\mh\big)$ such that to each $\pi\in R(\mc ,\mh )$
corresponds the integrated \rep $\int_G\pi$ of  
$C^*(\mc)$, which is determined by its values on elements of
$C_c(\mc)$: if $f\in C_c(\mc)$  
and $h\in\mh$, then
$(\int_G\pi)f\r{h}=\int_G\pi\big(f(t)\big)hdt$. Note as well that if
$\mc'=(C'_s)_{s\in 
H}$ is another \fb, then the map 
$R(\mc ,\mc',\mh)\to R\big(C^*(\mc),C^*(\mc'),\mh\big)$ such that
$(\pi,\pi')\mapsto (\int_G\pi ,\int_H\pi')$ is also a bijection,
because the corresponding integrands commute. 
On the other hand, by Proposition \ref{prop:maxreps} we have a
bijection between $R(\mc,\mc',\mh)$ and $R(\tmax{\mc}{\mc'},\mh)$, 
given by $(\pi_1,\pi_2)\mapsto\pi_1\times\pi_2$, where 
$(\pi_1\times\pi_2)(a\otimes b)=\pi_1(a)\pi_2(b)$. 
\par Let $j_c:\bats{C_c(\ma)}{C_c(\mb)}\to C_c(\tmax{\ma}{\mb})$ be the
map provided by Lemma \ref{lem:bolu2}. The comments above imply that 
$\tmax{C^*(\ma)}{C^*(\mb)}$ and $C^*(\tmax{\ma}{\mb})$ are
respectively the completions of \bats{C_c(\ma)}{C_c(\mb)} and 
$j_c(\bats{C_c(\ma)}{C_c(\mb)})$ with respect to the norms:  
\[\norm{\sum_if_i\odot g_i}
=\sup\{\norm{\sum_i\int_G\pi_1(f_i)\int_H\pi_2( g_i)}:\, (\pi_1,\pi_2)\in 
     R(\ma ,\mb ,\mh)\},\]
\[\norm{j_c(\sum_i\! f_i\odot g_i)}
=\sup\{\norm{\int_{G\times
H}\hspace*{-.3cm}(\pi_1\times\pi_2)(\sum_if_i\otimes g_i)} 
:(\pi_1,\pi_2)\in\! R(\ma,\mb,\mh)\}.\]
Now, if $h\in\mh$:
\begin{align*}
\left(\int_{G\times H}\big(\pi_1\times\pi_2\big)
                      \big(\sum_if_i\otimes g_i\big)\right)h
&=\int_{G\times H}\sum_i\pi_1\big(f_i(t)\big)\pi_2\big( g_i(s)\big)hd(t,s)\\
&=\sum_i\int_G\int_H\pi_1\big(f_i(t)\big)\pi_2\big( g_i(s)\big)hdsdt\\
&=\sum_i\int_G\pi_1\big(f_i(t)\big)\int_H\pi_2\big( g_i(s)\big)hdsdt\\
&=\sum_i\int_G\pi_1\big(f_i(t)\big)\left(\int_H\pi_2( g_i)\right)hdt\\
&=\sum_i\left(\int_G\pi_1(f_i)\int_H\pi_2( g_i)\right)h. 
\end{align*} 
Thus $j_c:\bats{C_c(\ma)}{C_c(\mb)}\to j_c(\bats{C_c(\ma)}{C_c(\mb)})$
is an isometry with these norms so it extends uniquely to an
isomorphism between $\tmax{C^*(\ma)}{C^*(\mb)}$ and the C*-algebra 
$\ov{j_c(\bats{C_c(\ma)}{C_c(\mb)}}$. Since by Lemma \ref{lem:bolu2} 
$j_c(\bats{C_c(\ma)}{C_c(\mb)}$ is dense in $C_c(\tmax{\ma}{\mb})$ in
the \ilt, then it is also dense in $C^*(\tmax{\ma}{\mb})$.
\end{proof}

\begin{prop}\label{prop:comp}
Let \fela and \felbh be \fbs, and suppose that $\al\geq\beta$ are 
$C^*$-norms on $\bats{\ma}{\mb}$. Then there exist a unique \hm   
$\sigma^{\al}_{\beta}:C^*(\bts{\ma}{\al}{\mb})\to  
C^*(\bts{\ma}{\beta}{\mb})$ such that $\sigma^{\al}_{\beta}(f\oslash g)=f\oslash g$, $\forall  
f\in C_c(\ma)$, $g\in C_c(\mb)$. Moreover $\sigma^{\al}_{\beta}$ is surjective.
\end{prop}
\begin{proof}
By Proposition \ref{prop:extmm2} there exists a surjective \hm of
\fbs $\sigma^{\al}_{\beta}:\bts{\ma}{\al}{\mb}\to    
\bts{\ma}{\beta}{\mb}$. Then there is an induced \hm $\sigma^{\al}_{\beta}: C^*(\bts{\ma}{\al}{\mb})\to C^*(\bts{\ma}{\beta}{\mb})$, which we still denote by $\sigma^{\al}_{\beta}$. If $f\in C_c(\ma)$, $g\in C_c(\mb)$, $\sigma^{\al}_{\beta}(f\oslash g)(t,s)= \sigma^{\al}_{\beta}\big(f(t)\otimes g(s)\big)=f(t)\otimes g(s)=f\oslash g (t,s)$, from where it follows that $\sigma^{\al}_{\beta}(f\oslash
g)=f\oslash g$. Since $\gen\{f\oslash g:\, f\in C_c(\ma),\, g\in
C_c(\mb)\}$ is dense in $C^*(\bts{\ma}{\al}{\mb})$, we conclude that
$\sigma^{\al}_{\beta}$ is surjective.    
\end{proof}

\par Consider two \fbs \fela and \felbh. Then $L^2(\ma)$ and
$L^2(\mb)$ are full right Hilbert modules over $A_e$ and $B_e$
respectively. If $\al$ is a $C^*$-norm on $\bats{\ma}{\mb}$, then $\al|_{(\bats{A_e}{B_e})}\in\mn(\bats{A_e}{B_e})$. Since $L^2(\ma)^r=A_e$ and $L^2(\mb)^r=B_e$, $\al|_{(\bats{A_e}{B_e})}$ defines a C*-norm $\tilde{\alpha}$ on $\bats{L^2(\ma)}{L^2(\mb)}$, given by \eqref{eqn:tildeh}, that is $\tilde{\alpha}(\mu):=\sqrt{\alpha(\pr{\mu}{\mu})}$, $\forall \mu\in \bats{L^2(\ma)}{L^2(\mb)}$. The completion $\bts{L^2(\ma)}{\tilde{\alpha}}{L^2(\mb)}$ of $\bats{L^2(\ma)}{L^2(\mb)}$ with respect to $\tilde{\alpha}$ is a full right Hilbert module over $\bts{A_e}{\alpha|_{A_e\otimes B_e}}{B_e}$, so we have $(\bts{L^2(\ma)}{\tilde{\alpha}}{L^2(\mb)})^r=\bts{A_e}{\alpha|_{A_e\otimes B_e}}{B_e}$, whose its corresponding inner
product is determined by $\pr{\xi_1\otimes\eta_1}{\xi_2\otimes\eta_2}
=\pr{\xi_1}{\xi_2}\otimes\pr{\eta_1}{\eta_2}$,  
$\forall \xi_1, \xi_2\in L^2(\ma)$, $\eta_1, \eta_2\in L^2(\mb)$  (see Theorem~\ref{thm:correspondence}).   
\begin{lem}\label{lem:l2}
Let \fela and \felbh be \fbs, and let $\al$ be a $C^*$-norm on 
$\bats{\ma}{\mb}$. If $\tilde{\alpha}$ is as above, then there exists a unique isomorphism   
$j_2:\bts{L^2(\ma)}{\tilde{\al}}{L^2(\mb)}\to  L^2(\bts{\ma}{\al}{\mb}),$ 
such that $j_2(\xi\otimes\eta)=\xi\oslash\eta$, 
$\forall\xi\in C_c(\ma)\subseteq L^2(\ma)$, $\eta\in C_c(\mb)\subseteq 
L^2(\mb)$. In particular we have
$\tmin{L^2(\ma)}{L^2(\mb)}\cong L^2(\tmin{\ma}{\mb})$ and $\tmax{L^2(\ma)}{L^2(\mb)}\cong L^2(\tmax{\ma}{\mb})$.
\end{lem}
\begin{proof}
Let $j_c$ be the map defined in Lemma \ref{lem:bolu2}.
If $\xi_1,\xi_2\in C_c(\ma)$,  $\eta_1,\eta_2\in C_c(\mb)$, then  
$j_c(\xi_1\otimes\eta_1)$, $j_c(\xi_2\otimes\eta_2)\in 
C_c(\bts{\ma}{\al}{\mb})\subseteq L^2(\bts{\ma}{\al}{\mb})$ and we
have  
\begin{gather*}
  \pr{j_c(\xi_1\otimes\eta_1)}{j_c(\xi_2\otimes\eta_2)} 
=\int_{G\times H}(\xi_1\otimes\eta_1)(t,s)^*(\xi_2\otimes\eta_2)(t,s)d(t,s)\\ 
=\int_G\int_H\xi_1(t)^*\xi_2(t)\otimes\eta_1(s)^*\eta_2(s)dsdt
=\pr{\xi_1}{\xi_2}\otimes\pr{\eta_1}{\eta_2}
\end{gather*}
On the other hand, if $a\in A_e$, $b\in B_e$, $\xi\in C_c(\ma)$, and 
$\eta\in C_c(\mb )$, we have 
$\big(j_c(\xi\otimes\eta)\big)(a\otimes b)(t,s)
=(\xi(t)\otimes\eta(s))(a\otimes b)=\xi(t)a\otimes\eta(s)b
=j_c\big((\xi\otimes\eta)(a\otimes b)\big)(t,s).$
Thus $j_c$ is a \hm of pre-Hilbert modules over 
$\bts{A_e}{\al}{B_e}$ which is injective by Lemma \ref{lem:bolu2} and 
has dense image in $L^2(\bts{\ma}{\al}{\mb})$: by Lemma
\ref{lem:bolu2}, the image of $j_c$ is dense in the
$C_c(\bts{\ma}{\al}{\mb})$ in the \ilt, and therefore is dense in 
$L^2(\bts{\ma}{\al}{\mb})$. Thus $j_c$ extends uniquely to 
an isomorphism $j_2:\bts{L^2(\ma)}{\tilde{\al}}{L^2(\mb)}\to
L^2(\bts{\ma}{\al}{\mb})$. The last two statements follow from the fact that if $\alpha=\norm{\ }_{\textrm{min}}$, then also  $\tilde{\alpha}=\norm{\ }_{\textrm{min}}$ and, similarly, if $\alpha=\norm{\ }_{\textrm{max}}$, then also  $\tilde{\alpha}=\norm{\ }_{\textrm{max}}$ (because $\alpha\mapsto \alpha_{\bats{A_e}{B_e}}\mapsto\tilde{\alpha}$ are isomorphisms of posets).  
\end{proof}

\par With the notation as above, we have inclusions
$C^*_r(\ma)\subseteq\mathcal{L}\big(L^2(\ma)\big)$  
and $C^*_r(\mb )\subseteq\mathcal{L}\big(L^2(\mb)\big)$, so  
$\bats{C^*_r(\ma )}{C^*_r(\mb )}$ is included in 
$\bats{\mathcal{L}\big(L^2(\ma)\big)}{\mathcal{L}\big(L^2(\mb)\big)}$, which in turn is included in $\mathcal{L}\big(\bts{L^2(\ma)}{\tilde{\al}}{L^2(\mb)}\big)$according to Corollary~\ref{cor:alfita}. Therefore we
have an inclusion $\bats{C^*_r(\ma )}{C^*_r(\mb )}\inc
\mathcal{L}\big(\bts{L^2(\ma)}{\tilde{\al}}{L^2(\mb)}\big).$  

\begin{df}\label{df:chal}
If $\al$ is a $C^*$-norm on $\bats{\ma}{\mb}$, we define  
\bts{C^*_r(\ma)}{\chal}{C^*_r(\mb)} to be the closure of 
$\bats{C^*_r(\ma)}{C^*_r(\mb)}$ in 
$\mathcal{L}\big(\bts{L^2(\ma)}{\tilde{\al}}{L^2(\mb)}\big)$ (that is: we call $\chal$ the norm on $\bats{C^*_r(\ma)}{C^*_r(\mb)}$ inherited by the inclusion above). Recall that, in particular, if
$\alpha=\norm{\ }_{\textrm{min}}$, then we also have $\chal=\norm{\ }_{\textrm{min}}$.
\end{df}

\par Suppose that $u:\mh_1\to\mh_2$ is a unitary operator between the
Hilbert modules $\mh_1$ and $\mh_2$. Then $u$ induces an isomorphism
$Ad_u:\adj{\mh_1}\to\adj{\mh_2}$, given by $Ad_u(T)=uTu^*$, $\forall
T\in\adj{\mh_1}$.  

\begin{prop}\label{prop:alpha}
Let \ma and \mb be \fbs over the locally compact groups $G$ and $H$
respectively, $j_2:\bts{L^2(\ma)}{\tilde{\al}}{L^2(\mb)}\to
L^2(\bts{\ma}{\al}{\mb})$ the isomorphism given by Lemma
\ref{lem:l2}, and $\chal$ the C*-norm given by Definition~\ref{df:chal}.  
Then $Ad_{j_2}\big(\bts{C^*_r(\ma)}{\chal}{C^*_r(\mb)}\big)= 
C^*_r(\bts{\ma}{\al}{\mb})$, and there is a unique isomorphism  
$j_r:\bts{C^*_r(\ma)}{\chal}{C^*_r(\mb)}\to
C^*_r(\bts{\ma}{\al}{\mb})$ such that 
$j_r\big(\La^{\ma}_f\otimes\La^{\mb}_g\big)=
\La^{\bts{\ma}{\al}{\mb}}_{(f\otimes g)}$,
$\forall f\in C_c(\ma)$, $g\in C_c(\mb)$. 
\par\noindent In particular $C^*_r(\tmin{\ma}{\mb})\cong\tmin{C^*_r(\ma)}{C^*_r(\mb)}$. 
\end{prop}
\begin{proof}
As usual, by the universal property of tensor products we see that
there exists a unique map $\La^{\ma}\odot\La^{\mb}:\bats{\ma}{\mb}\to 
\mathcal{L}\big(\bts{L^2(\ma)}{\al}{L^2(\mb)}\big)$ such that 
$(\La^{\ma}\odot\La^{\mb})(a\odot b)=\La^{\ma}_a\otimes\La^{\mb}_b$,  
$\forall a\in\ma$, $b\in\mb$. Writing just $\La$ instead of
$\La^{\bts{\ma}{\al}{\mb}}$ we have  
\begin{align*}
\La_{(a_t\otimes b_s)}\big(j_2(\xi\otimes\eta)\big)(t_0,s_0)
&=(a_t\otimes b_s)(\xi\otimes\eta)\big((t,s)^{-1}(t_0,s_0)\big)\\
&=a_t\xi(t^{-1}t_0)\otimes b_s\eta(s^{-1}s_0)\\
&=\big(\La^{\ma}_{a_t}\xi\big)(t_0)\otimes\big(\La^{\mb}_{b_s}\eta\big)(s_0)\\ 
&=\big(\La^{\ma}_{a_t}\xi\otimes\La^{\mb}_{b_s}\eta\big)(t_0,s_0)\\
&=j_2\big(\La^{\ma}_{a_t}\xi\otimes\La^{\mb}_{b_s}\eta\big){(t_0,s_0)}\\
&=j_2\big((\La^{\ma}\odot\La^{\mb})(a_t\odot
b_s)(\xi\otimes\eta)\big)(t_0,s_0),
\end{align*} 
It follows that $\La_x=j_2(\La^{\ma}\odot\La^{\mb})(x)j_2^*$, $\forall 
x\in\bats{\ma}{\mb}$, so $\La^{\ma}\odot\La^{\mb}$
extends uniquely to a \rep 
$\La^{\ma}\otimes\La^{\mb}:\bts{\ma}{\al}{\mb}\to
\adj{\bts{L^2(\ma)}{\al}{L^2(\mb)}}$ such that 
$\La^{\ma}\otimes\La^{\mb}(x)=j_2^*\La_xj_2$, $\forall x\in
\bts{\ma}{\al}{\mb}$. Taking the corresponding integrated \reps, we  
have that $\La_{(f\otimes g)}=j_2(\La^{\ma}\otimes\La^{\mb})(f\otimes
g)j_2^*$, $\forall f\in C_c(\ma)$, $g\in C_c(\mb)$. If $j_c$ is the
map given by Lemma \ref{lem:bolu2}, then 
$j_c(\bats{C_c(\ma)}{C_c(\mb)})$
is dense in $C^*_r(\bts{\ma}{\al}{\mb})$. Therefore we
conclude that $C^*_r(\bts{\ma}{\al}{\mb})=
j_2\big(\bts{C^*_r(\ma)}{\chal}{C^*_r(\mb)}\big)j_2^*$, 
as we wanted to prove. 
\par In particular,
$j_r:\bts{C^*_r(\ma)}{\chal}{C^*_r(\mb)}\to 
C^*_r(\bts{\ma}{\al}{\mb})$ given by $x\mapsto j_2xj_2^*$ is an
isomorphism satisfying $j_r\big(\La^{\ma}_f\otimes\La^{\mb}_g\big)=
\La_{(f\otimes g)}$, $\forall f\in C_c(\ma)$, $g\in C_c(\mb)$. The
uniqueness of such an isomorphism is clear. As for the last statement 
just recall that $\chal={\norm{\ }}_{\min}$ if $\alpha=\norm{\ }_{\min}$ (Corollary~\ref{cor:alfita}).  
\end{proof}
\par In functorial language, Proposition~\ref{prop:max} and the last statement of Proposition~\ref{prop:alpha} can be stated as follows. Let $\mb\mapsto C^*(\mb)$ and $\mb\mapsto C^*_r(\mb)$ be the functors sending a Fell bundle $\mb$ to its cross-sectional and reduced croseed sectional algebras respectively, and let $\otimes_{\textrm{max}}$ and $\otimes_{\textrm{min}}$ be the bifunctors of taking maximal and minimal tensor products respectively (of Fell bundles or of C*-algebras). Then we have:
\begin{gather*}
  C^*\circ\otimes_{\textrm{max}}=\otimes_{\textrm{max}}\circ (C^*\times C^*):\mathsf{F}\times \mathsf{F}\to \mathsf{C}^*\\
C^*_r\circ\otimes_{\textrm{min}}=\otimes_{\textrm{min}}\circ (C^*_r\times C^*_r):\mathsf{F}\times \mathsf{F}\to \mathsf{C}^*,
\end{gather*}
where $\mathsf{F}$ is the category of Fell bundles and $\mathsf{C}^*$ the category of C*-algebras. That is: taking full (reduced) cross-sectional algebras commute with taking maximal (respectively: minimal) tensor products.  
 
\begin{thm}\label{thm:main}
Let \ma and \mb be \fbs over the locally compact groups $G$ and $H$
respectively. Then for every $C^*$-norm $\al$ on \bats{\ma}{\mb} we
have the following commutative diagram $D_\alpha$:  
\[ 
\xymatrix
{C^*(\tmax{\ma}{\mb})\ar@{->>}[r]^-{\sigma^{\max}_{\al}} 
&C^*(\bts{\ma}{\al}{\mb})\ar@{->>}[r]^-{\La}
&C^*_r(\bts{\ma}{\al}{\mb})\\
\tmax{C^*(\ma)}{C^*(\mb)}
\ar@{->>}[r]_-{\La^{\ma}\otimes_{\max}\La^{\mb}}
 \ar[u]^-{j\cong}
&\tmax{C^*_r(\ma)}{C^*_r(\mb)}\ar@{->>}[r]_-{\tilde{\sigma}^{\max}
_{\chal}}
&\bts{C^*_r(\ma)}{\chal}{C^*_r(\mb )}\ar[u]_-{\cong j_r} }\]
where $\La=\La^{\bts{\ma}{\al}{\mb}}$, the map $\sigma^{\max}_{\al}$
is provided by Proposition~\ref{prop:comp}, $j$ is given by Proposition
\ref{prop:max}, $j_r$ by Proposition \ref{prop:alpha}, and  
$\La^{\ma}\otimes_{\max}\La^{\mb}$ is the tensor product 
of the \rrs of $C^*(\ma)$ and $C^*(\mb)$ respectively. Finally, the existence and the surjectivity of  $\tilde{\sigma}^{\max}_{\tilde{\al}}$ is obvious.  
\end{thm}
\begin{proof}
Let $f\in C_c(\ma)$, $g\in C_c(\mb)$. Then, by Proposition
\ref{prop:max} and Lemma \ref{lem:bolu2} we have  
$\La\sigma^{\max}_{\al}j(f\otimes g)
=\La\sigma^{\max}_{\al}(f\oslash g)
=\La_{(f\oslash g)}$. On the other hand Lemma
\ref{lem:bolu2} and Proposition \ref{prop:alpha} imply that 
$j_r\tilde{\sigma}^{\max}_{\al}
(\La^{\ma}\otimes_{\max}\La^{\mb})(f\otimes g)
=j_r\tilde{\sigma}^{\max}_{\al}
               \big(\La^{\ma}_f\otimes\La^{\mb}_g\big)
=j_r\big(\La^{\ma}_f\otimes\La^{\mb}_g\big)
=\La_{(f\oslash g)}$. 
Since $\bats{C_c(\ma)}{C_c(\mb)}$ is dense in
$\tmax{C^*(\ma)}{C^*(\mb)}$, we conclude that 
$\La_{\bts{\ma}{\al}{\mb}}\sigma^{\max}_{\al}j(x)=
j_r\tilde{\sigma}^{\max}_{\al}
(\La^{\ma}\otimes_{\max}\La^{\mb})(x)$, $\forall x\in  
\tmax{C^*(\ma)}{C^*(\mb)}$, and therefore the diagram commutes.  
\end{proof}
\begin{cor}\label{cor:amenamax}
The \fb $\tmax{\ma}{\mb}$ is amenable if and only if $\ma$ and 
$\mb$ are amenable and $\tmax{C^*(\ma )}{C^*(\mb)} =
\bts{C^*(\ma )}{\overline{\max}}{C^*(\mb )}$.  
\end{cor}
\begin{proof}
For $\al =\max$, the diagram $D_{\max}$ becomes:  
\[ \xymatrix
{C^*(\tmax{\ma}{\mb})
\ar@{->>}[rr]^-{\La}
&&C^*_r(\tmax{\ma}{\mb})\\ 
\tmax{C^*(\ma)}{C^*(\mb)}
\ar@{->>}[rr]_-{\tilde{\sigma}^{\max}_{\overline{\max}}\circ 
                (\La^{\ma}\otimes_{\max}\La^{\mb})}
\ar[u]^-{j\cong}
&&\bts{C^*_r(\ma)}{\overline{\max}}{C^*_r(\mb)}
\ar[u]_-{\cong j_r}  }\]
If $\La$ is an isomorphism, then so is $\tilde{\sigma}^{\max}_{\overline{\max}}\circ(\La^{\ma}\otimes\La^{\mb})$, and therefore also $\tilde{\sigma}^{\max}_{\overline{\max}}$, because $\La^{\ma}\otimes\La^{\mb}$ is surjective. Moreover the injectivity $\La$ implies that of $\La^{\ma}\otimes\La^{\mb}$, and therefore that of $\La^{\ma}$ and $\La^{\mb}$,  and also that $\norm{\cdot}_{\max}=\norm{\cdot}_{\overline{\max}}$. In other words, the amenability of $\tmax{\ma}{\mb}$ implies the amenability of $\ma$ and of $\mb$,  and also that $\norm{\cdot}_{\max}=\norm{\cdot}_{\overline{\max}}$. The converse is clear. 
\end{proof}

\subsection{Some applications}\label{subsec:appls}
\par Suppose that $\mb =(B_t)_{t\in G}$ is a \fb over the locally
compact group $G$, and for $\phi,\psi\in C_c\big(G,M(B_e)\big)$ and
$b\in B_t$ let $\phi\cdot
b\cdot\psi:=\int_G\phi_i(s)^*b\psi_i(t^{-1}s)ds$. Since every $x\in
M(B_e)$ defines a multiplier of $\mb$ of order $e$, then we have that  
$\phi\cdot b\cdot \psi\in B_t$, $\forall b\in B_t$. So we have a map
$\Phi_{\phi,\psi}:\mb\to\mb$ defined 
by $b\mapsto\phi\cdot b\cdot\psi$. 
For $b_t\in B_t$ we have 
\begin{gather*}
  \phi\cdot b_t\cdot\psi =\int_{(\supp\phi)\cap (t\supp\psi)}\phi(s)^*b_t\psi(t^{-1}s)ds,\textrm{ so if $m$ is Haar measure:}\\
  \norm{\phi\cdot b_t\cdot\psi}\leq m\big((\supp\phi)\cap (t\supp\psi)\big)\norm{\phi}_\infty\norm{\psi}_\infty\norm{b_t}.
\end{gather*}
Besides, if $f\in C_c(\mb)$, we have $\phi\cdot f\cdot\psi\in C_c(\mb)$, with $\supp(\phi\cdot f\cdot\psi)\subseteq\supp(f)$ and $\norm{\phi\cdot f\cdot\psi}_\infty\leq m\big((\supp\phi)\cap\supp(f)(\supp\psi)\big)\norm{\phi}_\infty \norm{\psi}_\infty\norm{f}_\infty$.
By Lemma
\ref{lem:multis} the map $(s,t)\mapsto\phi(s)^*f(t)\psi(t^{-1}s)$ is 
continuous.
Then \cite[II-15.19]{fd} implies that $\phi\cdot
f\cdot\psi\in C_c(\mb)$.
It follows that the map $\mb\to\mb$ such that $b\mapsto \phi\cdot b\cdot \psi$ is a continuous map on the bundle $\mb$, and that $\Phi_{\phi,\psi}:C_c(\mb)\to C_c(\mb)$ given by $f\mapsto \phi\cdot f\cdot \psi$ is continuous in the inductive limit topology. In fact, in \cite[Lemma~3.2]{exeng} is shown that
\begin{equation}\label{eqn:ap}
  \norm{\Phi_{\phi,\psi}(b)}\leq\norm{\phi}\,\norm{\psi}\norm{b},
\end{equation}
where $\norm{\phi}$ and $\norm{\psi}$ are the norms of $\phi$ and $\psi$ as elements of $L^2(G,M(B_e))$. Hence we also have $\norm{\Phi_{\phi,\psi}(f)}_\infty\leq\norm{\phi}\,\norm{\psi}\norm{f}_\infty$ and $\norm{\Phi_{\phi,\psi}(f)}_1\leq\norm{\phi}\,\norm{\psi}\norm{f}_1$ $\forall f\in C_c(\mb)$, and therefore $\Phi_{\phi,\psi}$ extends to a bounded map on $L^1(\mb)$.        

\begin{df}\label{df:app}(cf. \cite[Definition~3.6]{exeng})
Let \mb be a \fb over the locally compact group $G$, and $M\geq 0$. 
\be
\item We say that $\mb$ has the pointwise $M$-approximation
      property if there exist nets 
      $(\phi_i)_{i\in I},(\psi_i)_{i\in I}\subseteq
      C_c\big(G,M(B_e)\big)$ such that:\\  
      (i) $\sup_{i\in I}\{\norm{\phi_i}\,\norm{\psi_i}\}\leq M$ (as
          elements of $\bts{L^2(G)}{}{M(B_e)}$), and \\ 
      (ii) $\phi_i\cdot b\cdot\psi_i$ converges to $b$, $\forall 
           b\in\mb$.\\ 
      If $I=\nt$ we say that $\mb$ has the countable $M$-pointwise
      approximation property. We say that $\mb$ has the (countable) pointwise
      approximation property if $\mb$ has the (countable) $M$-pointwise
      approximation property for some $M>0$.  
    \item $\mb$ is said to have the $M$-\ap 
      if there are 
      nets $(\phi_i),(\psi_i)$ as in (1) such that $\phi_i\cdot f\cdot\psi_i$
      converges uniformly to $f$, $\forall f\in C_c(\mb)$. 
      It is said to have the \textit{\ap} if it has the $M$-\ap for some $M\geq 0$. 
\item We say that $\mb$ has the $L^1$-\ap if there are 
      nets $(\phi_i),(\psi_i)$ as in (1) such that $\phi_i\cdot f\cdot\psi_i$
      converges to $f$ in $L^1(\mb)$, $\forall f\in C_c(\mb)$.   
      \ee
      In all the cases above we say that $\mb$ has the \textit{positive} corresponding approximation property if we can choose $\phi_i=\psi_i$, $\forall i$.  
\end{df} 
\medskip
\par The fact that we allow the approximating nets $(\phi_i)_{i\in I}$ and $(\psi_i)_{i\in I}$ to take values on the multiplier algebra $M(B_e)$ rather than in $B_e$ is not an essential change in relation to the original definition of \ap, but it allows some more flexibility (it is enough to multiply the approroximating nets by an approximate unit of $B_e$ to obtain nets as in \cite[Definition~3.6]{exeng}).     
\par It was proved in \cite{exeng} that if $G$ is an amenable group then the 
\fb has the positive 1-\ap. 
\par For a Fell bundle $\mb$ over a discrete group it is currently customary to say that $\mb$ has the approximation property when it has the positive 1-approximation property. The corresponding net is called a Cesaro net for $\mb$ by Exel in \cite[Definition~20.4]{pds}.
\medskip
\par Since $L^2(G)$ is a Hilbert space, it is a nuclear \ct, so there
is a unique tensor product $\bts{L^2(G)}{}{M(B_e)}$. On the other hand 
$\bts{L^2(G)}{}{M(B_e)}$ is naturally identified with 
$L^2\big(G,M(B_e)\big)$, the completion of $C_c\big(G,M(B_e)\big)$
with respect to the inner product: $\pr{f}{g}=\int_Gf(t)^*g(t)dt$.
Thus we also have that $L^2\big(G,M(B_e)\big)=L^2\big(G\times M(B_e)\big)$, where 
$G\times M(B_e)$ is the \fb over $G$ with the product topology and 
pointwise defined operations.
  
\begin{prop}\label{prop:casos}
Let $\mb$ be a \fb over the locally compact group $G$. 
We have:  
\be
 \item If $G$ is discrete, the three next statements are equivalent to each other: $\mb$ has the $M$-pointwise \ap; $\mb$ has the $M$-\ap; $\mb$ has the $M$- $L^1$-\ap.  
 \item If \mb has the \ap then it also has the $L^1$-appro\-ximation property.         
 \item If $\mb$ has the countable pointwise \ap, then $\mb$ has the $L^1$-\ap. 
\ee
\end{prop}
\begin{proof}
The first statement easily follows by observing that, if $G$ is discrete and $f\in C_c(\mb)$, then $\textrm{supp}(f)$ is finite.
Suppose now that $(\phi_i)$, $(\psi_i)\subseteq C_c(G,M(B_e))$ are nets
such that $\phi_i\cdot f\cdot\psi_i$ converges uniformly to $f$,
$\forall f\in C_c(\mb)$, with $\sup_i\norm{\phi_i}\,\norm{\psi_i}\leq
M<\infty$. Thus $\phi_i\cdot f\cdot\psi_i$ converges 
to $f$ in the \ilt, because $\textrm{supp}(\phi\cdot
f\cdot\psi)\subseteq\textrm{supp}(f)$. Therefore the net $\phi_i\cdot
f\cdot\psi_i$ converges to $f$ in $L^1(\mb)$. 
\par Suppose now that $\mb$ has the countable pointwise \ap: there
exist sequences $(\phi_n),(\psi_n)\subseteq C_c(G,M(B_e))$ with 
$\sup_{n\in \nt}\{\norm{\phi_n}\,\norm{\psi_n}\}=M<\infty$ and $\phi_n\cdot
b\cdot\psi_n\to b$, $\forall b\in \mb$. Let
$\Phi_n:=\Phi_{\phi_n,\psi_n}:C_c(\mb)\to C_c(\mb)$ be the
corresponding induced map. Then, since $\norm{\Phi_n(f)-f}_\infty\leq (M+1)\norm{f}_\infty$, we have that
$\norm{(\Phi_n(f)-f)}_1=\int_{\textrm{supp}(f)}\norm{\Phi_n(f)(t)-f(t)}dt\to 
0$ by the dominated convergence theorem.  
\end{proof}

\par The following theorem is a direct generalization of the
corresponding result \cite[Theorem~4.6]{examen} for discrete groups, so we omit
the proof here, although for the convenience of the reader we have provided its details in Appendix~\ref{pf:pal1}.    

\begin{thm}\label{thm:pal1}
If $\mb$ is a \fb with the $L^1$-\ap, then $\mb$ is amenable. In particular if $\mb$ has the \ap, then $\mb$ is amenable. 
\end{thm}

\par Let \fela and \felbh be \fbs, and $\al$ a $C^*$-norm on
$\bats{\ma}{\mb}$. Let $\chal$ be the $C^*$-norm on
$\bats{M(A_e)}{M(B_e)}$ as a subalgebra of $M(\bts{A_e}{\alpha}{B_e})$ (see \cite[T.6.3]{w}, or alternatively use Corollary~\ref{cor:alfita} and \cite[Theorem 2.4]{l} for the right Hilbert modules $A_e$ and $B_e$ over themselves), and $\bts{M(A_e)}{\overline{\al}}{M(B_e)}$ the corresponding tensor
product.
If $\phi\in C_c(G,M(A_e))$ and $\phi'\in
C_c(H,M(B_e))$, we have a section $\phi\otimes\phi'$ $\in 
C_c\big(G\times H,\bts{M(A_e)}{\chal}{M(B_e)}\big)\subseteq 
L^2(G\times H,M(\bts{A_e}{\al}{B_e}))$ such that
$\phi\otimes\phi'(t,s)=\phi(t)\otimes\phi'(s)$, $\forall (t,s)\in
G\times H$. Moreover: 
$\norm{\phi\otimes\phi'}=\norm{\phi}\,\norm{\phi'}$ by Remark~\ref{rk:crossnorms}.   

\begin{prop}\label{prop:tensap}
Let \fela and \felbh be \fbs, and $\al$ a $C^*$-norm on
$\bats{\ma}{\mb}$. Suppose that $(\phi_i)_{i\in I}, (\psi_i)_{i\in
I}\subseteq C_c(G,M(A_e))$ and $(\phi_j')_{j\in J}, (\psi_j')_{j\in
J}\subseteq C_c(H,M(B_e))$. Consider
$(\phi_i\otimes\phi_j')_{(i,j)\in I\times J}$,
$(\psi_i\otimes\psi_j')_{(i,j)\in I\times J}\subseteq C_c\big(G\times
H,M(\bts{A_e}{\al}{B_e})\big)$. Then:   
\be
\item If $\phi_i\cdot a\cdot\psi_i\to a$, $\forall a\in\ma$ and
      $\phi_j'\cdot b\cdot\psi_j'\to b$ and $\sup_{i\in
      I}\{\norm{\phi_i}\,\norm{\psi_i}\}\leq M<\infty$,  
      $\sup_{j\in J}\{\norm{\phi_j'}\,\norm{\psi_j'}\}\leq
      N<\infty$, then $(\phi_i\otimes\phi_j')\cdot x\cdot
      (\phi_i\otimes\phi_j')\to x$, $\forall x\in\bts{\ma}{\al}{\mb}$,
      and $\sup_{(i,j)\in I\times J}\{\norm{\phi_i\otimes 
      \phi_j'}\,\norm{\psi_i\otimes\psi_j'}\}\leq MN<\infty$.  
\item If $\ma$ and \mb have the (positive, countable) pointwise \ap,
      then $\bts{\ma}{\al}{\mb}$ also has the (respectively: positive, 
      countable) pointwise \ap.
\item If $\ma$ and $\mb$ have the $L^1$-\ap then
      $\bts{\ma}{\al}{\mb}$ also has the $L^1$-\ap, and therefore it
      is amenable.
\ee 
\end{prop}
\begin{proof}
\par Note first that 2) follows from 1). To prove 1), let
$\Phi_i:\ma\to\ma$ and $\Phi_j':\mb\to\mb$ be the maps 
induced by the pairs $(\phi_i,\psi_i)$ and $(\phi_j',\psi_j')$,
$\forall (i,j)\in I\times J$.  
Let $a_{t_0}\in\ma$ and $b_{s_0}\in\mb$. Since $\Phi_i$ and $\Phi_j'$ 
converge pointwise to the identity maps respectively on $\ma$ and
$\mb$, there exist $i_0\in I$, $j_0\in J$ such that $\forall i\geq
i_0$, $j\geq j_0$ we have $\norm{\Phi_i(a_{t_0})-a_{t_0}}<\epsilon
/N(1+\norm{a_{t_0}}+\norm{b_{s_0}})$ and $\norm{\Phi_j'(b_{s_0})-b_{s_0}}<\epsilon /(1+\norm{a_{t_0}}+\norm{b_{s_0}})$. 
Consider $\Phi_{i,j}:\bts{\ma}{\al}{\mb}\to\bts{\ma}{\al}{\mb}$ such
that 
\[\Phi_{i,j}(x_{(t,s)})=\int_{G\times H}
(\phi_i\otimes\phi_j')(t',s')x_{(t,s)}
(\psi_i\otimes\psi_j')(t^{-1}t',s^{-1}s')d(t',s').\]    
Then for $(i,j)\geq (i_0,j_0)$ we have:  
\begin{gather*}
  \norm{\Phi_{i,j}(a_{t_0}\otimes b_{s_0})-(a_{t_0}\otimes b_{s_0})}
=\norm{\Phi_i(a_{t_0})\otimes\Phi_j'(b_{s_0})-(a_{t_0}\otimes
b_{s_0})}\\
\leq\norm{\big(\Phi_i(a_{t_0})-a_{t_0}\big)\otimes\Phi_j'(b_{s_0})}+
     \norm{a_{t_0}\otimes\big(\Phi_j'(b_{s_0})-b_{s_0}\big)}\\
\leq \norm{\Phi_i(a_{t_0})-a_{t_0}}\,\norm{\Phi_j'(b_{s_0})}+
      \norm{a_{t_0}}\, \norm{\Phi_j'(b_{s_0})-b_{s_0}}\\
      \leq  \frac{\epsilon}{N(1+\norm{a_{t_0}}+\norm{b_{s_0}})} N\norm{b_{s_0}}+\frac{\epsilon}{(1+\norm{a_{t_0}}+\norm{b_{s_0}})} \norm{a_{t_0}}
      <\epsilon.
\end{gather*}
By \eqref{eqn:ap} 
we have $\norm{\Phi_{i,j}(x)}\leq MN\norm{x}$, $\forall 
x\in\bts{\ma}{\al}{\mb}$, and consequently $\Phi_{i,j}(x)\to x$, 
$\forall x\in\bts{\ma}{\al}{\mb}$. This proves 1) and therefore also
2). 
\par To see that 3) holds, suppose now that for the maps $\Phi_i$ and
$\Phi_j'$ above and every $f\in C_c(\ma)$, $g\in C_c(\mb)$ we have
that $\norm{\Phi_i(f)-f}_1 \to 0$ and $\norm{\Phi_j'(g)-g}_1\to 0$.   
Note that if $f\in C_c(\ma)$ and $g\in C_c(\mb)$, then
$\Phi_{i,j}(f\oslash g)=\Phi_i(f)\oslash\Phi_j'(g)$, and therefore   
\begin{align*}
\norm{\Phi_{i,j}(f\oslash g)-f\oslash g}_1
&=\norm{\Phi_i(f)\oslash\Phi_j'(g)-f\oslash g}_1\\
&\leq\norm{\Phi_i(f)\oslash\big(\Phi_j'(g)-g\big)}_1
 +\norm{\big(\Phi_i(f)-f\big)\oslash g}_1\\
&\leq M\,\norm{f}_1\,\norm{\Phi_j'(g)-g}_1
 +\norm{\Phi_i(f)-f}_1\, \norm{g}_1\\
&\to 0\resp\text{when $i,j\to\infty$}
\end{align*}
It follows that $\Phi_{i,j}(l)\to l$ in $L^1(\bts{\ma}{\al}{\mb})$,
 $\forall l\in L=\{\sum_kf_k\oslash g_k\}$. Since $L$ is dense in  
$C_c(\bts{\ma}{\al}{\mb})$ in the \ilt, it is also dense in 
$L^1(\bts{\ma}{\al}{\mb})$. Since $\norm{\Phi_{i,j}}\leq MN$, $\forall
i\in I,j\in J$, then $\norm{\Phi_{i,j}(h)-h}_1\to 0$,
 $\forall h\in L^1(\bts{\ma}{\al}{\mb})$. 
\end{proof}
The last statement of the next result was first proved by the author
in the case of discrete groups in the previous preprint version of the
present paper mentioned at the end of the introduction, and was later proved for arbitrary locally compact groups in \cite{favthesis} and in \cite{exeng}.  
\begin{cor}\label{cor:amenuclear}
If  $\fela$ and $\felbh$ are \fbs with the $L^1$-\ap, then 
$\bats{A_e}{B_e}$ admits exactly one $C^*$-norm if and only if 
$\bats{C^*(\ma)}{C^*(\mb)}$ admits exactly one $C^*$-norm. In   
particular, if $\ma$ is a \fb with the $L^1$-\ap (this is
automatically true if $G$ is amenable) and nuclear unit
fiber $A_e$, then $C^*(\ma)$ is also nuclear. 
\end{cor}
\begin{proof}
Since $\ma$ and $\mb$ are \fbs with the $L^1$-\ap, then 
$\tmin{\ma}{\mb}$ also has the $L^1$-\ap by Proposition
\ref{prop:tensap}, so the diagram $D_{\min}$ becomes:  
\[ \xymatrix
{C^*(\tmax{\ma}{\mb})
\ar@{->>}[r]^-{\sigma^{\max}_{\min}}
&C^*(\tmin{\ma}{\mb})\\ 
\tmax{C^*(\ma)}{C^*(\mb)}
\ar@{->>}[r]_-{\tilde{\sigma}^{\max}_{\min}}\ar[u]^-{j\cong}
&\tmin{C^*(\ma)}{C^*(\mb)}
\ar[u]_-{\cong j_r}  }\]
If \bats{A_e}{B_e} admits just one $C^*$-norm, then 
$\tmax{\ma}{\mb}=\tmin{\ma}{\mb}$ which implies 
$\sigma^{\max}_{\min}=id$, and therefore  
$\tilde{\sigma}^{\max}_{\min}=id$, from where it follows that 
$\tmax{C^*(\ma)}{C^*(\mb)}=\tmin{C^*(\ma)}{C^*(\mb)}$. 
Conversely, suppose now that \bats{C^*(\ma)}{C^*(\mb)} admits  
just one $C^*$-norm. Then 
$\tilde{\sigma}^{\max}_{\min}=id$, and therefore 
$\sigma^{\max}_{\min}=id$. Thus  
$\tmax{\ma}{\mb}=\tmin{\ma}{\mb}$, so 
$\tmax{A_e}{B_e}=\tmin{A_e}{B_e}$. 
\par As for the last assertion, notice that every \cs $B$ may be
considered as a \fb over the trivial group, and it is clear that this
\fb has the (positive, countable) $L^1$-\ap: it is enough to take
$\phi:G\to M(B)$ such that $\phi(t)=1$, $\forall t\in G$.    
Consequently, by the first part of this Corollary we have 
$\tmax{C^*(\ma)}{B}=\tmin{C^*(\ma)}{B}$, that is,  
$C^*(\ma)$ is a nuclear \cs.
\end{proof}

\begin{cor}
If $\fela$ is a \fb with the $L^1$-\ap and nuclear unit fiber, and if 
$\felbh$ is an amenable \fb, then $\bts{\ma}{}{\mb}$ also is
amenable. 
\end{cor}
\begin{proof}
By Corollary \ref{cor:amenuclear}, our assumptions on 
$\ma$ imply that $C^*(\ma)$ is nuclear. Therefore we have 
$\tmax{C^*(\ma)}{C^*(\mb)}=\bts{C^*(\ma)}{\overline{\textrm{max}}}{C^*(\mb)}$,  
and then the result follows from Corollary \ref{cor:amenamax}.
\end{proof}

\begin{cor}\label{cor:tpcp}
Any twisted partial crossed product of a nuclear \cs by an amenable
group is nuclear. In particular, the partial \cs $C^*_p(G)$ of an
amenable discrete group $G$ is nuclear.
\end{cor}

\begin{df}\label{df:seqfell}
Let $\ma$, $\mb$ and $\mc$ be \fbs over the locally compact group $G$.   
We say that a sequence 
$\xymatrix@1
{0\ar[r]&\ma\ar[r]^-{\phi}
           &\mb\ar[r]^-{\psi}
           &\mc\ar[r]&0} $ 
is exact if $\phi$ is injective, $\psi$ is surjective, and $\ker\psi =
\text{Im}\phi$, where $\ker\psi:=\{b\in \mb:\psi(b)\textrm{ is a zero
element}\}$.  
\end{df}

\begin{prop}\label{prop:seqfell}
The functors $\ma\mapsto L^1(\ma)$ and $\ma\mapsto C^*(\ma)$ are
exact. That is, if 
 $ \xymatrix@1
   {0\ar[r]&\ma\ar[r]^-{\phi}
           &\mb\ar[r]^-{\psi}
           &\mc\ar[r]&0} $ is an \es of \fbs over the locally compact
group $G$, then:
\be
 \item $\xymatrix@1
        {0\ar[r]&L^1(\ma)\ar[r]^-{\phi^1}
           &L^1(\mb)\ar[r]^-{\psi^1}
           &L^1(\mc)\ar[r]&0}$ is exact, and 
 \item  $\xymatrix@1
        {0\ar[r]&C^*(\ma)\ar[r]^-{C^*(\phi)}
           &C^*(\mb)\ar[r]^-{C^*(\psi)}
           &C^*(\mc)\ar[r]&0}$ also is exact.  
\ee 
\end{prop}
\begin{proof}
Since every non-degenerate \rep of $L^1(\ma)$ has a unique extension
to a \rep of $L^1(\mb)$, we have that $C^*(\ma)$ is the closure of
$L^1(\ma)$ in $C^*(\mb)$, so it is enough to prove 1), because then 2)
follows from \cite[2.29]{zm} and the fact that $L^1(\mathcal{F})$ has
an approximate unit, for every \fb $\mathcal{F}$.
\par Since $\norm{\phi(a)}=\norm{a}$, $\forall a\in\ma$, it follows
that $\phi^1$ is an isometry. 
\par Let see that $\ker\psi^1=\text{Im}\phi^1$. The inclusion
$\text{Im}\phi^1\subseteq\ker\psi^1$ is clear. In order to see the
converse inclusion let $g\in\ker\psi^1$. Then $\norm{\psi^1(g)}_1=0$,
so $\psi\big(g(t)\big)=0$ almost everywhere in $G$. Without loss of
generality we may suppose that $\psi\big(g(t)\big)=0$, $\forall t$ in
$G$. Thus, $g(t)\in\ker\psi=\text{Im}\phi$, $\forall t\in G$, and
therefore there exists a unique   
$f(t)\in A_t$ such that $\phi\big(f(t)\big)=g(t)$, $\forall t\in
G$. Since $\phi$ is a continuous and isometric isomorphism between
$\ma$ and $\phi(\ma)$ (by \cite[II-13.17]{fd}), and we have
$f=(\phi^1)^{-1}(g)$, then $f\in L^1(\ma)$. Thus
$g\in\text{Im}\phi^1$. 
\par Finally, we show that $\psi^1$ is surjective. We will suppose, as
we can, that $\ma\subseteq\mb$. Note that $L^1(\ma)$ is a closed
*-ideal of $L^1(\mb)$. Thus there exists an isomorphism of   
*-algebras $\frac{L^1(\mb
)}{L^1(\ma)}\stackrel{\bar{\psi}}{\to}{\psi^1\big(L^1(\mb)\big)}\subseteq
L^1(\mc)$. The image of $\bar{\psi}$ contains
$\psi^1\big(C_c(\mb)\big)$, which is dense in $C_c(\mc)$ in the \ilt:
since $\psi$ is surjective then we may apply \cite[5.1]{env} to   
$\psi^1\big(C_c(\mb)\big)$, thus concluding that $\text{Im}\bar{\psi}$
is dense in $L^1(\mc)$. Then it is sufficient to prove that
$\bar{\psi}$ is an isometry, where $\frac{L^1(\mb )}{L^1(\ma)}$ is
endowed with the quotient norm. Let $f\in C_c(\mb)$ and $\bar{f}$ its
projection into the quotient space. Then
$\norm{\bar{\psi}}=\norm{\psi^1}$, and therefore
$\norm{\bar{\psi}(\bar{f})}\leq\norm{\psi^1}\,\norm{\bar{f}} 
\leq\norm{\bar{f}}$. To prove the converse inequality consider an
arbitrary $\epsilon >0$, and let $M$ be the measure of a compact 
neighborhood $V$ of $\textrm{supp}(f)$. For each 
$s\in V$, there exists $g_s\in C_c(\ma )$ such that 
$\norm{f(s)-g_s(s)}<\norm{\ov{f(s)}}+\epsilon /M$, and we may suppose
that $\textrm{supp}(g_s)\subseteq V$. Since $f$, $g_s$, and
$t\mapsto\norm{\ov{f(t)}}=\norm{\psi(f(t))}$ are continuous, for every
$s\in\textrm{supp}(f)$ 
must exist an open neighborhood $V_s$ of $s$, which we may suppose to 
be contained in $V$, such that
$\norm{f(t)-g_s(t)}<\norm{\ov{f(t)}}+\epsilon /M$, $\forall t\in V_s$.   
Now, $\{ V_s:\ s\in \sup(f)\}$ is an open covering of the compact set 
$\textrm{supp}(f)$. Let $V_{s_1}, \ldots,V_{s_n}$ be a finite
subcovering. Let $G_{\star}$ be the one point compactification of $G$,
and define $s_{n+1}:=\star$,
$V_{s_{n+1}}=G_{\star}\setminus\textrm{supp}(f)$ and   
$g_{s_{n+1}}=0$, where $\star$ represents the adjoined point at
infinity. Then $\{V_{s_i}\}_{i=1}^{n+1}$ is an open covering of
$G_{\star}$. Let $(\phi_i)_{i=1}^{n+1}$ be a partition of the unit of
$G_{\star}$, subordinated to $\{V_{s_i}\}_{i=1}^{n+1}$, and define   
$g(t)=\sum_{i=1}^{n+1} \phi_i(t)g_{s_i}(t)$, $\forall t\in G$. Then
$g\in C_c(\ma)$, $\textrm{supp}(g)\subseteq V$, and 
  \begin{gather*}
   \norm{\bar{f}}
   \leq\int_G\norm{f(t)-g(t)}dt
   =\int_{V}\norm{\sum_{i=1}^{n+1}
        \big(\phi_i(t)f(t)-\phi_i(t)g_{s_i}(t)\big)}dt\\ 
   \leq\sum_{i=1}^{n+1}\int_{V_{s_i}}\phi_i(t)\norm{f(t)-g_{s_i}(t)}dt
   \leq\int_{V}\sum_{i=1}^{n+1}\phi_i(t)\big(\norm{\ov{f(t)}}+
                                                     \epsilon /M\big)dt\\ 
   \leq\int_G\norm{\psi\big(f(t)\big)}dt + \epsilon
   =\norm{\psi^1(f)}_1 +\epsilon . 
  \end{gather*}
Since $\epsilon$ was arbitrary, we conclude that 
$\norm{\bar{f}}\leq\norm{\psi^1(f)}_1$, and therefore
$\norm{\bar{f}}=\norm{\bar{\psi}(f)}_1$. Moreover $\bar{\psi}$ has
dense image in $L^1(\mc)$, so $\text{Im}(\bar{\psi})=L^1(\mc)$ and, 
since $\text{Im}(\psi^1)=\text{Im}(\bar{\psi})$, we conclude that
$\psi^1$ is surjective.
\end{proof}
\par For our next result, recall from \cite[5.3]{trings} that the definition of exact C*-algebra extends to C*-trings, and that a C*-tring $E$ is exact if and only if $E^r$ is exact. In particular, if $\ma=(A_t)_{t\in G}$ is a Fell bundle, its unit fiber is an exact C*-algebra if and only if each fiber $A_t$ is an exact C*-tring.  
\begin{thm}\label{thm:mainexact}
Let $\ma =(A_t)_{t\in G}$ be a \fb with exact unit fiber and the  
$L^1$-\ap. Then $C^*(\ma )$ is an exact \cs.
\end{thm}
\begin{proof}
Let $B$ be a \cs and $I\id B$. Since 
$A_t$ is exact $\forall t\in G$, the sequence of \fbs 
$ \xymatrix@1
   {0\ar[r]&\bts{\ma}{}{I}\ar[r]
           &\bts{\ma}{}{B}\ar[r]
           &\bts{\ma}{}{(B/I)}\ar[r]&0}$
is exact, and each one of the bundles in the sequence has the
   $L^1$-\ap by Proposition \ref{prop:tensap} (here
$\bigotimes=\bigotimes_{\min}$). Since $C^*$ is an exact functor from   
the category of Banach *-algebras with approximate unit to the
category of \css, and in this case we have $C^*=C^*_r$, the sequence
of \css  
\[ \xymatrix
   {0\ar[r]&C^*_r\big(\bts{\ma}{}{I}\big)\ar[r]
           &C^*_r\big(\bts{\ma}{}{B}\big)\ar[r]
           &C^*_r\big(\bts{\ma}{}{(B/I)}\big)\ar[r]&0} \] 
also is exact. 
\par Now Proposition \ref{prop:alpha} provides a natural
isomorphism between $C^*_r\big(\bts{\ma}{}{C}\big)$ and
$\bts{C^*_r(\ma )}{}{C}$, for every \cs $C$. Thus we obtain the
following commutative diagram: 
\[ \xymatrix
   {0\ar[r]&C^*_r(\bts{\ma}{}{I})\ar[r]\ar[d]_{\cong}
           &C^*_r(\bts{\ma}{}{B})\ar[r]\ar[d]_{\cong}
           &C^*_r(\bts{\ma}{}{(B/I)})\ar[r]\ar[d]_{\cong}&0\\
    0\ar[r]&\bts{C^*_r(\ma )}{}{I}\ar[r]
           &\bts{C^*_r(\ma )}{}{B}\ar[r]
           &\bts{C^*_r(\ma )}{}{(B/I)}\ar[r]&0       
                       }\] 
Since the first row is exact and the diagram is commutative, then the
second row also is exact. Hence it follows that $C^*(\ma )$, which is
equal to $C^*_r(\ma )$, is an exact \cs. 
\end{proof}

\par \par Since any Fell bundle over an amenable locally compact group has the approximation property, from Theorem~\ref{thm:mainexact} we obtain the following generalization of \cite[Proposition~7.1]{k} (see also \cite[Proposition~7.5]{bew}): 

\begin{cor}\label{cor:crossprods}
Any twisted partial crossed product of an exact \cs by an amenable
group $G$ is also exact.
\end{cor}



\section{Appendix}
\begin{proof}[Proof of Theorem~\ref{thm:pal1}]\label{pf:pal1}
Let $\pi:\mb\to B(H)$ be a non-degenerate \rep such that  
$\pi\r{B_e}$ is faithful. We also call $\pi$ the integrated \rep of $\pi$, 
and to its unique extension to $M(\mb)$ as well. 
Since $\pi\r{B_e}$ is faithful, then so is $\pi\r{M(B_e)}$. 
Given $\phi\in C_c\big(G,M(B_e)\big)\subseteq \bts{L^2(G)}{}{M(B_e)}$, 
consider the operator $V_{\phi}:H\to\bts{L^2(G)}{}{H}$ such that  
$V_{\phi}h\r{t}=\pi\big(\phi(t)\big)h$. 
We have:
\begin{gather*}
\norm{V_{\phi}}^2
=\sup_{\norm{h}=1}\int_G\pr{\pi\big(\phi(t)\big)h}{\pi\big(\phi(t)\big)h}dt
=\sup_{\norm{h}=1}\int_G\pr{\pi\big(\phi(t)^*\phi(t)\big)h}{h}dt\\
=\sup_{\norm{h}=1}\pr{\int_G\pi\big(\phi(t)^*\phi(t)\big)dt\, h}{h}
=\sup_{\norm{h}=1}\pr{\pi(\pr{\phi}{\phi})h}{h}\\
=\norm{\pi(\pr{\phi}{\phi})}
=\norm{\pr{\phi}{\phi}}.
\end{gather*}
where we used that $\pi(\pr{\phi}{\phi})$ is positive and $\pi\r{M(B_e)}$ is faithful to obtain the last two equalities.
We compute $V_{\phi}^*$: if $h\in H$, $\xi\in L^2(G)$, then 
\[\pr{V_{\phi}h}{\xi}
=\int_G\pr{\pi\big(\phi(t)\big)h}{\xi(t)}dt
=\int_G\pr{h}{\pi\big(\phi(t)^*\big)\xi(t)}dt
=\pr{h}{\int_G\pi\big(\phi(t)^*\big)\xi(t)dt},\] 
so 
$V_{\phi}^*(\xi)=\int_G\pi\big(\phi(t)^*\big)\xi(t)dt$. Note that if  
$\psi\in C_c\big(G,M(B_e)\big)$, then
\[V_{\phi}^*V_{\psi}h=
\int_G\pi\big(\phi(t)^*\big)V_{\psi}h\r{t}dt
=\int_G\pi\big(\phi(t)^*\big)\pi\big(\psi(t)\big)hdt
=\pi(\pr{\phi}{\psi})h.\] 
Moreover, if $\phi_1$, $\phi_2$, $\phi_3\in C_c\big(G,M(B_e)\big)$, 
$h\in H$, we have:
\begin{gather*}
V_{\phi_1\pr{\phi_2}{\phi_3}}h\r{t}
=\pi(\phi_1(t)\pr{\phi_2}{\phi_3})h
=\pi(\phi_1(t))\pi(\pr{\phi_2}{\phi_3})h
=V_{\phi_1}V_{\phi_2}^*V_{\phi_3}h\r{t}
\end{gather*}
Therefore, since $\phi\mapsto V_{\phi}$ is an isometry on the dense 
subspace $C_c\big(G,M(B_e)\big)$ of the C*-tring $L^2(G,M(B_e))$, it extends to a \hm of positive C*-trings  
$\pi_2:\bts{L^2(G)}{}{M(B_e)}\to\pi_2\big(\bts{L^2(G)}{}{M(B_e)}\big)\subseteq 
B\big(H,\bts{L^2(G)}{}{H}\big)$, which is consequently an isomorphism of C*-trings. 
\par Consider now the \rep $\pi_{\la}:\mb\to B\big(\bts{L^2(G)}{}{H}\big)$, 
such that  $\pi_{\la}(b_t)=\la_t\otimes\pi(b_t)$ and its integrated \rep, 
which we continue to call $\pi_{\la}:C^*(\mb)\to B(\bts{L^2(G)}{}{H})$ (here $\lambda$ is the left regular representation of $G$). 
Define, for $\phi$, $\psi\in\bts{L^2(G)}{}{M(B_e)}$, the completely bounded map $\Psi:\pi_{\la}\big(C^*(\mb)\big)\to B(H)$, given by  
$\Psi(x)=V_{\phi}^*xV_{\psi}$, $\forall x\in \pi_{\la}\big(C^*(\mb)\big)$. 
We have $\norm{\Psi(x)}\leq\norm{\phi}\,\norm{\psi}\,\norm{x}$, so  
$\norm{\Psi}\leq\norm{\phi}\,\norm{\psi}$. 
Consider also, for $f\in C_c(\mb)$, the function $\Phi (f):G\to\mb$ 
such that $\Phi(f)\r{t}=\int_G\phi(s)^*f(t)\psi(t^{-1}s)ds$. Let $F(t,s)=
\phi(s)^*f(t)\psi(t^{-1}s)$. By \ref{lem:multis}, 
$F:G\times G\to\mb$ is continuous and of compact support, and such that 
$F(t,s)\in B_t$, $\forall t\in G$. Then by \cite[II-15.19]{fd}, 
the function $t\mapsto\int_GF(t,s)ds$ is a compactly supported continuous section of $\mb$. In other words, $\Phi (f)\in C_c(\mb)$. In fact, it is clear that  
$\supp\big(\Phi(f)\big)\subseteq\supp(f)$. Besides, we have   
$\pi\big(\Phi(f)\big)=\Psi\big(\pi_{\la}(f)\big)$, for if $h,k\in H$: 
\begin{gather*}
\pr{\pi\big(\Phi(f)\big)h}{k}
=\pr{\int_G\pi\big(\Phi(f)\r{t}\big)h}{k}
=\int_G\pr{\pi\big[\int_G\phi(s)^*f(t)\psi(t^{-1}s)ds\big]h}{k}\\
=\int_G\int_G\pr{\pi\big(\phi(s)^*f(t)\psi(t^{-1}s)\big) h}{k}dsdt\\
=\int_G\int_G\pr{\pi\big(\phi(s)^*\big)\pi\big(f(t)\big)
                  \pi\big(\psi(t^{-1}s)\big)h}{k}dsdt\\
=\int_G\pr{\pi\big(\phi(s)^*\big)\int_G\pi\big(f(t)\big)
            V_{\psi}h\r{t^{-1}s}dt}{k}ds\\
=\int_G\pr{\pi\big(\phi(s)^*\big)
     \big[\int_G(\la_t\otimes\pi)\big(f(t)\big)(V_{\psi}h)dt\big]\r{s}}{k}ds\\
=\int_G\pr{\pi\big(\phi(s)^*\big)\pi_{\la}(f)(V_{\psi}h)(s)ds}{k}
=\pr{\int_G\pi\big(\phi(s)^*\big)\big(\pi_{\la}(f)V_{\psi}h(s)\big)ds}{k}\\
=\pr{(V_{\phi}^*\pi_{\la}(f)V_{\psi})h}{k}
=\pr{\Psi\big(\pi_{\la}(f)\big)h}{k},
\end{gather*}
whence $\pi\big(\Phi(f)\big)=\Psi\big(\pi_{\la}(f)\big)$. 
In particular we have $\Psi\big(\pi_{\la}\big(C_c(\mb)\big)\big)
\subseteq\pi\big(C_c(\mb)\big)\subseteq\pi\big(C^*(\mb)\big)$, 
which is closed, and therefore  
$\Psi\big(\pi_{\la}\big(C^*(\mb)\big)\big)\subseteq\pi\big(C^*(\mb)\big)$. 
Then we have $\Psi:\pi_{\la}\big(C^*(\mb)\big)\to 
\pi\big(C^*(\mb)\big)$. 
Suppose now that $(\phi_i)$, $(\psi_i)$ are approximating nets as in (3) of Definition~\ref{df:app}, with  
$\norm{\phi_i}\,\norm{\psi_i}\leq M$, $\forall i$, so we have $\Phi_i:C_c(\mb)\to C_c(\mb)$ and $\Phi_i(f)$ converges to $f$ in $L^1(\mb)$, for all $f\in C_c(\mb)$. Let $\Psi_i:\pi_{\la}\big(C^*(\mb)\big)\to\pi\big(C^*(\mb)\big)$ be the 
correspondending completely bounded maps, that is:   
$\Psi_i(x)=V_{\phi_i}^*xV_{\psi_i}$, 
$\forall x\in \pi_{\la}\big(C^*(\mb)\big)$. 
Since $\Phi_i(f)\to f$ in $L^1(\mb)$, then  
$\Phi_i(f)\to f$ also in $C^*(\mb)$, thus $\pi\big(\Phi_i(f)\big)\to\pi(f)$ 
in $\pi\big(C^*(\mb)\big)$. Consequently,  
$\norm{\pi\big(\Phi_i(f)\big)}\to\norm{\pi(f)}$. On the other hand, 
$\pi\big(\Phi_i(f)\big)=\Psi_i\big(\pi_{\la}(f)\big)$, whence 
\begin{gather*}
  \norm{\pi(f)}
=\lim_i\norm{\pi\big(\Phi_i(f)\big)}
=\lim_i\norm{\Psi_i\big(\pi_{\la}(f)\big)}\\ 
\leq\limsup_i\norm{\Psi_i}\,\norm{\pi_{\la}(f)}
\leq M\norm{\pi_{\la}(f)}.
\end{gather*}
\par Since $C_c(\mb)$ is dense in $L^1(\mb)$, it follows that  
$\norm{\pi(y)}\leq M\norm{\pi_{\la}(y)}$, $\forall y\in C^*(\mb)$. In particular, if $\pi$ is a faithful \rep of $C^*(\mb)$, we conclude that  $\pi_{\la}$ is also faithful. On the other hand, it is proved in \cite[Proposition~2.3]{exnote} that the \rep $\Lambda\otimes id:\mb\to B(L^2(\mb)\otimes_{B_e} H)$, given by $(\Lambda\otimes id)_b(\xi\otimes h):=\Lambda_b\xi\otimes h$, is equivalent to a subrepresentation of $\pi_\lambda$, so it is faithful as well. This implies that $\Lambda$ is faithful, which is to say that $C^*(\mb)=C^*_r(\mb)$.
\end{proof}

\bigskip\bigskip\bigskip


\end{document}